\documentclass[11pt]{article}
\usepackage{amssymb}
\usepackage{amsmath}
\usepackage{bbm}
\usepackage{multirow}	  
\usepackage{enumerate}
\usepackage[dvips]{graphics}
\usepackage{color}
\usepackage{arydshln}	 
\graphicspath{{images/}} 
\usepackage{ifpdf}
\ifpdf
\usepackage{epstopdf}
\fi
\def\vv#1{\text{\rm\bfseries #1}}
\def\vvm#1{\text{\boldmath $#1$}}
\def\inn#1{\langle{#1}\rangle}
\def\qin{\quad{\rm in}\quad}
\def\qon{\quad{\rm on}\quad}
\def\qan{{\qquad\hbox{and}\qquad}}
\def\mat#1{\left(#1\right)}
\def\td{\text{\tt d}}
\def\tra{\text{\tt t}}
\def\tr{\text{\rm tr}}
\def\div{\text{\rm div}}
\def\bdiv{\textrm{\bf div}}
\def\rot{\text{\rm rot}}
\def\brot{\textrm{\bf rot}}
\def\kron{\textrm{kron}}
\def\R{\text{\rm R}}
\def\Lv{\vv{L}^2}
\def\Lm{\mathbb{L}^2}
\def\Hv{\vv{H}}
\def\Hm{\mathbb{H}}
\def\P{\textrm{P}}
\def\Pv{\vv{P}}
\def\Pm{\mathbb{P}}
\def\Th{\mathcal{T}_h}
\def\Gk{\mathcal{G}_k^{\perp}(K)}
\def\Gkb{\vvm{\mathcal{G}}_k^{\perp}(K)}
\def\psiH{\Psi^K}
\def\PsiH{\vec{\vvm{\Psi}}^K}
\def\psiV{\psi^K}
\def\PsiV{\vec{\vvm{\psi}}^K}
\def\Rk{\mathcal{R}_k^K}
\def\Pk{\mathcal{P}_k^K}
\def\bPk{\vvm{\mathcal{P}}_k^K}
\def\nv{\vvm{n}}
\def\Id{\vv{I}}
\def\bsig{\vvm{\sigma}}
\def\btau{\vvm{\tau}}
\def\brho{\vvm{\rho}}
\def\bzeta{\vvm{\zeta}}
\def\bu{\vv{u}}
\def\bv{\vv{v}}
\def\bw{\vv{w}}
\def\bz{\vv{z}}
\def\fv{\vv{f}}
\def\bg{\vv{g}}
\def\AaK{\vv{A}_{\text{dev}}^K}
\def\AbK{\vv{A}_{\text{sta}}^K}
\def\AcK{\vv{A}_{\text{div}}^K}
\def\AdK{\vv{A}_{\text{tra}}^K}
\def\CK{\vv{a}_{\,\text{tra}}^K}
\def\BK{\vv{B}^K}
\def\DaK{\vv{D}_{\text{gra}}^K}
\def\DbK{\vv{D}_{\text{sta}}^K}
\def\DcK{\vv{D}_{\text{bou}}^K}
\def\EK{\vv{C}^K}
\def\mo{m_0}
\def\m1{m_1}
\def\Mk{\vv{M}_{\text{mass},e}}
\def\iMk{\Mk^{-1}}
\def\Me{\vv{M}_{e,K}}
\def\MASS{\vv{M}_{\text{mass},K}^{(k+1)}}
\def\Mass{\vv{M}_{\text{mass},K}}

\def\tM{\vv{M}_{\text{grad},K}}
\def\Mo{\vv{M}_{0,K}}
\def\Ao{\vv{A}_0}
\def\CLe#1{\vv{C}_L^{e_{#1}}}
\def\NLe#1{\vv{N}_L^{e_{#1}}}
\def\ML{\vv{M}_{\text{lap},K}}
\def\MCBe{\vv{M}_{\text{Lag},e}}
\def\MassL{\vv{M}_{\text{mass},L}}
\def\Pb{\vv{P}^K}
\def\ZZ{\vv{P}_{\text{grad}}}
\def\Be{\vv{B}_e}
\def\Di{\vv{M}_{\text{div}}}
\def\Pr{\vv{P}_1}
\def\IMASS{{\vv{I}}_{\text{Mass}}}
\def\Pbt{\widehat{\vv{P}}_0}
\def\Rb{\vv{R}^K}
\def\BR{\vv{B}^R}
\def\tR{\widehat{\vv{R}}_0}
\def\Ha{\vv{H}_A}
\def\Ee#1{\vv{H}_2^{#1}}
\def\Hd{\vv{H}_D}
\def\PV{\vv{P}_{\text{eval}}}
\def\PU{\vv{P}_{V}}
\def\PGU{\mathcal{D}\vv{P}_{V}}
\def\DxDy{\vv{P}_{\text{div}}}
\def\PLe{\vv{P}_{L}^e}
\def\Pg#1{\vv{p}_{g_{#1}}^e}
\def\Pf#1{\vv{p}_{f_{#1}}^K}
\def\MZ#1{\vv{M}_{z_{#1}}^{(s)}}
\def\MdPU{\vv{M}_u^{(s)}}
\def\IMd0I{\vv{I}_{\text{per}}}
\numberwithin{equation}{section}

\usepackage{amsthm}

\def\thanksFila{Escuela de Matem\'atica, Universidad Nacional, Campus Omar Dengo,
				Heredia, Costa Rica, email: {\tt filander.sequeira@una.cr}.}
\def\thanksHelen{Escuela de Matem\'atica, Universidad Nacional, Campus Omar Dengo,
				Heredia, Costa Rica, email: {\tt hellen.guillen.oviedo@una.ac.cr}.}

\begin{document}

\title{Some aspects on the computational implementation of diverse terms arising
in mixed virtual element formulations}

\author{{\sc Fil\'ander A. Sequeira}\thanks{\thanksFila}\qquad
        {\sc Helen Guill\'en-Oviedo}\thanks{\thanksHelen}}

\date{}

\maketitle

\begin{abstract}
\noindent
In the present paper we describe the computational implementation of
some integral terms that arise from mixed virtual element methods
(mixed-VEM) in two-dimensional pseudostress-velocity formulations.
The implementation presented here consider any polynomial degree
$k \geq 0$ in a natural way by building several local matrices of
small size through the matrix multiplication and the Kronecker
product. In particular, we apply the foregoing mentioned matrices to
the Navier-Stokes equations with Dirichlet boundary conditions, whose
mixed-VEM formulation was originally proposed and analyzed in 
a recent work using virtual element subspaces for $H(\div)$ and $H^1$,
simultaneously. In addition, an algorithm is proposed for the
assembly of the associated global linear system for the Newton's
iteration. Finally, we present a numerical example in order to
illustrate the performance of the mixed-VEM scheme and confirming
the expected theoretical convergence rates.
\end{abstract}

\noindent
{\bf Key words:} mixed virtual element method, high-order approximations,
computational implementation, Navier-Stokes problem, pseudostress-velocity
formulation, augmented formulation


\section{Introduction}\label{sec:intro}

The virtual element method (VEM) was introduced in \cite{firstVEM} for the
Poisson equation, which is one of the high-order discretization schemes for
the approximation of solutions to partial differential equations, that can be
seen as a generalization of the standard finite element method. In other words, 
the method extends the classical finite element technique to general polygonal
and polyhedral meshes. Moreover, according to \cite{bbmr-2014}, additional
advantages of VEM schemes, when compared with finite volume methods, mimetic
finite difference method, and related techniques, are given by its solid mathematical
ground, the simplicity of the respective computational coding, and the quality
of the numerical results provided. On the other hand, regarding to purely mixed
virtual element techniques, that is based on dual-mixed variational formulations,
the method was initially developed in \cite{bfm-M2AN-2014}, and more recently 
extended in \cite{bbmr-NM-2016}, \cite{bbmr-2016-ESAIM}, \cite{cg-IMANUM-2017},
\cite{cgs-1}, and \cite{gms-vem-ns}. In particular, edge and face VEM spaces in
2D and 3D, which together with the nodal and volume spaces constitute a discrete
complex, were developed in \cite{bbmr-NM-2016}, whereas \cite{bbmr-2016-ESAIM}
generalizes the results of \cite{bbmr-NM-2016} to the case of variable coefficients.
In turn, \cite{cg-IMANUM-2017} and \cite{gms-vem-ns} provide the first analysis of
a virtual element method for a mixed variational formulation of the Stokes and
Navier-Stokes problems, respectively, in which the pseudostress and the velocity
are the only unknowns, whereas the pressure is computed via a postprocessing formula. 
Additionally, the analysis presented in \cite{gms-vem-ns} allows to study problems
of the same nature as Navier-Stokes, such as the Boussinesq problem, where a mixed
method of virtual elements was introduced and analyzed in \cite{gms-brinkman}.
For several other contributions on VEM and mixed-VEM we refer for instance to
\cite{abmv-SINUM-2014}, \cite{bbm-SINUM-2013}, \cite{bm-CMAME-2013},
\cite{gtp-CMAME-2014}, and \cite{ms-aposteriori-2020}.

The previously mentioned references omit to present aspects related to computational
implementation of mixed-VEM schemes. For this reason, the main goal of this paper is to
describe a computational approach for mixed-VEM methods in order to obtain high-order
approximations, without imposing a certain programming language. More precisely, the
reason for this contribution is that there is a few literature that explains how to program
the subspaces of virtual elements accurately (particularly the mixed ones). Up to the authors
knowledge, some works such as \cite{firstVEM}, \cite{bbmr-NM-2016}, \cite{bfm-M2AN-2014},
\cite{bbmr-2014}, \cite{gtp-CMAME-2014}, \cite{blv-2018} and \cite{tesisErnesto}, 
explain specific aspects of the computational implementation, but in general they do not
describe the structures employed. Conversely, in \cite{fiftylines-2017} the authors
present the first paper concerning a detailed implementation of virtual element
method. However, this contribution focus only on the lowest order approximations for Poison
equation. Next, in \cite{helen-fila} the authors describe in spanish some specific aspects
on the computational implementation of the a mixed-VEM method for the 2D linear Brinkman
model proposed and analyzed in \cite{cgs-1}. In fact, the present paper extends the
approach used in \cite{cgs-1}, in order to obtain implementation techniques for several
mixed-VEM schemes, including those with nonlinearities (see, e.g., \cite{gms-vem-ns}).


The paper is organized as follows. In Section \ref{sec:vem-spaces} we
introduce the virtual element subspaces for $H(\div)$-conforming and
$H^1$-conforming that will be employed. This includes the main ingredients
for the polygonal mesh structure, the definitions of the local degrees
of freedom, and the projections to be employed, along with a description
about the explicit calculation of each projector. Next, in Sections \ref{sec:projectionMat}
and \ref{sec:operators} we present the main contributions of this work.
Indeed, in Section \ref{sec:projectionMat} we describe the assemble
of the local matrices associated with the projectors respect to the
local virtual spaces, whereas in Section \ref{sec:operators} the
computational aspects required for the construction of some local terms
arising in mixed-primal virtual element formulations are described. 
In addition, we remark in advance that each discrete operator is built
for an arbitrary polynomial degree $k\geq 0$, which means that we will
develop a high-order computational approach.
Finally, in Section \ref{sec:navier-stokes}, in order to illustrate the
use of the matrices introduced in previous sections, we recall the boundary
value problem and its mixed-VEM formulation introduced and analyzed in \cite{gms-vem-ns}.
More precisely, we present a mixed virtual element method for the two-dimensional
pseudostress-velocity formulation of the Navier-Stokes equations with Dirichlet
boundary conditions. Therein, the continuous and discrete formulations are presented.
Furthermore, we propose an algorithm for the assembly of the associated global linear
system for the Newton's iteration, and then, a numerical example illustrating
the performance of the mixed-VEM scheme and confirming the expected theoretical
convergence rates is presented.

\subsection*{Notations}

We end the present section by providing some notations to be used along the paper.
Indeed, in what follows we consider a bounded domain $\Omega\subseteq\R^2$ with
boundary $\Gamma$. Moreover, standard terminology for Lebesgue and Sobolev spaces
will be adopted, where given a generic scalar functional space $\mathrm H$, 
we denote by $\mathbf{H}$ and $\mathbb{H}$ be the corresponding vectorial and
tensorial counterparts, respectively. For example, given an integer $\ell\geq 0$
and $U\subseteq\R^2$, we let $\P_{\ell}(U)$ be the space of polynomials on $U$
of degree up to $\ell$, whereas $\Pv_{\ell}(U)$ stands for its vectorial version,
that is, $\Pv_{\ell}(U) := [\P_{\ell}(U)]^2$. In addition, $\Pm_{\ell}(U) := [\P_{\ell}(U)]^{2\times 2}$
corresponds to its tensorial version.

Now, we employ $\vv{v}\otimes\vv{w}$ to stand the usual dyadic product for
two column vectors $\vv{v},\vv{w}\in\R^2$, that is, $\vv{v}\otimes\vv{w} := \vv{v}\,\vv{w}^{\tra} \in\R^{2\times 2}$.
On the other hand, given two matrices $\vv{A}\in\R^{m\times n}$ and $\vv{B}\in\R^{p\times q}$,
we denote the matrix concatenation of $\vv{A}$ and $\vv{B}$ as follows:
$$\left[\begin{array}{c:c}
\vv{A} & \vv{B}
\end{array}\right]\quad\text{if } m=p\,,\qquad\text{or}\qquad \left[\begin{array}{c}
\vv{A}\\ \hdashline \vv{B}
\end{array}\right]\quad\text{if } n=q\,.$$
In addition, we let $\kron(\cdot,\cdot)$ be the usual Kronecker product,
that is:
$$\kron(\vv{A},\vv{B})\ :=\ \left[\begin{array}{c:c:c}
a_{11}\vv{B} & \cdots & a_{1n}\vv{B}\\[-2ex]
\\ \hdashline
\vdots & \ddots & \vdots \\[-2ex]
\\ \hdashline
\\[-2ex]
a_{m1}\vv{B} & \cdots & a_{mn}\vv{B}
\end{array}\right]\in\R^{(mp)\times(nq)}\,.$$
Finally, when we write $\vv{A} := [a_{ij}]\in\R^{m\times n}$,
it means that $i=1,2,\ldots,m$ and $j=1,2,\ldots,n$.

\section{The virtual element subspaces}\label{sec:vem-spaces}

In this section we recall two local element subspaces, usually used in the
design of virtual element schemes. More precisely, we introduce the spaces $\Hv_k^K$
(see \eqref{eqn:defHhv}) and $V_k^K$ (see \eqref{eqn:defVhv}) for  $H(\div)$-conforming
and $H^1$-conforming elements, respectively. In order to do that, we let $\{\Th\}_{h>0}$
be a family of decompositions of $\Omega$ in polygonal elements, where $h$ denotes the
largest of its diameters. For theoretical purposes, it is assumed that $\Th$ satisfies
some conditions described at the beginning of \cite[Section 3.1]{gms-vem-ns}.

In what follows, we consider $K$ an arbitrary element of $\Th$, which it is important
to realize that the decomposition $\Th$ needs a quite more sophisticated computational
structure than those used in classical finite element methods. Indeed, we recall here
that $\Th$ can contain elements with several shapes that affect important aspects,
such as: the number of edges and the calculation of its diameter. In particular,
any structure that is implemented for mesh management, from a connectivity
point of view, must be able to indicate:
	\begin{itemize}
	\item the number of nodes (points)
	\item the number of edges
	\item the number of elements
	\item the number of boundary edges
	\item for any element:
		\begin{itemize}
		\item the number of nodes of the element
		\item the global index (identifier) of a node of the element
		\item the local indexes of the nodes in a specific edge
		\item the global indexes of the nodes in a specific edge
		\item the global index of an edge of the element
		\end{itemize}
	\item for any edge:
		\begin{itemize}
		\item the global index of a node of the edge
		\item the global index of the edge
		\item the global index of the neighbor element shared by the edge
		\item the global indexes of the elements that contain the edge
		\item the orientation in a specific elements that contain the edge
		\item an identifier that establishes if the edge is in the boundary
		\end{itemize}
	\end{itemize}
Here, the orientation of an edge corresponds to a boolean identifier that indicates
the ordering (independent of the element that contains it) of its extreme nodes.
In addition, in the case of geometric aspects, for all $K\in\Th$, it
must be able to calculate:
	\begin{itemize}
	\item the number of vertices, or equivalently, the number of edges
	\item the coordinates of a specific vertex
	\item the coordinates of the barycenter of $K$
	\item the area of $K$
	\item the diameter of $K$
	\item the midpoints, normal vectors and lengths of each of its edges
	\end{itemize}

\subsection{$\Hv(\div;K)$-conforming subspace, associated bilinear form and projection}\label{sec:dofH}

Let $e$ be an edge of $\Th$ with midpoint $x_e$ and length $h_e$. Thus, given
an integer $\ell \geq 0$, we consider the following set of $k+1$ normalized
monomials on $e$:
\begin{equation}\label{eqn:defBe}
\mathcal{B}_{\ell}(e)\ :=\ \left\{\left(\frac{x-x_e}{h_e}\right)^j\right\}_{0\,\leq\, j\,\leq\, \ell}\,,
\end{equation}
which constitutes a basis of $\P_{\ell}(e)$. Similarly, given an element $K\in\Th$ with barycenter
$\vv{x}_K$ and diameter $h_K$, we define the following set of $\frac{1}{2}(\ell+1)(\ell+2)$
normalized monomials on $K$:
\begin{equation}\label{eqn:defBk}
\mathcal{B}_{\ell}(K)\ :=\ \left\{\left(\frac{\vv{x}-\vv{x}_K}{h_K}\right)^{\vvm{\alpha}}\right\}_{0\,\leq\, 
|\vvm{\alpha}|\,\leq\, \ell}\,,
\end{equation}
which is a basis of $\P_{\ell}(K)$. It is important to remark that in \eqref{eqn:defBk}
we use the multi-index notation, where given $\vv{x} := (x_1,x_2)^{\tra}\in\R^2$
and $\vvm{\alpha}:=(\alpha_1,\alpha_2)^{\tra}$, with non-negative integers $\alpha_1$,
$\alpha_2$, we let $\vv{x}^{\vvm{\alpha}} := x_1^{\alpha_1}x_2^{\alpha_2}$ and
$|\vvm{\alpha}| := \alpha_1+\alpha_2$.

Next, we introduce the auxiliary local virtual element space of order $k\geq 0$
(see, e.g., \cite{bbmr-NM-2016,bbmr-2016-ESAIM})
\begin{equation}\label{eqn:defHhv}
\begin{array}{l}
\Hv_k^K\ :=\ \Big\{\tau := (\tau_1,\tau_2)^{\tra} \in\Hv(\div;K)\cap\Hv(\rot;K)\,:\quad \tau\cdot\nv|_e\in\P_k(e)\\[1ex]
\phantom{\Hv_k^K\ :=\ \Big\{}\forall \mbox{ edge } e\in\partial K\,,\quad
\div(\tau)\in\P_k(K)\,, \;\text{ and }\; \rot(\tau)\in\P_{k-1}(K)\Big\}\,,
\end{array}
\end{equation}
where $\rot(\tau) := \frac{\partial \tau_{2}}{\partial{x_1}}-\frac{\partial\tau_1}{\partial{x_2}}$
and $\P_{{-}1}(K) := \{0\}$. 
Moreover, the local degrees of freedom for $\tau\in\Hv_k^K$ are
given by (see \cite{bbmr-NM-2016,bbmr-2016-ESAIM})
\begin{eqnarray}\label{eqn:dof-Hkv}
\begin{array}{rcl}
m_{q,\nv}^H(\tau)  & := & \displaystyle\int_e\btau\cdot\nv\,q\qquad\;\; \forall\ q\in\mathcal{B}_k(e)\,,\quad \forall\ \mbox{edge}\ e\in\partial K\,,\\[3ex]
m_{q,\div}^H(\tau) & := & \displaystyle\int_K\tau\cdot \nabla q\qquad\; \forall\ q\in\mathcal{B}_k(K)\setminus\{1\}\,,\\[3ex]
m_{\vv{q},\rot}^H(\tau) & := & \displaystyle\int_K\tau\cdot\vv{q}\qquad\quad \forall\ \vv{q}\in\Gk\,,
\end{array}
\end{eqnarray}
where $\Gk$ is a basis of $(\nabla\P_{k+1}(K))^{\perp}\cap\Pv_k(K)$,
which corresponds to the $\Lv(K)$-orthogonal of $\nabla\P_{k+1}(K)$
in $\Pv_k(K)$. Then, according to the cardinalities of $\mathcal{B}_k(e)$
and $\mathcal{B}_k(K)$, and the dimensions of $\Pv_k(K)$ and $\nabla\P_{k+1}(K)$,
it follows that the cardinality of $\Gk$ is $\frac{1}{2}k(k+1)$. Thus, the number
of local degrees of freedom defined in \eqref{eqn:dof-Hkv} (i.e. the dimension of
$\Hv_k^K$) is given by:
\begin{eqnarray}
n_k^H\ =\ n_k^H(K) & := & (k+1)\,d_K\ +\ \left\{\frac{(k+1)(k+2)}{2} - 1\right\}\ +\ \frac{k(k+1)}{2}\nonumber\\
& = & (k+1)(d_K + k + 1) - 1\,,\label{eqn:n1}
\end{eqnarray}
where $d_K$ corresponds to the number of edges in $K$. Furthermore,
it was proved in \cite[Section 3.4]{bbmr-NM-2016} that, for every
$K\in\Th$, these $n_k^H$ local degrees of freedom are unisolvent
in $\Hv_k^K$.

On the other hand, we employ the space $\Hv_k^K$ to define a tensor virtual element
space $\Hm_k^K$ as:
\begin{equation}\label{eqn:defHh}
\Hm_k^K\ =\ \Big\{\btau\in\Hm(\bdiv;K)\,:\quad
(\tau_{i1},\, \tau_{i2})^{\tra}\in\Hv_k^K\quad\forall\ i\in\{1,2\}\Big\}\,.
\end{equation}
In other words, $\Hm_k^K$ is a subspace of $\Hm(\bdiv;K)$, where each row of
$\btau\in\Hm_k^K$ belongs to $\Hv_k^K$. According to this, it is natural
to consider for each $\btau\in\Hm_k^K$, the following $2n_k^H$ local degrees
of freedom:
\begin{eqnarray}\label{eqn:dof-Hk}
\begin{array}{rcl}
\vvm{m}_{\vv{q},\nv}^H(\btau)     & := & \displaystyle\int_e\btau\nv\cdot\vv{q}\qquad\;\; \forall\ \vv{q}\in\vvm{\mathcal{B}}_k(e)\,,\quad \forall\ \mbox{edge}\ e\in\partial K\,,\\[3ex]
\vvm{m}_{\vv{q},\bdiv}^H(\btau)   & := & \displaystyle\int_K\btau : \nabla\vv{q}\qquad \forall\ \vv{q}\in\vvm{\mathcal{B}}_k(K)\setminus\{(1,0)^{\tra},(0,1)^{\tra}\}\,,\\[3ex]
\vvm{m}_{\brho,\brot}^H(\btau) & := & \displaystyle\int_K\btau : \brho\qquad\quad \forall\ \brho\in\Gk\,,
\end{array}
\end{eqnarray}
where
\begin{eqnarray*}
\vvm{\mathcal{B}}_{\ell}(e) & := & \big\{ (q,0)^{\tra}\,:\, q\in \mathcal{B}_{\ell}(e)\big\} \,\cup\,
\big\{ (0,q)^{\tra}\,:\, q\in \mathcal{B}_{\ell}(e)\big\}\,,\nonumber\\[2ex]
\vvm{\mathcal{B}}_{\ell}(K) & := & \big\{ (q,0)^{\tra}\,:\, q\in \mathcal{B}_{\ell}(K)\big\} \,\cup\,
\big\{ (0,q)^{\tra}\,:\, q\in \mathcal{B}_{\ell}(K)\big\}\,,\label{eqn:defBvk}
\end{eqnarray*}
and
$$\Gkb\ :=\ \left\{\mat{\begin{array}{c}
\vv{q} \\[0.2ex] \vv{0}
\end{array}}\,:\, \vv{q} \in \Gk\right\} \,\cup\,
\left\{\mat{\begin{array}{c}
\vv{0} \\[-0.5ex] \vv{q}
\end{array}}\,:\, \vv{q} \in \Gk\right\}\,.$$
Furthermore, according to the degrees of freedom defined in \eqref{eqn:dof-Hkv}
and \eqref{eqn:dof-Hk}, it is possible to define a bilinear form
$\mathcal{S}_H^K:\Hm_k^K\times\Hm_k^K\to\R$ based on these.
Indeed, let $K\in\Th$ and consider the union of all local degrees of freedom
(cf. \eqref{eqn:dof-Hkv}) of a given $\tau\in\Hv^1(K)$, in a set
$\{m_{i,K}^H(\tau)\}_{i=1}^{n_k^H}$. Then, as usual, let $\{\psiH_j\}_{j=1}^{n_k^H}$
be the canonical basis of $\Hv_k^K$. That is, given $i=1,2,\ldots, n_k^H$,
$\psiH_i$ is the unique element in $\Hv_k^K$ such that:
\begin{equation}\label{eqn:basisH}
m_{j,K}^H(\psiH_i)\ =\ \delta_{ij}\qquad\forall\ j=1,2,\ldots,n_k^H\,,
\end{equation}
where, in particular, there holds:
$$\tau \ =\ \sum_{j=1}^{n_k^H}m_{j,K}^H(\tau)\,\psiH_j \qquad \forall\,\tau \in \Hv^K_k\,.$$

Now, let $s_H^K : \Hv_k^K\times\Hv_k^K\rightarrow\R$ be
the bilinear form associated with the identity matrix in $\R^{n_k^H\times n_k^H}$,
respect to the basis $\{\psiH_j\}_{j=1}^{n_k^H}$ of $\Hv_k^K$. More precisely, we
have:
\begin{equation*}\label{eqn:defskh}
s_H^K(\zeta,\tau)\ :=\ \sum_{i=1}^{n_k^H}m_{i,K}^H(\zeta)\,m_{i,K}^H(\tau)\qquad\forall \, \zeta, \, \tau\in\Hv_k^K\,.
\end{equation*}
Then, we define $\mathcal{S}_H^K:\Hm_k^K\times\Hm_k^K\rightarrow\R$ as the bilinear
form associated with the degrees of freedom of $\Hm_k^K$ as follows:
\begin{equation}\label{eqn:defSH}
\mathcal{S}_H^K(\bzeta,\,\btau)\ :=\ \sum_{i=1}^2s_H^K\big((\zeta_{i1},\zeta_{i2})^{\tra}\,, (\tau_{i1},\tau_{i2})^{\tra}\big)\,,
\end{equation}
for all $\bzeta:=(\zeta_{ij}), \, \btau := (\tau_{ij})\in\Hm_k^K$.

In addition, let $P_k^K : L^2(K)\rightarrow\P_k(K)$ be corresponding orthogonal
projection, such that, for $v\in L^2(K)$, it is characterized by:
\begin{equation}\label{eqn:defL2Proj}
P_k^K(v) \in \P_k(K) \qan
\int_KP_k^K(v)\,q\ =\ \int_K v\,q \qquad\forall\, q \in\P_k(K)\,.
\end{equation}
In turn, let $\Pk : \Lv(K)\rightarrow\Pv_k(K)$ its corresponding vectorial version,
such that, for $\bv\in\Lv(K)$ there hold:
\begin{equation}\label{eqn:defL2Proj-v}
\Pk(\bv) \in \Pv_k(K) \qan
\int_K\Pk(\bv)\cdot\vv{q}\ =\ \int_K\bv\cdot\vv{q}\qquad\forall\, \vv{q}\in\Pv_k(K)\,,
\end{equation}
where, notice that $\Pk(\bv) = \big( P_k^K(v_1), P_k^K(v_2)\big)^{\tra}$
for all $\bv := (v_1,v_2)^{\tra} \in \Lv(K)$.

Next, for the following sections, it is important to mention that, according
to \cite[Section 3.2]{bbmr-2016-ESAIM}, the degrees of freedom given in
\eqref{eqn:dof-Hkv} allow us the explicit calculation of $\Pk(\tau)$ for every
$\tau\in\Hv_k^K$. That is, it is possible to determine the $\Lv(K)$-orthogonal projector
for elements in the virtual space. Indeed, it is sufficient to verify that the
right-hand side in the second expression of \eqref{eqn:defL2Proj-v} can be calculated in these cases.
To do that, notice from the definitions of $m_{q,\nv}^H(\tau)$ and $m_{q,\div}^H(\tau)$
(cf. \eqref{eqn:dof-Hkv}), that it is possible to determine the value of $\div(\tau)\in\P_k(K)$
using the identity:
\begin{equation}\label{eqn:compute-div}
\int_K\div(\tau)\,q \ =\ -\int_K\tau\cdot\nabla q \ +\ \int_{\partial K} \tau \cdot\nv\,q \qquad \forall \, q\in\P_k(K)\,.
\end{equation}
Moreover, given $\vv{q}\in\Pv_k(K)$, it is well known that there exist unique
$\vv{q}^{\perp}\in (\nabla\P_{k+1}(K))^{\perp}\cap\Pv_k(K)$ and $\widetilde{q}\in\P_{k+1}(K)$,
such that: $\vv{q} = \vv{q}^{\perp} + \nabla\widetilde{q}$.
In this sense, it follows that:
\begin{equation}\label{eqn:proyP}
\int_K\tau\cdot\vv{q}\ =\ \int_K\tau \cdot\vv{q}^{\perp} \,+\, 
\int_K \tau \cdot\nabla\widetilde{q}\ =\ \int_K\tau \cdot\vv{q}^{\perp} \,-\, 
\int_K\widetilde{q}\,\div(\tau) \,+\, \int_{\partial K} \tau \cdot \nv\,\widetilde{q}\,,
\end{equation}
which, in accordance with \eqref{eqn:compute-div} and the definition of $m_{\vv{q},\rot}^H(\tau)$
(see \eqref{eqn:dof-Hkv}), allow the required calculation.

Finally, we also consider the $\Lm(K)$-orthogonal projection $\bPk : \Lm(K)\rightarrow\Pm_k(K)$.
In other words, $\bPk$ is the operator $\Pk$ (cf. \eqref{eqn:defL2Proj-v}) acting
on each row of a tensor of $\Lm(K)$, which according to the foregoing discussion
it quite simple to see that $\bPk(\btau)$ can be explicitly calculated for each
$\btau\in\Hm_k^K$.

\subsection{$H^1(K)$-conforming subspace, associated bilinear form and projection}\label{sec:dofV}

We follow the previous section by defining the following local virtual
element space of order $k\geq 0$ (see, e.g., \cite{aabmr-CMA-2013})
\begin{equation}\label{eqn:defVhv}
\begin{array}{l}
V_k^K\ :=\ \bigg\{v\in H^1(K)\,:\quad v|_{\partial K}\in E_{k+1}(K)\,,\quad \Delta v\in\P_{k+1}(K)\,,\\[1ex]
\phantom{\vv{V}_k^K\ :=\ \Big\{}\text{and }\;\displaystyle\int_K\big\{R_k^K(v)-v\big\}q\, =\, 0 \quad\forall\ q\in\widetilde{\mathcal{B}}_k(K)\bigg\}\,,
\end{array}
\end{equation}
where $E_{k+1}(K) := \big\{v\in C(\partial K)\,:\, v|_e\in\P_{k+1}(e)\quad\forall\ \text{edge }e\subseteq\partial K\big\}$,
$\widetilde{\mathcal{B}}_0(K) := \mathcal{B}_1(K)$, and $\widetilde{\mathcal{B}}_k(K) := \mathcal{B}_{k+1}(K)\setminus\mathcal{B}_{k-1}(K)$
for $k \geq 1$. In addition, $R_k^K : H^1(K)\to\P_{k+1}(K)$ is the projection
operator defined for each $v\in H^1(K)$ as the unique polynomial $R_k^K(v)\in\P_{k+1}(K)$
such that (see \cite{bbmr-M3AS-2016})
\begin{equation}\label{eqn:defRh-0}
\begin{array}{rcl}
\displaystyle\int_K\nabla R_k^K(v)\cdot\nabla q & = & \displaystyle\int_K\nabla v\cdot\nabla q\qquad\forall\ q\in\P_{k+1}(K)\,,\\[3ex]
\displaystyle\int_{U}R_k^K(v) & = & \displaystyle\int_{U}v\,,
\end{array}
\end{equation}
where $U = \partial K$ if $k = 0$, and $U = K$ if $k\geq1$.
Now, recalling from \cite{aabmr-CMA-2013} the following degrees of freedom
for a given $v\in V_k^K$
\begin{equation}\label{eqn:dof-Vkv}
\begin{array}{rcl}
m^{V}_{i,v}(v) & := & \text{value of } v \text{ at the } i\text{th vertex of }K\,,\quad\forall\ i\text{ vertex of } K\,, \\[2ex]
m^{V}_{e}(v)   & := & \text{values of } v \text{ at } k\text{ uniformly spaced points on } e\,, \,\, \forall\ e\in\partial K\,, \,\, \text{for } k\geq 1\,,\\[1ex]
m^{V}_{q,K}(v) & := & \text{value of }\displaystyle\int_Kvq\,,\,\, \forall\ q\in\mathcal{B}_{k-1}(K)\,,
\,\, \text{for } k\ge 1\,,
\end{array}
\end{equation}
it easy to check that the dimension of $V_k^K$ is given by
\begin{eqnarray}
n_k^V\ =\ n_k^V(K) & := & d_K\ +\ k\,d_K\ +\ \frac{k(k+1)}{2}\nonumber\\
& = & (k+1)d_K\ +\ \frac{k(k+1)}{2}\,.\label{eqn:n2}
\end{eqnarray}
In addition, from \cite[Propositions 1 and 2]{aabmr-CMA-2013} we know
that \eqref{eqn:dof-Vkv} are unisolvent in $V_k^K$.

Next, we now let $\vv{V}_k^K$ be the vectorial version of $V_k^K$
given by:
\begin{equation}\label{eqn:defVh}
\vv{V}_k^K\ =\ \Big\{\bv:=(v_1,v_2)^{\tra}\in\Hv^1(K)\,:\quad v_i\in V_k^K\quad\forall\ i\in\{1,2\}\Big\}\,,
\end{equation}
which, in particular, satisfies that $\dim\vv{V}_k^K = 2n_k^V$. Moreover,
let $\Rk : \Hv^1(K)\to\Pv_{k+1}(K)$ be the vectorial version of the
operator $R_k^K$ (cf. \eqref{eqn:defRh-0}) as
\begin{equation*}\label{eqn:defRh}
\begin{array}{rcl}
\displaystyle\int_K\nabla\Rk(\bv):\nabla\vv{q} & = & \displaystyle\int_K\nabla\bv:\nabla\vv{q}\qquad\forall\ \vv{q}\in\Pv_{k+1}(K)\,,\\[3ex]
\displaystyle\int_{U}\Rk(\bv) & = & \displaystyle\int_{U}\bv\,,
\end{array}
\end{equation*}
where $U = \partial K$ if $k = 0$, and $U = K$ if $k\geq1$, which
allow us to rewrite $\vv{V}_k^K$ in the form:
\begin{equation*}
\begin{array}{l}
\vv{V}_k^K\ :=\ \bigg\{\bv\in\Hv^1(K)\,:\quad \bv|_{\partial K}\in\vv{E}_{k+1}(K)\,,\quad \Delta\bv\in\Pv_{k+1}(K)\,,\\[1ex]
\phantom{\vv{V}_k^K\ :=\ \Big\{}\text{and }\;\displaystyle\int_K\big\{\Rk(\bv)-\bv\big\}\cdot\vv{q}\, =\, 0 \quad\forall\ \vv{q}\in\widetilde{\vvm{\mathcal{B}}}_k(K)\bigg\}\,,
\end{array}
\end{equation*}
where $\vv{E}_{k+1}(K) := [E_{k+1}(K)]^2$, $\widetilde{\vvm{\mathcal{B}}}_0(K) := \vvm{\mathcal{B}}_1(K)$, and
$\widetilde{\vvm{\mathcal{B}}}_k(K) := \vvm{\mathcal{B}}_{k+1}(K)\setminus\vvm{\mathcal{B}}_{k-1}(K)$
for $k \geq 1$.

Furthermore, we now denote by $\{m^V_{j,K}(v)\}_{j=1}^{n^V_k}$ the degrees of freedom 
defined by \eqref{eqn:dof-Vkv}, and let $s_V^K : V_k^K\times V_k^K\rightarrow\R$ be the
associated bilinear form:
\begin{equation*}\label{eqn:defskv}
s_V^K(w,v)\ :=\ \sum_{i=1}^{n_k^V}m_{i,K}^V(w)\,m_{i,K}^V(v)\qquad\forall \, w, \, v \in V_k^K\,,
\end{equation*}
which allows us to define the bilinear form $\mathcal{S}_V^K:\vv{V}_k^K\times\vv{V}_k^K\rightarrow\R$
as follows:
\begin{equation}\label{eqn:defSV}
\mathcal{S}_V^K(\bw,\bv)\ :=\ \sum_{i=1}^2s_V^K\big(w_i\,, v_i\big)\,,
\end{equation}
for all $\bw := (w_1,w_2)^{\tra}, \, \bv:=(v_1,v_2)^{\tra}\in\vv{V}_k^K$.
On the other hand, we introduce $\{\psiV_j\}_{j=1}^{n_k^V}$ as the canonical
basis of $V_k^K$ such that
\begin{equation}\label{eqn:basisV}
m_{j,K}^V(\psiV_i)\ =\ \delta_{ij}\qquad\forall\ j=1,2,\ldots,n_k^V\,,
\end{equation}
for a given $i=1,2,\ldots, n_k^V$. In particular, there holds:
$$v \ =\ \sum_{j=1}^{n_k^V}m_{j,K}^V(v)\,\psiV_j \qquad \forall\,v \in V^K_k\,.$$

We end this section by clarifying that, for each $\bv\in\vv{V}_k^K$, its projections $\Rk(\bv)$, $\Pk(\bv)$ and $\bPk(\nabla\bv)$
can be computed explicitly by using the degrees of freedom defined in \eqref{eqn:dof-Vkv}.
Indeed, using \eqref{eqn:defVh} it is enough, for each $v\in V_k^K$, to describe
how use the degrees of freedom \eqref{eqn:dof-Vkv} to compute $R_k^K(v)$ (cf. \eqref{eqn:defRh-0}),
$P_k^K(v)$ (cf. \eqref{eqn:defL2Proj}), and $\Pk(\nabla v)$ (cf. \eqref{eqn:defL2Proj-v}),
respectively. Indeed, we begin by noticing that, for $v\in V_k^K$ and $q\in\P_k(K)$,
the right-hand side of the first equation of \eqref{eqn:defRh-0} can be integrated
by parts to yield
\begin{equation}\label{eqn:proyR}
\int_K\nabla v\cdot \nabla q\ =\ -\int_K v\,\Delta q\ +\ \int_{\partial K}(\nabla q\cdot\nv)\,v\,,
\end{equation}
where, since $\Delta q\in\P_{k-2}(K)$ and $\nabla q\cdot\nv\in\P_{k-1}(K)$, the first
integral on the right-hand side can be computed by using the degrees of freedom
$m^{V}_{q,K}(v)$, whereas for the second one using $m^{V}_{i,v}(v)$ and $m^{V}_{e}(v)$.
Finally, for the right-hand side of the second equation of \eqref{eqn:defRh-0},
it is straightforward to see that $\int_{\partial K}v = \int_{\partial K}v\cdot 1$
can be calculated using again $m^{V}_{i,v}(v)$ and $m^{V}_{e}(v)$, whereas 
$\int_{K}v = \int_{K}v\cdot 1$ utilizing $m^{V}_{q,K}(v)$ for $k\geq 1$.

Similarly, integrating by parts we observe that
\begin{equation}\label{eqn:proyPDU}
\begin{array}{rcl}
\displaystyle\int_K\nabla v\cdot \vv{q} & = & \displaystyle -\int_K v\,\div(\vv{q})\ +\ \int_{\partial K}(\vv{q}\cdot\nv)\,v\\[3ex]
& = & \displaystyle -\int_KP_k^K(v)\,\div(\vv{q})\ +\ 
\int_{\partial K}(\vv{q}\cdot\nv)\,v\qquad\forall\ \vv{q}\in\Pv_k(K)\,,
\end{array}
\end{equation}
which yields the explicit computation of $\Pk(\nabla v)$ for all $v\in V_k^K$.

Finally, for each $v\in V_k^K$, the right-hand side of \eqref{eqn:defL2Proj}
can be computed using the degrees of freedom given by $m^{V}_{q,K}(v)$ (cf. \eqref{eqn:dof-Vkv}).
Indeed, given $q\in\P_k(K)$ we can write $q = \widehat{q} + \widetilde{q}$ such
that $\widehat{q}\in\mathcal{B}_{k-1}$ and $\widetilde{q}\in\widetilde{\mathcal{B}}_k$.
Thus, $\int_Kv\,\widehat{q}$ can be computed using $m^{V}_{q,K}(v)$, whereas recalling
from \eqref{eqn:defVhv} that
\begin{equation}\label{eqn:proyPU}
\int_Kv\,\widetilde{q}\ =\ \int_KR_k^K(v)\,\widetilde{q}\,,
\end{equation}
we can compute $\int_Kv\,\widetilde{q}$, since $R_k^K(v)$ is explicitly
computable for each $v\in V_k^K$.

\section{Matrices associated with the projectors}\label{sec:projectionMat}

We now aim to describe the explicit calculation of some projections of
elements of the virtual subspaces defined in previous section. We beginning by
remarking that the bases $\{\psiH_j\}_{j=1}^{n_k^H} $ and $\{\psiV_j\}_{j=1}^{n_k^V}$
of $\Hv_k^K$ (cf. \eqref{eqn:defHhv}) and $V_k^K$ (cf. \eqref{eqn:defVhv}), respectively,
are called ``virtual" since they do not really known. Indeed, for example, we
do not know precisely if both sets are contained in the polynomial space. We
only know some conditions that satisfy their elements on $K$. 
More precisely, from the definition of $\Hv_k^K$,
it is quite clear that for every $\psiH_j$, its normal components, divergence and 
rotational are known. However, this is not entirely accurate. Indeed, what is really
known for each $\psiH_j$ are the values of their moments \eqref{eqn:dof-Hkv}, which
are given in \eqref{eqn:basisH}. In fact, these information is enough to determine,
through usual calculations, its normal components, divergence and rotational,
explicitly. We remark in advance that for the following computational implementation,
the identities \eqref{eqn:basisH} (resp. \eqref{eqn:basisV}) are the only ones required.

In order to perform the implementation of the matrices associated with the projectors,
as well as the future matrices associated with the terms arising from virtual schemes,
that we follow the methodology employed in \cite{helen-fila}, in which each matrix
is assemble through the previous construction of auxiliary matrices intrinsically related
to the mixed-VEM method defined by the subspaces $\Hv_k^K$ and $V_k^K$ (see \eqref{eqn:defHhv}
and \eqref{eqn:defVhv}, respectively).

\subsection{Preliminaries}

According to the previous discussion, we need to know polynomial functions that allow
us to compute the local operators in a clear way. Thus, consider $K\in\Th$ and $k\geq 0$.
Then, let $\phi_1^e,\phi_2^e,\ldots,\phi_{k+1}^e$ be the basis on the edge $e\in\partial K$
defined in \eqref{eqn:defBe}. That is, it follows that
$$\phi_i^e(x)\ :=\ \left(\frac{x-x_e}{h_e}\right)^{i-1}\,,\qquad\text{for } i = 1,2,\ldots, k+1\,,$$
where $x_e$ is the midpoint of $e$ and $h_e$ its length. In turn, consider $\{\vvm{\phi}_j^e\}_{j=1}^{2(k+1)}$
the vectorial version of the basis $\{\phi_j^e\}_{j=1}^{k+1}$ defined by:
$$\vvm{\phi}_1^e\ :=\ \mat{\begin{array}{c}
\phi_1^e\\ 0
\end{array}},\quad
\vvm{\phi}_2^e\ :=\ \mat{\begin{array}{c}
\phi_2^e\\ 0
\end{array}}\,,\quad\ldots\,,\quad 
\vvm{\phi}_{k+1}^e\ :=\ \mat{\begin{array}{c}
\phi_{k+1}^e\\ 0
\end{array}},$$
$$\vvm{\phi}_{k+2}^e := \mat{\begin{array}{c}
0\\ \phi_1^e
\end{array}},\;\;
\vvm{\phi}_{k+3}^e := \mat{\begin{array}{c}
0\\ \phi_2^e
\end{array}},\;\;\ldots\,,\;\;
\vvm{\phi}_{2(k+1)}^e := \mat{\begin{array}{c}
0\\ \phi_{k+1}^e
\end{array}}.$$
Now, let $\varphi_1^K,\varphi_2^K,\ldots,\varphi_m^K$, with $m:=\frac{(k+1)(k+2)}{2}$,
be the basis \eqref{eqn:defBk} given by:
\begin{equation}\label{def:varphi}
\varphi_i^K(x,y)\ :=\ \left(\frac{x-x_K}{h_K}\right)^{\alpha}\left(\frac{y-y_K}{h_K}\right)^{\beta}\,,\qquad\text{for } i = 1,2,\ldots, m\,,
\end{equation}
where $\alpha+\beta\in\{0,1,\ldots,k\}$, $(x_K,y_K)$ is the barycenter of $K$ and $h_K$
the diameter of $K$. In particular, we consider 
$$i\ :=\ \frac{(\alpha+\beta+1)(\alpha+\beta+2)}{2}-\beta\,,$$
where $\alpha$, $\beta$ are positive integers, and $\alpha+\beta \leq k$.
Using the previous ordering, we can guarantee the hierarchy of the basis $\{\varphi_j^K\}_{j=1}^{m}$,
which indicates that $\varphi_1^K$ is a constant polynomial.
Next, a basis of $\Pv_k(K)$ is given by:
$$\vvm{\varphi}_1^K\ :=\ \mat{\begin{array}{c}
\varphi_1^K\\ 0
\end{array}},\quad
\vvm{\varphi}_2^K\ :=\ \mat{\begin{array}{c}
\varphi_2^K\\ 0
\end{array}},\quad\ldots\,,\quad 
\vvm{\varphi}_m^K\ :=\ \mat{\begin{array}{c}
\varphi_m^K\\ 0
\end{array}},$$
$$\vvm{\varphi}_{m+1}^K\ :=\ \mat{\begin{array}{c}
0\\ \varphi_1^K
\end{array}},\quad
\vvm{\varphi}_{m+2}^K\ :=\ \mat{\begin{array}{c}
0\\ \varphi_2^K
\end{array}},\quad\ldots\,,\quad 
\vvm{\varphi}_{2m}^K\ :=\ \mat{\begin{array}{c}
0\\ \varphi_m^K
\end{array}},$$
whereas, we consider the following basis of $\Pm_k(K)$:
$$\Phi_1^K\ :=\ \mat{\begin{array}{cc}
\varphi_1^K & 0\\ 0 & 0
\end{array}},\quad
\Phi_2^K\ :=\ \mat{\begin{array}{cc}
\varphi_2^K & 0\\ 0 & 0
\end{array}},
\quad\ldots\,,\quad
\Phi_m^K\ :=\ \mat{\begin{array}{cc}
\varphi_m^K & 0\\ 0 & 0
\end{array}},$$
$$\Phi_{m+1}^K\ :=\ \mat{\begin{array}{cc}
0 & \varphi_1^K\\ 0 & 0
\end{array}},\quad
\Phi_{m+2}^K\ :=\ \mat{\begin{array}{cc}
0 & \varphi_2^K\\ 0 & 0
\end{array}},
\quad\ldots\,,\quad
\Phi_{2m}^K\ :=\ \mat{\begin{array}{cc}
0 & \varphi_m^K\\ 0 & 0
\end{array}},$$
$$\Phi_{2m+1}^K\ :=\ \mat{\begin{array}{cc}
0 & 0\\ \varphi_1^K & 0
\end{array}},\quad
\Phi_{2m+2}^K\ :=\ \mat{\begin{array}{cc}
0 & 0\\ \varphi_2^K & 0
\end{array}},\quad\ldots\,,\quad
\Phi_{3m}^K\ :=\ \mat{\begin{array}{cc}
0 & 0\\ \varphi_m^K & 0
\end{array}},$$
$$\Phi_{3m+1}^K\ :=\ \mat{\begin{array}{cc}
0 & 0\\ 0 & \varphi_1^K
\end{array}},\quad
\Phi_{3m+2}^K\ :=\ \mat{\begin{array}{cc}
0 & 0\\ 0 & \varphi_2^K
\end{array}},\quad\ldots\,,\quad
\Phi_{4m}^K\ :=\ \mat{\begin{array}{cc}
0 & 0\\ 0 & \varphi_m^K
\end{array}}.$$
It is important to clarify that the bases $\{\phi_i^e\}_{i=1}^{k+1}$ and $\{\varphi_i^K\}_{i=1}^m$
are related to the element $K$. That is, these must be constructed for each element $K\in\Th$. This
is a consequence of the fact that the geometry of every $K$ is not necessarily the same. Then, notice
that, from the hierarchical property of the basis $\{\varphi_i^K\}_{i=1}^m$, we can easily extend
this basis to degree $k+1$, which will be used throughout this section.

Furthermore, for each $e\in\partial K$, we also require the corresponding
Lagrange basis $\mathcal{L}_1^e,\mathcal{L}_2^e,\ldots,\mathcal{L}_{k+2}^e$
on $k+2$ uniformly spaced points of $e$. More precisely, for each $v\in V_k^K$,
using the notation $\alpha_{1,v}^e$, $\alpha_{2,v}^e$, $\ldots$, $\alpha_{k+2,v}^e$ by
the $k+2$ uniformly spaced points on $e$ with the corresponding orientation, we
have that $m^{V}_{a,v}(v) = \alpha_{1,v}^e$, $m^{V}_{e}(v) \in \{\alpha_{2,v}^e,\alpha_{3,v}^e,\ldots,\alpha_{k+1,v}^e\}$,
and $m^{V}_{b,v}(v) = \alpha_{k+2,v}^e$ (cf. \eqref{eqn:dof-Vkv}), where $a$
and $b$ are the vertices of $K$ that delimit $e$. Hence, it is well known that
there holds
\begin{equation}\label{eqn:interpLagrange}
v(x)\ =\ \sum_{i=1}^{k+2}\alpha_{i,v}^e\,\mathcal{L}_i^e(x)\,,\qquad\forall\ x\in e\,.
\end{equation}
In turn, we define
$$\vvm{\mathcal{L}}_1^e\ :=\ \mat{\begin{array}{c}
\mathcal{L}_1^e\\ 0
\end{array}},\quad
\vvm{\mathcal{L}}_2^e\ :=\ \mat{\begin{array}{c}
\mathcal{L}_2^e\\ 0
\end{array}}\,,\quad\ldots\,,\quad 
\vvm{\mathcal{L}}_{k+2}^e\ :=\ \mat{\begin{array}{c}
\mathcal{L}_{k+2}^e\\ 0
\end{array}},$$
$$\vvm{\mathcal{L}}_{k+3}^e := \mat{\begin{array}{c}
0\\ \mathcal{L}_1^e
\end{array}},\;\;
\vvm{\mathcal{L}}_{k+4}^e := \mat{\begin{array}{c}
0\\ \mathcal{L}_2^e
\end{array}},\;\;\ldots\,,\;\;
\vvm{\mathcal{L}}_{2(k+2)}^e := \mat{\begin{array}{c}
0\\ \mathcal{L}_{k+2}^e
\end{array}}.$$
At this point we remark in advance that we do not need to compute
the basis $\{\mathcal{L}_i^e\}_{i=1}^{k+2}$ explicitly, since we
only require the Lagrange basis defined on $k+2$ uniformly spaced points
of the interval $[0,1]$, which is denoted as:
\begin{equation}\label{def:lagrange}
\widehat{\mathcal{L}}_i(t)\ :=\ \prod_{\substack{{j=1}\\{j\neq i}}}^{k+2}\frac{t-t_j}{t_i-t_j}\,,\qquad\forall\; t\in[0,1]\,,\qquad\text{for } i = 1,2,\ldots, k+2\,,
\end{equation}
where $t_j := \frac{j-1}{k+1}$, with $j = 1,2,\ldots, k+2$,
when positive orientation is considered, and $t_j := 1-\frac{j-1}{k+1}$,
with $j = 1,2,\ldots, k+2$, otherwise.

On the other hand, we now introduce some matrices in order to facilitate
the construction of the discrete operator below. For simplicity, according
to the notation introduced at the end of Section \ref{sec:intro}, if we
write $\vv{A} := [a_{ij}]\in\R^{p\times q}$, then there hold that $i=1,2,\ldots,p$
and $j=1,2,\ldots,q$.
\begin{itemize}
\item The mass matrix on an edge $e$:
$$\Mk\ :=\ \left[\int_e\phi^e_i\phi^e_j\right]\in\R^{(k+1)\times(k+1)},$$
which, in order to compute this matrix, note that:
$$\int_e\phi^e_i\phi^e_j\ =\ h_e\int_0^1\widehat{\phi}_i(x)\widehat{\phi}_j(x)\,\text{d}x\,,$$
where $\widehat{\phi}_i(x) := \left(x-\frac{1}{2}\right)^{i-1}$, for $i=1,2,\ldots,k+1$.
Thus, the entries in $\Mk$ can be determined using a sufficiently precise
quadrature rule over $[0,1]$. More precisely, there holds:
$$\Mk\ :=\ h_e\left[\int_0^1\widehat{\phi}_i(x)\widehat{\phi}_j(x)\,\text{d}x\right]\in\R^{(k+1)\times(k+1)},$$
where the matrix on the right-hand side can be precomputed independent of
the edge $e$, which is important since we also require the matrix:
$$\iMk\ :=\ (\Mk)^{-1}\ =\ \frac{1}{h_e}\left[\int_0^1\widehat{\phi}_i(x)\widehat{\phi}_j(x)\,\text{d}x\right]^{-1}.$$

\item The following matrix:
\begin{equation*}\label{mat:Me}
\Me\ :=\ \left[\int_e\phi^e_i\varphi^K_j\right]\in\R^{(k+1)\times\m1},
\end{equation*}
with the basis $\{\varphi_i^K\}$ until degree $k+1$, and then
$\m1 := \dim\P_{k+1}(K)\frac{(k+2)(k+3)}{2}$. Next, notice that
$$\int_e\phi^e_i\varphi^K_j\ =\ h_e\int_0^1\widehat{\phi}_i(x)\,\varphi^K_j\big((1-x)\,\vv{v}_1 + x\,\vv{v}_2\big)\,\text{d}x\,,$$
where $\vv{v}_1$ and $\vv{v}_2$ are the vertices of $e$. Here, the ordering
of $\vv{v}_1$ and $\vv{v}_2$ is according to the orientation of $e$. For
example, we consider $\vv{v}_1$ the vertex of $e$ that has the lowest global
index in $\Th$. This selection guarantees to follow a unique orientation
on the same edge, independent of the element to which it belongs.

\item The mass matrix on an element $K$:
$$\MASS\ :=\ \left[\int_K\varphi_i^K\varphi_j^K\right]\in\R^{\m1\times\m1},$$
which, it is calculated for degree $k + 1$. In addition, we also require
the main submatrix $\Mass\in\R^{m\times m}$, whose index range is $[1,m]\times[1,m]$.
In other words, $\MASS$ is the mass matrix for the basis of $\P_{k+1}(K)$,
whereas $\Mass$ is the mass matrix for the basis of $\P_k(K)$.

Now, in order to compute the entries of $\MASS$, a quadrature rule over $K$
is not used. Indeed, employing the divergence theorem, we replace the area
integral to a sum of line integrals. More precisely, using Gauss's divergence
theorem, note that:
\begin{equation}\label{eqn:formulaIntegral}
\begin{array}{l}
\displaystyle\int_K\left(\frac{x-x_K}{h_K}\right)^{\alpha}\left(\frac{y-y_K}{h_K}\right)^{\beta}
\ =\ \frac{h_K}{\alpha+\beta+2}\int_K\div\left(\begin{array}{c}
\left(\frac{x-x_K}{h_K}\right)^{\alpha+1}\left(\frac{y-y_K}{h_K}\right)^{\beta}\\[1.5ex]
\left(\frac{x-x_K}{h_K}\right)^{\alpha}\left(\frac{y-y_K}{h_K}\right)^{\beta+1}
\end{array}\right)\\[4ex]
\quad =\ \displaystyle\sum_{e\in\partial K}\frac{h_K}{\alpha+\beta+2}\int_e\left(\frac{x-x_K}{h_K}\right)^{\alpha}\left(\frac{y-y_K}{h_K}\right)^{\beta}\bigg[\left(\frac{x-x_K}{h_K}\right)n_1^e\, +\, \left(\frac{y-y_K}{h_K}\right)n_2^e\bigg]\,,
 \end{array}
\end{equation}
where $\nv^e := (n_1^e,n_2^e)^{\tra}$ is the unit outward normal at $e$. Thus,
since all entries in $\MASS$ have the form $\int_K\big(\frac{x-x_K}{h_K}\big)^{\alpha}\big(\frac{y-y_K}{h_K}\big)^{\beta}$,
it is enough to use \eqref{eqn:formulaIntegral} to find them (see \eqref{def:varphi}).

\item Continuing with the basis $\{\varphi_i^K\}_{i=1}^{\m1}$
of $\P_{k + 1}(K)$, the gradient matrix is defined as:
\begin{equation}\label{mat:tM}
\begin{array}{rcl}
\tM & := & \displaystyle \left[\int_K\nabla\varphi_{i+1}^K\cdot\nabla\varphi_{j+1}^K\right]\\[2ex]
    & = & \displaystyle \left[\int_K\partial_x\varphi_{i+1}^K\partial_x\varphi_{j+1}^K + \partial_y\varphi_{i+1}^K\partial_y\varphi_{j+1}^K\right]\in\R^{(\m1-1)\times(\m1-1)}.
\end{array}
\end{equation}
Note that the basis $\{\varphi_i^K\}$ is extended again to degree
$k+1$, but eliminating the first constant element. Also, the entries of $\tM$ can
be calculated by using the formula \eqref{eqn:formulaIntegral} twice.
\end{itemize}

Next, for the following matrices, it is important to note that
$\{\nabla\varphi_{i+1}^K\}_{i=1}^{\m1-1}$ is a basis of $\nabla\P_{k+1}(K)$.
Thus, we now aim to obtain a basis $\{\vv{q}_i^K\}_{i=1}^{\mo}$ of $\Gk$,
with $\mo := \dim\P_{k-1}(K) = \frac{1}{2}k(k+1)$. To do that, it is required to find, for
$i=1,2,\ldots,\mo$, the constants $\{\alpha_j^{(i)}\}_{j=1}^{2m}$ such that:
$$\vv{q}_i^K\ :=\ \sum_{j=1}^{2m}\alpha_j^{(i)}\,\vvm{\varphi}_j^K\qan \int_K\nabla\varphi^K_{i+1}\cdot\vv{q}_j^K\ =\ 0\,,\quad\text{for } i = 1,2,\ldots,\m1-1\,.$$
Equivalently, it is required to solve local rectangular linear systems:
\begin{equation}\label{eqn:systemImplemtation}
\Mo\Ao\ =\ \vv{0}_{(\m1-1)\times \mo}\,,
\end{equation}
where $\vv{0}_{(\m1-1)\times \mo}$ is the zero matrix of $\R^{(\m1-1)\times \mo}$, and the
matrices $\Mo\in\R^{(\m1-1)\times(2m)}$ and $\Ao\in\R^{(2m)\times \mo}$
are defined by:
$$\Mo\ :=\ \left[\int_K\nabla\varphi_{i+1}^K\cdot\vvm{\varphi}_j^K\right]$$
and
$$\Ao\ :=\ \left[\alpha_j^{(\ell)}\right]\ =\ \mat{\begin{array}{cccc}
\alpha_1^{(1)} & \alpha_1^{(2)} & \cdots & \alpha_1^{(\mo)}\\
\alpha_2^{(1)} & \alpha_2^{(2)} & \cdots & \alpha_2^{(\mo)}\\
\vdots & \vdots & \ddots & \vdots\\
\alpha_{2m}^{(1)} & \alpha_{2m}^{(2)} & \cdots & \alpha_{2m}^{(\mo)}
\end{array}}\,,$$
for all $i=1,2,\ldots,\m1-1$,\, $j=1,2,\ldots,2m$ and $\ell=1,2,\ldots,\mo$.
Regarding $\Mo$, its construction is similar to $\tM$ (cf. \eqref{mat:tM}).
Indeed, it is sufficient to see that $\Mo$ can be assembled at the following
block level:
\begin{equation*}\label{mat:Mo}
\Mo\ =\ \left[\begin{array}{c:c}
\displaystyle\left[\int_K\partial_x\varphi_{i+1}^K\varphi_j^K\right] \; &\;
\displaystyle\left[\int_K\partial_y\varphi_{i+1}^K\varphi_j^K\right]
\end{array}\right],
\end{equation*}
where both blocks have size $(\m1-1)\times m$. On the other hand, to determine
$\Ao$ we need to solve the system \eqref{eqn:systemImplemtation}, where $QR$
decomposition is a good choice for that.

Alternatively, it is important to mention that the degrees of freedom given
by $m_{\vv{q},\rot}^H$ (cf. \eqref{eqn:dof-Hkv}) can also be defined employing,
instead of $\Gk$, the basis of any polynomial space $\widetilde{\Pv}_k(K)$ that
satisfies:
$$\Pv_k(K)\ =\ \nabla\P_{k+1}(K)\,\oplus\,\widetilde{\Pv}_k(K)\,.$$
In this way, it is possible to choose a space $\widetilde{\Pv}_k(K)$ such that
the $QR$ decomposition is not required, and then we can calculate these degrees
of freedom more efficiently. For more details, we recommend to check
\cite[Section 2.1]{serendipityVEM}, and particularly \cite[eq. (2.10)]{serendipityVEM}.

At this point we already introduced some matrices that allow us to construct
the operators associated with elements in $\Hv_k^K$ and $\Hm_k^K$. However,
we will consider discrete schemes that also require the operators associated
with elements in $V_k^K$ and $\vv{V}_k^K$. According to this, in similar way
as before, we now introduce the following matrices to use below.
\begin{itemize}
\item Let $\{e_1,e_2,\ldots,e_{d_K}\}$ be the edges of $K$. Then, the matrices
$\CLe{\ell}\in\R^{(k+2)\times n_k^V}$, with $\ell = 1,2,\ldots,d_K$, are defined
at block level as:
\begin{equation*}\label{mat:CLe}
\CLe{\ell}\ :=\ \left[\begin{array}{c:c:c:c}
\vv{C}_1^{(\ell)}\;\; & \;\vv{0}_{(k+2)\times (k(\ell-1))}\; & \;\;\vv{C}_2\;\; & \;\vv{0}_{(k+2)\times (k(d_K-\ell) + \mo)}
\end{array}\right],
\end{equation*}
where $n_k^V$ is defined in \eqref{eqn:n2}. In addition, matrix
$\vv{C}_1^{(\ell)}\in\R^{(k+2)\times d_K}$ is identical to the zero matrix
$\vv{0}_{(k+2)\times d_k}$, except that it has $1$ in the entries $(1,\ell)$
and $(k+2,\widehat{d}\,)$, where $\widehat{d} = \ell+1$ if $\ell \neq d_K$,
and $\widehat{d} = 1$ otherwise. Furthermore, matrix $\vv{C}_2\in\R^{(k+2)\times k}$
corresponds to the identity matrix $\Id_{k+2}$, whose first and last column
were removed.

\item For $\ell = 1,2,\ldots,d_K$, we introduce the following matrices:
\begin{equation*}\label{mat:NLe}
\begin{array}{rcl}
\NLe{\ell} & := & \displaystyle\left[\int_{e_{\ell}}\Big\{n_1^{e_{\ell}}\partial_x\varphi_{i+1}^K + n_2^{e_{\ell}}\partial_y\varphi_{i+1}^K\Big\}\,\mathcal{L}_j^{e_{\ell}}\right]\\[2ex]
& = & \displaystyle \left[h_{e_{\ell}}\int_0^1\Big\{n_1^{e_{\ell}}\partial_x\varphi^K_{i+1}(\vv{p}_{t,\ell}) + n_2^{e_{\ell}}\partial_y\varphi^K_{i+1}(\vv{p}_{t,\ell})\Big\}\,\widehat{\mathcal{L}}_j(t)\,\text{d}t\right]\in\R^{(\m1-1)\times(k+2)},
\end{array}
\end{equation*}
where $\nv^{e_{\ell}} := (n_1^{e_{\ell}},n_2^{e_{\ell}})^{\tra}$ is the unit
outward normal at $e_{\ell}$, and $\vv{p}_{t,\ell}:=(1-t)\,\vv{v}_1^{(\ell)} + t\,\vv{v}_2^{(\ell)}$.
Finally, $\vv{v}_1^{(\ell)}$ and $\vv{v}_2^{(\ell)}$ are the vertices of $e_{\ell}$, sorted
according to the orientation of $e_{\ell}$.	Moreover, we recall here that
$\{\widehat{\mathcal{L}}_j\}_{j=1}^{k+2}$ (cf. \eqref{def:lagrange}) depends of the orientation
of $e_{\ell}$. One again, the entries in $\NLe{\ell}$ can be determined using a
sufficiently precise quadrature rule over $[0,1]$.

\item The matrix associated with the Laplacian:
\begin{equation*}\label{mat:ML}
\begin{array}{rcl}
\ML & := & \displaystyle \left[\int_K\Delta\varphi_{i+1}^K\varphi_j^K\right]\\[2ex]
    & = & \displaystyle \left[\int_K\partial_{xx}\varphi_{i+1}^K\varphi_j^K + \partial_{yy}\varphi_{i+1}^K\varphi_j^K\right]\in\R^{(\m1-1)\times\mo}.
\end{array}
\end{equation*}
We remark here that the entries of $\ML$ can be calculated by using the formula
\eqref{eqn:formulaIntegral} twice.

\item The change of basis matrix on $e$:
$$\MCBe\ :=\ \left[\int_e\mathcal{L}_i^e\,\phi^e_j\right]\ =\ h_e\left[\int_0^1\widehat{\mathcal{L}}_i(t)\,\widehat{\phi}_j(t)\,\text{d}t\right]\in\R^{(k+2)\times(k+1)}\,,$$
where, it is important to note that the matrix on the right-hand side can be
precomputed independent of the edge $e$.

\item Now, we consider the mass matrix for the Lagrange basis:
$$\MassL\ :=\ \left[\int_0^1\widehat{\mathcal{L}}_i(t)\,\widehat{\mathcal{L}}_j(t)\,\text{d}t\right]\in\R^{(k+2)\times(k+2)},$$
which in fact has two possibilities, one for each possible orientation.
\end{itemize}

Finally, in order to define the foregoing local matrices, we first recalling from the
Section \ref{sec:dofH} (see also \eqref{eqn:basisH}) that $\{\psiH_j\}_{j=1}^{n_k^H}$
denotes the canonical basis of $\Hv_k^K$. Then, let $\{\PsiH_j\}_{j=1}^{2n_k^H}$
be the canonical basis of $\Hm_k^K$ given by
$$\PsiH_1\ :=\ \mat{\begin{array}{c}
\psiH_1\\[0.5ex] \vv{0}
\end{array}}\,,\quad
\PsiH_2\ :=\ \mat{\begin{array}{c}
\psiH_2\\[0.5ex] \vv{0}
\end{array}}\,,\quad\ldots\,,\quad \PsiH_{n_k^H}\ :=\ \mat{\begin{array}{c}
\psiH_{n_k^H}\\[0.5ex] \vv{0}
\end{array}}\,,$$
$$\PsiH_{n_k^H+1}\ :=\ \mat{\begin{array}{c}
\vv{0}\\ \psiH_1
\end{array}}\,,\quad \PsiH_{n_k^H+2}\ :=\ \mat{\begin{array}{c}
\vv{0}\\ \psiH_2
\end{array}}\,,\quad\ldots\,,\quad
\PsiH_{2n_k^H}\ :=\ \mat{\begin{array}{c}
\vv{0}\\ \psiH_{n_k^H}
\end{array}}.$$
Similarly, from Section \ref{sec:dofV} (see also \eqref{eqn:basisV}),
we use the basis $\{\psiV_j\}_{j=1}^{n_k^V}$ of $V_k^K$ to define
a basis $\{\PsiV_j\}_{j=1}^{2n_k^V}$ of $\vv{V}_k^K$ as
$$\PsiV_1\ :=\ \mat{\begin{array}{c}
\psiV_1\\[0.5ex] 0
\end{array}}\,,\quad
\PsiV_2\ :=\ \mat{\begin{array}{c}
\psiV_2\\[0.5ex] 0
\end{array}}\,,\quad\ldots\,,\quad \PsiV_{n_k^V}\ :=\ \mat{\begin{array}{c}
\psiV_{n_k^V}\\[0.5ex] 0
\end{array}}\,,$$
$$\PsiV_{n_k^V+1}\ :=\ \mat{\begin{array}{c}
0\\ \psiV_1
\end{array}}\,,\quad \PsiV_{n_k^V+2}\ :=\ \mat{\begin{array}{c}
0\\ \psiV_2
\end{array}}\,,\quad\ldots\,,\quad
\PsiV_{2n_k^V}\ :=\ \mat{\begin{array}{c}
0\\ \psiV_{n_k^V}
\end{array}}.$$

\subsection{The $L^2(K)$-orthogonal projection for elements of $\Hv_k^K$}\label{sec:Pb}

Proceeding similar to \cite[Section 5.2]{helen-fila}, we now aim to
describe the implementation of the $\Lm(K)$-orthogonal projection
$\bPk : \Lm(K)\rightarrow\Pm_k(K)$. More precisely, we introduce
the matrix $\Pb\in\R^{(4m)\times (2n_k^H)}$, which performs that
projection. In other words, $\Pb$ is a matrix that allows us to
approximate the elements of the tensor virtual basis $\{\PsiH_i\}_{i=1}^{2n_k^H}$
through the linear combinations of the elements of the tensor polynomial
basis $\{\Phi_i^K\}_{i=1}^{4m}$ by using $\Pb$. It is important
to recall that $n_k^H$ was defined in \eqref{eqn:n1}. More precisely,
given $\btau\in\Hm_k^K$ such that
$$\btau\ =\ \sum_{j=1}^{2n_k^H}\alpha_j\,\PsiH_j\qan \bPk(\btau)\ =\ \sum_{j=1}^{4m}\beta_j\,\Phi_j^K\,,$$
we seek a matrix $\Pb\in\R^{(4m)\times (2n_k^H)}$ such that
$\vvm{\beta} = \Pb\vvm{\alpha}$. In order to do that, we first
consider the following sequential list of matrices:
\begin{itemize}
\item $\ZZ\, :=\, (\Mo)^{\tt t}\,\tM^{-1}\in\R^{(2m)\times(\m1-1)}$,
with $m := \frac{1}{2}(k+1)(k+2)$ and $\m1 := \frac{1}{2}(k+2)(k+3)$.

\item $\Be\, :=\, s^e\,(\Me)^{\tra}\,\iMk\in\R^{\m1\times(k+1)}$, for each $e\in\partial K$,
where if $I_{v_1}$ and $I_{v_2}$ are the global indexes of the
vertices of $e$, it follows:
\begin{equation}\label{eqn:signo}
s^e\ :=\ \left\{\begin{array}{rl}
 1 & \text{if } I_{v_1} < I_{v_2},\\[0.5ex]
-1 & \text{otherwise},
\end{array}\right.
\end{equation}
represents an orientation for $e$, which allows the integrals involved
in $\Me$ to have the same value independent of the element $K$. Moreover,
the matrix $\Be^{(r,s)}\in\R^{(s-r+1)\times(k+1)}$ matrix is defined as
the submatrix of $\Be$ that contains all of its columns, but only the rows
from $r$ to $s$, with $1\leq r < s \leq \m1$.

\item The divergence matrix:
\begin{equation}\label{mat:Di}
\Di\ :=\ \Mass^{-1}\,\widetilde{\vv{B}}\in\R^{m\times n_k^H},
\end{equation}
with $n_k^H$ defined in \eqref{eqn:n1} and $\widetilde{\vv{B}}\in\R^{m\times n_k^H}$
is defined at block level as:
$$\widetilde{\vv{B}}\ :=\ \left[\begin{array}{c:c:c:c:c:c}
\vv{B}_{e_1}^{(1,m)} & \vv{B}_{e_2}^{(1,m)} & \cdots & \vv{B}_{e_{d_k}}^{(1,m)} & \;-\widetilde{\vv{C}}\;\; & \;\;\vv{0}_{m\times\mo}
\end{array}\right]\,,$$
where $\mo := \frac{1}{2}k(k+1)$. Furthermore, the matrix $\widetilde{\vv{C}}\in\R^{m\times(m-1)}$
is identical to the identity matrix $\Id_m$, except that its first column
was removed.

\item The matrix $\Pr\in\R^{(2m)\times n_k^H}$, which is associated to the
right-hand side of \eqref{eqn:defL2Proj-v}. More precisely, using \eqref{eqn:proyP}
we define $\Pr$ as follows:
\begin{eqnarray*}
\Pr & := & \left[\int_K\vvm{\varphi}_i^K\cdot\psiH_j\right]\ =\ \left[-\int_K\widetilde{\varphi_i^K}\,\div(\psiH_j)\right]\, +\,
           \left[\int_{\partial K}\widetilde{\varphi_i^K}\,(\psiH_j\cdot \nv)\right]\, +\,
           \left[\int_K(\vvm{\varphi}_i^K)^{\perp}\cdot\psiH_j\right] \\[1ex]
     & = & -\ZZ\,\widetilde{\vv{M}}_{\text{mass}}\,\Di\, +\, \ZZ\left[\begin{array}{c:c:c:c:c}
\vv{B}_{e_1}^{(2,\m1)} & \vv{B}_{e_2}^{(2,\m1)} & \cdots & \vv{B}_{e_{d_k}}^{(2,\m1)} & \;\;\vv{0}_{(\m1-1)\times(m-1+\mo)}
\end{array}\right]\\
    &    & +\ \Big[\begin{array}{c:c}
\vv{0}_{(2m)\times((k+1)d_K + m-1)}\; & \; (\IMASS - \ZZ\,\Mo)\,\Ao\,(\Ao^{\tra}\,\IMASS\,\Ao)^{-1}
\end{array}\Big]\,,
\end{eqnarray*}
where $\widetilde{\vv{M}}_{\text{mass}}\in\R^{(\m1-1)\times m}$ is the
submatrix of $\MASS$, which considers the rows from $2$ to $\m1$ and columns
from $1$ to $m$. In turn, it follows that:
\begin{equation*}\label{mat:IMASS}
\IMASS\ :=\ \kron(\Id_2,\, \Mass)\in\R^{(2m)\times(2m)}\,,
\end{equation*}
where $\kron(\cdot,\cdot)$ corresponds to the usual Kronecker product.
\end{itemize}
Finally, using the foregoing matrices, we assemble the following matrix:
\begin{equation}\label{mat:Pb-0}
\Pbt\ :=\ \kron(\Id_2,\, \Mass^{-1})\,\Pr\in\R^{(2m)\times n_k^H}\,,
\end{equation}
which corresponds to the coefficients of the elements of $\{\psiH_i\}_{i=1}^{n_k^H}$
under the operator $\Pk:\Lv(K)\to\Pv_k(K)$. On the other hand, the matrix of coefficients
of the elements of $\{\PsiH_i\}_{i=1}^{2n_k^H}$ under the operator $\bPk:\Lm(K)\to\Pm_k(K)$
is given by:
$$\Pb\ :=\ \kron(\Id_2,\, \Pbt)\in\R^{(4m)\times (2n_k^H)}\,.$$

\subsection{The projection operator $R_k^K$ for elements of $V_k^K$}\label{sec:Rb}

Concerning the projection $R_k^K : H^1(K)\to\P_{k+1}(K)$ (cf. \eqref{eqn:defRh-0}),
we now aim to define the matrix $\Rb$ associated to the coefficients of the
elements of $\{\psiV_i\}_{i=1}^{n_k^V}$ under the operator $R_k^K$. Firstly,
we define the matrices $\BR_1,\BR_2\in\R^{(\m1-1)\times n_k^V}$ as follow:
$$\BR_1\ :=\ \displaystyle\left[-\int_K\Delta\varphi_{i+1}^K\,\psiV_j\right]\ =\ \displaystyle -\ML\,\big(\vv{M}_{\text{mass},K}^{(k-1)}\big)^{-1}\,\left[\begin{array}{c:c}\vv{0}_{\mo\times((k+1)d_K)}\; & \;\;\Id_{\mo}
\end{array}\right],$$
and
$$\BR_2\ :=\ \displaystyle\left[\int_{\partial K}(\nabla\varphi_{i+1}^K\cdot\nv)\,\psiV_j\right]\ =\ \displaystyle\sum_{\ell=1}^{d_K}\NLe{\ell}\,\CLe{\ell},$$
where we use \eqref{eqn:interpLagrange} with $v = \psiV_j$. In addition,
$\vv{M}_{\text{mass},K}^{(k-1)}\in\R^{\mo\times\mo}$ is the submatrix
of $\MASS$, which considers the rows and columns from $1$ to $\mo$.

Now, from \eqref{eqn:proyR} we easily realize that the first equation
in \eqref{eqn:defRh-0} is associated to the matrix
\begin{equation}\label{mat:Rb-0}
\tR\ :=\ \tM^{-1}\,(\BR_1 + \BR_2)\in\R^{(\m1-1)\times n_k^V}\,,
\end{equation}
whereas the second one is implemented as:
$$\vvm{\eta}_1\ :=\ \left\{\begin{array}{rl}
\displaystyle\bigg(\sum\limits_{e\in\partial K}h_e\bigg)^{-1}\sum_{e\in\partial K}\left\{h_e\,\vv{c}_0\,\CLe{}\, -\, \vv{c}_1\tR\right\} & \text{if } k=0,\\[4ex]
\displaystyle\frac{1}{|K|}\left\{\left[\begin{array}{c:c:c}\vv{0}_{1\times((k+1)d_K)} & 1 & \vv{0}_{1\times(\mo-1)}\end{array}\right]\,
-\, \vv{c}_2\,\tR\right\} & \text{otherwise},
\end{array}\right.$$
where $\vv{c}_1,\vv{c}_2\in\R^{1\times \m1-1}$ are the first row of $\Me\in\R^{(k+1)\times\m1}$
and $\MASS\in\R^{\m1\times\m1}$, respectively, which their first entry has been removed.
Moreover, the vector $\vv{c}_0 := \left[\int_0^1\widehat{\mathcal{L}}_j(t)\,\text{d}t\right]\in\R^{1\times(k+2)}$
can be precomputed depending only of the orientation of $e$.

We end this section by introducing the matrix of coefficients of the
elements of $\{\psiV_i\}_{i=1}^{n_k^V}$ under the operator $R_k^K : H^1(K)\to\P_{k+1}(K)$
is given by:
$$\Rb\ :=\ \left[\begin{array}{c}
\vvm{\eta}_1\\
\\[-2ex] \hdashline
\\[-2ex]
\tR
\end{array}\right]\in\R^{\m1\times n_k^V}\,.$$
It is important to remark here that the matrix of coefficients of the
elements of $\{\PsiV_i\}_{i=1}^{2n_k^V}$ under the operator $\Rk : \Hv^1(K)\to\Pv_{k+1}(K)$
can be obtained by $\kron(\Id_2,\,\Rb)\in\R^{(2\m1)\times(2n_k^V)}$.

\subsection{The $L^2(K)$-orthogonal projection for elements of $V_k^K$}\label{sec:PL2u}

Similarly to the previous section, we now aim to introduce the matrices
associated with the projection $P_k^K:L^2(K)\to\P_k(K)$ for elements of
the canonical basis of $V_k^K$, and the projection $\bPk:\Lm(K)\to\Pm_k(K)$
for the gradient of elements of the canonical basis of $\vv{V}_k^K$. Indeed,
we begin by considering the matrix of coefficients of the elements of
$\{\psiV_i\}_{i=1}^{n_k^V}$ under the operator $P_k^K$, which employing
\eqref{eqn:proyPU} is given by:
\begin{eqnarray}
\PU & := & \left[\int_K\varphi_i^K\,\varphi_j^K\right]^{-1}\left[\int_K\varphi_i^K\,\psiV_j\right]\nonumber\\[1ex]
& = & \Mass^{-1}\,\left[\begin{array}{c}
\left[\begin{array}{c:c}\vv{0}_{\mo\times((k+1)d_K)}\; & \; \Id_{\mo}\end{array}\right]\\
\\[-2ex] \hdashline
\\[-2ex]
\widetilde{\vv{U}}\,\Rb
\end{array}\right]\in\R^{m\times n_k^V}\,,\label{eqn:matPV}
\end{eqnarray}
where the matrix $\widetilde{\vv{U}}\in\R^{(m-\mo)\times\m1}$ is the submatrix
of $\MASS$, whose index range is $[\mo+1,m]\times[1,\m1]$. Now, we follow \eqref{eqn:proyPDU}
and let $\PGU\in\R^{(4m)\times(2n_k^V)}$ be the matrix of coefficients of the elements of
$\{\nabla\PsiV_i\}_{i=1}^{2n_k^V}$ under the operator $\bPk$. Hence, we deduce that:
\begin{eqnarray*}
\PGU\!\!\! & := & \left[\int_K\Phi_i^K:\Phi_j^K\right]^{-1}\left[\int_K\Phi_i^K:\nabla\PsiV_j\right]\\[1ex]
     & = & \kron\left(\Id_4,\, \Mass^{-1}\right)\kron\left(\Id_2,\, \left[\int_K\vvm{\varphi}_i^K\cdot\nabla\psiV_j\right]\right)\\[1ex]
     & = & \kron\left(\Id_4,\, \Mass^{-1}\right)\kron\left(\Id_2,\, {-}\left[\int_K\div(\vvm{\varphi}_i^K)\,P_k^K(\psiV_j)\right]\, +\, \left[\int_{\partial K}(\vvm{\varphi}_i^K\cdot\nv)\,\psiV_j\right]\right)\\[1ex]
     & = & \kron\left(\Id_4,\, \Mass^{-1}\right)\kron\left(\Id_2,\, {-}\left[\int_K\div(\vvm{\varphi}_i^K)\,\varphi_j^K\right]\PU\, +\! \sum_{e\in\partial K}\kron\left(\nv^e,\, \left[\int_e\varphi_i^K\,\mathcal{L}_j^e\right]\!\CLe{}\right)\right)\!,
\end{eqnarray*}
where $\nv^e := (n_1^e,n_2^e)^{\tra}$ is the unit outward normal at $e$.
In addition, we introduce the following matrices:
$$\DxDy\ :=\ \left[\int_K\div(\vvm{\varphi}_i^K)\,\varphi_j^K\right]\ =\ \left[\begin{array}{c}
\vv{0}_{1\times m}\\[-2ex]
\\ \hdashline
\\[-2ex]
\vv{D}_x\\[-2ex]
\\ \hdashline
\\[-2ex]
\vv{0}_{1\times m}\\[-2ex]
\\ \hdashline
\\[-2ex]
\vv{D}_y
\end{array}\right]\in\R^{(2m)\times m}\,,$$
and
$$\PLe\ :=\ \left[\int_e\varphi_i^K\,\mathcal{L}_j^e\right]\ =\ \left[h_e\int_0^1\varphi_i^K((1-t)\,\vv{v}_1 + t\,\vv{v}_2)\,\widehat{\mathcal{L}}_j(t)\,\text{d}t\right]\in\R^{m\times(k+2)}\,,$$
where $\vv{D}_x\in\R^{(m-1)\times m}$ is the submatrix of $\Mo$ with
index range $[1,m-1]\times[1,m]$, whereas $\vv{D}_y\in\R^{(m-1)\times m}$
is the submatrix of $\Mo$ with index range $[1,m-1]\times[m+1,2m]$.
Moreover, once again, $\vv{v}_1$ and $\vv{v}_2$ are the oriented vertices
of $e$. Then, according to the previous analysis, we conclude that
$$\PGU\ =\ \kron\left(\Id_4,\, \Mass^{-1}\right)\kron\left(\Id_2,\, {-}\DxDy\,\PU\, +\, \sum_{e\in\partial K}\kron\left(\nv^e,\, \PLe\,\CLe{}\right)\right)\in\R^{(4m)\times(2n_k^V)}\,.$$

\section{Local discrete operators arising from VEM schemes}\label{sec:operators}

We begin by describing in detail a way to assemble diverse local terms
arising in mixed-primal virtual element formulations as in \cite{cgs-1,cgs-2,cgs-3,gms-vem-ns,gms-brinkman,gms-boussinesq}.
In particular, we are interested in those coming from the scheme proposed in \cite{gms-vem-ns}
for the Navier-Stokes system, which is recalling in Section \ref{sec:model} below.
For simplicity, we define three categories to separate the terms (or operators)
related to the two virtual subspaces $\Hm_k^K$ and $\vv{V}_k^K$, as well as those
that combine both. Therein, matrices and vectors are described locally to eventually
be assembled to the respective global linear system. A fourth category related to
particular nonlinear schemes (such as Navier-Stokes or Boussinesq) is presented
later in the Section \ref{sec:opNoLineal}.

\begin{figure}[h!t]\centering
	\scalebox{0.69}{\includegraphics{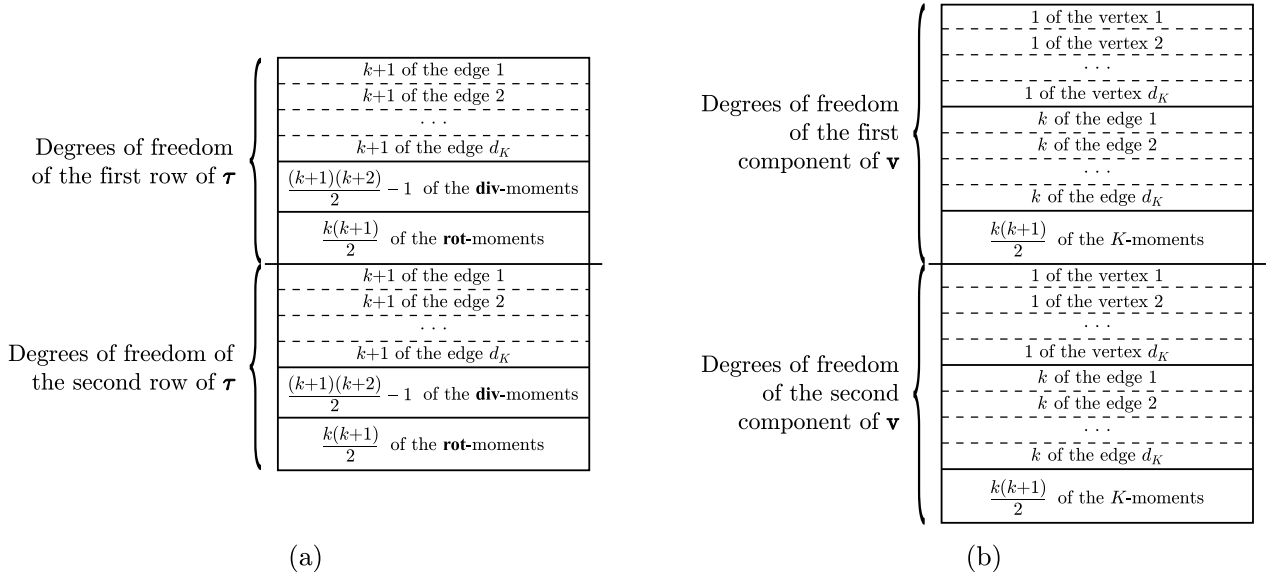}}
	\vspace{-2ex}
	\caption{Ordering of local degrees of freedom for each $K\in\Th$.}\label{fig:dof}
	\end{figure}

On the other hand, for each $K\in\Th$, in what follows we sort the degrees of freedom
of each element $\btau\in\Hm_k^K$ as in Figure \ref{fig:dof}(a), whereas the moments of
each element $\vv{v}\in\vv{V}_k^K$ follow the ordering described in Figure \ref{fig:dof}(b).

\subsection{Operators related to the elements of $\Hm_k^K$}

Given a volume force $\fv\in\Lv(\Omega)$ and a Dirichlet datum $\bg\in\Hv^{1/2}(\Gamma)$,
we consider the following discrete operators:
	\begin{itemize}
	\item $\displaystyle\AaK\ :=\ \left[\int_K\big[\bPk(\PsiH_i)\big]^{\td}:\big[\bPk(\PsiH_j)\big]^{\td}\right]\in\R^{2n_k^H\times 2n_k^H}$
	\item $\displaystyle\AdK\ :=\ \left[\int_K\tr\big(\bPk(\PsiH_i)\big)\,\tr\big(\bPk(\PsiH_j)\big)\right]\in\R^{2n_k^H\times 2n_k^H}$
	\item $\displaystyle\AcK\ :=\ \left[\int_K\bdiv(\PsiH_i)\cdot\bdiv(\PsiH_j)\right]\in\R^{2n_k^H\times 2n_k^H}$
	\item $\displaystyle\AbK\ :=\ \left[\mathcal{S}_H^K\big(\PsiH_i - \bPk(\PsiH_i), \PsiH_j - \bPk(\PsiH_j)\big)\right]\in\R^{2n_k^H\times 2n_k^H}$
	\item $\displaystyle\CK\ :=\ \left[\int_K\tr(\PsiH_i)\right]\in\R^{2n_k^H\times 1}$
	\item $\displaystyle\vv{b}_1^K\ :=\ \left[\int_{\partial K\cap\Gamma}\PsiH_i\nv\cdot\vv{g}\right]\in\R^{2n_k^H\times 1}$
	\item $\displaystyle\vv{b}_2^K\ :=\ \left[\int_K\bdiv(\PsiH_i)\cdot\vv{f}\right]\in\R^{2n_k^H\times 1}$
	\end{itemize}
where the bilinear form $\mathcal{S}_H^K(\cdot,\cdot)$ is defined in \eqref{eqn:defSH}.
Now, the operator $\AaK\in\R^{2n_k^H\times 2n_k^H}$ is defined by:
\begin{eqnarray*}
\AaK & := & \left[\int_K\big[\bPk(\PsiH_i)\big]^{\td}:\big[\bPk(\PsiH_j)\big]^{\td}\right]\ =\ (\Pb)^{\tra}\,\left[\int_K(\Phi_i^K)^{\td}:(\Phi_j^K)^{\td}\right]\,\Pb\\[1ex]
& = & (\Pb)^{\tra}\,\kron(\vv{M}_{\text{dev}},\, \Mass)\,\Pb,
\end{eqnarray*}
where
\begin{equation}\label{mat:Mdev}
\vv{M}_{\text{dev}}\ :=\ \mat{\begin{array}{rrrr}
\frac{1}{2} & 0 & 0 & -\frac{1}{2}\\ 0 & 1 & 0 & 0\\ 0 & 0 & 1 & 0\\ -\frac{1}{2} & 0 & 0 & \frac{1}{2}
\end{array}},
\end{equation}
is the matrix associated to the deviator operator.
Similarly, for the operator $\AdK\in\R^{2n_k^H\times 2n_k^H}$,
it follows:
\begin{eqnarray*}
\AdK & := & \left[\int_K\tr\big(\bPk(\PsiH_i)\big)\,\tr\big(\bPk(\PsiH_j)\big)\right]\ =\ (\Pb)^{\tra}\,\left[\int_K\tr(\Phi_i^K)\,\tr(\Phi_j^K)\right]\,\Pb\\[1ex]
& = & (\Pb)^{\tra}\,\kron\left(\mat{\begin{smallmatrix}1 & \;0\; & 0 & \;1\\[0.5ex]0 & 0 & 0 & \;0\\[0.5ex]0 & 0 & 0 & \;0\\[0.5ex]1 & 0 & 0 & \;1\end{smallmatrix}},\, \Mass\right)\,\Pb,
\end{eqnarray*}
whereas, the operator $\AcK\in\R^{2n_k^H\times 2n_k^H}$ can be computed as:
\begin{eqnarray*}
\AcK & := & \left[\int_K\bdiv(\PsiH_i)\cdot\bdiv(\PsiH_j)\right]\ =\ \kron\left(\Id_2,\, \left[\int_K\div(\psiH_i)\,\div(\psiH_j)\right]\right)\\[1ex]
& = & \kron\left(\Id_2,\, (\Di)^{\tra}\,\left[\int_K\varphi_i^K\varphi_j^K\right]\,\Di\right)\ =\ \kron\left(\Id_2,\, (\Di)^{\tra}\,\Mass\,\Di\right).
\end{eqnarray*}

On the other hand, the operator $\AbK\in\R^{2n_k^H\times 2n_k^H}$
is assemble in a different way. More precisely, from \eqref{eqn:defSH},
we have that
$$\AbK\ :=\ \left[\mathcal{S}_H^K\big(\PsiH_i - \bPk(\PsiH_i), \PsiH_j - \bPk(\PsiH_j)\big)\right]\ =\ \kron(\Id_2,\, \Ha^{\tra}\;\Ha)\,,$$
where $\Ha\in\R^{n_k^H\times n_k^H}$ can be decomposed as: $\Ha = \vv{H}_1 - \vv{H}_2\,\Pbt$.
Here, $\Pbt\in\R^{(2m)\times n_k^H}$ is defined in \eqref{mat:Pb-0},
whereas $\vv{H}_1\in\R^{n_k^H\times n_k^H}$ is a diagonal matrix given
at block level by:
\begin{equation*}\label{mat:H1}
\vv{H}_1\ :=\ \left[\begin{array}{c:c:c:c:c}
s^{e_1}\,\Id_{k+1} & & & & \\
& & & & \\[-2ex] \hdashline
& & & & \\[-2ex]
& s^{e_2}\,\Id_{k+1} & & & \\
& & & & \\[-2ex] \hdashline
& & & & \\[-2ex]
& &\ddots & & \\
& & & & \\[-2ex] \hdashline
& & & & \\[-2ex]
& & & s^{e_{d_K}}\,\Id_{k+1} & \\
& & & & \\[-2ex] \hdashline
& & & & \\[-2ex]
& & & & \vv{0}_{m-1+\mo}\\
& & & & \\[-2ex]
\end{array}\right],
\end{equation*}
where the constant $s^e$ is defined in \eqref{eqn:signo}, and $\mo := \frac{1}{2}k(k+1)$.
Furthermore, for each $e\in\partial K$, let $\Ee{e}\in\R^{(k+1)\times(2m)}$ be the matrix given by:
$$\Ee{e}\ :=\ \kron(\nv^e,\, \widetilde{\vv{M}}_{e,K})\,,$$
where $\nv^e\in\R^{1\times 2}$ is the unit outward row normal at $e$
and $\widetilde{\vv{M}}_{e,K}\in\R^{(k+1)\times m}$ is the submatrix
of $\Me$, which contains all its rows, but only the columns from $1$ to $m$.
In this way, it follows that $\vv{H}_2\in\R^{n_k^H\times(2m)}$ is defined as:
\begin{equation*}\label{mat:H2}
\vv{H}_2\ :=\ \left[\begin{array}{c}
\Ee{e_1}\\
\\[-2ex] \hdashline
\\[-2ex]
\Ee{e_2}\\
\\[-2ex] \hdashline
\\[-2ex]
\vdots\\
\\[-2ex] \hdashline
\\[-2ex]
\Ee{e_{d_K}}\\
\\[-2ex] \hdashline
\\[-2ex]
\vv{0}_{(m-1+\mo)\times(2m)}
\end{array}\right].
\end{equation*}

Furthermore, the operator $\CK\in\R^{2n_k^H\times 1}$ satisfies that:
$$\CK\ :=\ \left[\int_K\tr(\PsiH_i)\right]\ =\ (\Pb)^{\tra}\left[\int_K\tr(\Phi_i^K)\right]\ =\ (\Pb)^{\tra}\,\kron\left(\mat{\begin{smallmatrix}1\\[0.5ex]0\\[0.5ex]0\\[0.5ex]1\end{smallmatrix}},\, \vv{c}_3\right),$$
where the vector $\vv{c}_3\in\R^{m\times 1}$ is the first column of $\Mass$.

Finally, letting $\vv{g} := (g_1, g_2)^{\tra}$ be the Dirichlet datum,
and given a boundary edge $e$, we define the following vectors for the
coefficients of the $L^2(e)$-orthogonal projection:
\begin{equation}\label{eqn:pg}
\begin{array}{rcl}
\Pg{\ell} & := & \displaystyle\left[\int_e\phi_i^e\,\phi_j^e\right]^{-1}\left[\int_e\phi_i^e\,g_{\ell}\right]\\[3ex]
& = & \displaystyle\iMk\left[h_e\int_0^1\widehat{\phi}_i(x)\,g_{\ell}((1-x)\,\vv{v}_1 + x\,\vv{v}_2)\,\text{d}x\right]\in\R^{(k+1)\times 1}\,,
\end{array}
\end{equation}
for $\ell=1,2$, where $\vv{v}_1$ and $\vv{v}_2$ are the oriented vertices of $e$.
Next, note that the operator $\vv{b}_1^K\in\R^{2n_k^H\times 1}$ is given by:
\begin{eqnarray*}
\vv{b}_1^K & := & \left[\int_{\partial K\cap\Gamma}\PsiH_i\nv\cdot\vv{g}\right]\ =\ \sum_{e\in\partial K\cap\Gamma}\left[\begin{array}{c}
\displaystyle\left[\int_e(\psiH_i\cdot\nv^e)\,g_1\right]\\[-2ex]
\\ \hdashline
\\[-2ex]
\displaystyle\left[\int_e(\psiH_i\cdot\nv^e)\,g_2\right]
\end{array}\right]\\[1ex]
& = & \sum_{e\in\partial K\cap\Gamma}\left[\begin{array}{c}
\displaystyle\left[\int_e(\psiH_i\cdot\nv^e)\,\phi_j^e\right]\Pg{1}\\[-2ex]
\\ \hdashline
\\[-2ex]
\displaystyle\left[\int_e(\psiH_i\cdot\nv^e)\,\phi_j^e\right]\Pg{2}
\end{array}\right]\ =\ \left[\begin{array}{c}
\vv{b}_{1,1}^{e_1}\\[-2ex]
\\ \hdashline
\vdots\\[-2ex]
\\ \hdashline
\\[-2ex]
\vv{b}_{1,1}^{e_{d_K}}\\[-2ex]
\\ \hdashline
\\[-2ex]
\vv{0}_{(m-1+\mo)\times 1}\\[-2ex]
\\ \hdashline\hdashline
\\[-2ex]
\vv{b}_{1,2}^{e_1}\\[-2ex]
\\ \hdashline
\vdots\\[-2ex]
\\ \hdashline
\\[-2ex]
\vv{b}_{1,2}^{e_{d_K}}\\[-2ex]
\\ \hdashline
\\[-2ex]
\vv{0}_{(m-1+\mo)\times 1}
\end{array}\right],
\end{eqnarray*}
where $\vv{b}_{1,\ell}^e\in\R^{(k+1)\times 1}$, for each $e\in\partial K$,
is defined as:
$$\vv{b}_{1,\ell}^e\ :=\ \left\{\begin{array}{cl}
s^e\,\Pg{\ell} & \text{if } e\in\Gamma,\\[1ex]
\vv{0}_{(k+1)\times 1} & \text{otherwise},
\end{array}\right.$$
for $\ell=1,2$, where $s^e$ is defined in \eqref{eqn:signo}.
Now, in a similar way, we consider $\vv{f} := (f_1, f_2)^{\tra}$ the source term,
and then introduce the following vectors for the coefficients of the $L^2(K)$-orthogonal
projection:
\begin{equation}\label{eqn:pf}
\Pf{\ell}\ :=\ \left[\int_K\varphi_i^K\,\varphi_j^K\right]^{-1}\left[\int_K\varphi_i^K\,f_{\ell}\right]\ =\
  \Mass^{-1}\left[\int_K\varphi_i^K\,f_{\ell}\right]\in\R^{m\times 1}\,,
\end{equation}
where the integral $\int_K\varphi_i^K\,f_{\ell}$ requires a suitable quadrature for
polygonal domains. Some examples can be found in \cite{polyquad-1,polyquad-2,polyquad-3}.
Hence, it follows that
\begin{eqnarray*}
\vv{b}_2^K & := & \left[\int_K\bdiv(\PsiH_i)\cdot\vv{f}\right]\ =\ \left[\begin{array}{c}
\displaystyle\left[\int_K\div(\psiH_i)\,P_k^K(f_1)\right]\\
\\[-2ex] \hdashline
\\[-2ex]
\displaystyle\left[\int_K\div(\psiH_i)\,P_k^K(f_2)\right]
\end{array}\right]\\[1ex]
& = & \left[\begin{array}{c}
\displaystyle (\Di)^{\tra}\left[\int_K\varphi_i^K\,\varphi_j^K\right]\Pf{1}\\
\\[-2ex] \hdashline
\\[-2ex]
\displaystyle (\Di)^{\tra}\left[\int_K\varphi_i^K\,\varphi_j^K\right]\Pf{2}
\end{array}\right]\ =\ \left[\begin{array}{c}
\displaystyle (\Di)^{\tra}\,\Mass\,\Pf{1}\\
\\[-2ex] \hdashline
\\[-2ex]
\displaystyle (\Di)^{\tra}\,\Mass\,\Pf{2}
\end{array}\right]\\[1ex]
& = & \left[\begin{array}{c}
\displaystyle (\Di)^{\tra}\left[\int_K\varphi_i^K\,f_1\right]\\
\\[-2ex] \hdashline
\\[-2ex]
\displaystyle (\Di)^{\tra}\left[\int_K\varphi_i^K\,f_2\right]
\end{array}\right]\in\R^{(2n_k^H)\times 1}\,,
\end{eqnarray*}
where $\Di\in\R^{m\times n_k^H}$ is defined in \eqref{mat:Di}.

\subsection{Operators related to the elements of $\vv{V}_k^K$}

Considering again $\fv\in\Lv(\Omega)$ and $\bg\in\Hv^{1/2}(\Gamma)$,
together with the bilinear form $\mathcal{S}_V^K(\cdot,\cdot)$ (see
\eqref{eqn:defSV}), we now consider the operators:
	\begin{itemize}
	\item $\displaystyle\DaK\ :=\ \left[\int_K\nabla\Rk(\PsiV_i):\nabla\Rk(\PsiV_j)\right]\in\R^{2n_k^V\times 2n_k^V}$
	\item $\displaystyle\DbK\ :=\ \left[\mathcal{S}_V^K\big(\PsiV_i - \Rk(\PsiV_i), \PsiV_j - \Rk(\PsiV_j)\big)\right]\in\R^{2n_k^V\times 2n_k^V}$
	\item $\displaystyle\DcK\ :=\ \left[\int_{\partial K\cap\Gamma}\PsiV_i\cdot\PsiV_j\right]\in\R^{2n_k^V\times 2n_k^V}$
	\item $\displaystyle\vv{b}_3^K\ :=\ \left[\int_{\partial K\cap\Gamma}\PsiV_i\cdot\vv{g}\right]\in\R^{2n_k^V\times 1}$
	\item $\displaystyle\vv{b}_4^K\ :=\ \left[\int_K\PsiV_i\cdot\Pk(\vv{f})\right]\in\R^{2n_k^V\times 1}$
	\end{itemize}
In turn, we begin by describing the operator $\DaK\in\R^{2n_k^V\times 2n_k^V}$
as follows:
\begin{eqnarray*}
\DaK & := & \left[\int_K\nabla\Rk(\PsiV_i):\nabla\Rk(\PsiV_j)\right]\ =\ \kron\left(\Id_2,\, \left[\int_K\nabla R_k^K(\psiV_i)\cdot\nabla R_k^K(\psiV_j)\right]\right)\\[1ex]
 & = & \kron\left(\Id_2,\, (\Rb)^{\tra}\,\left[\int_K\nabla\varphi_i^K\cdot\nabla\varphi_j^K\right]\,\Rb\right)\\[1ex]
 & = & \kron\left(\Id_2,\, (\Rb)^{\tra}\,\left[\begin{array}{c:c}
0 & \vv{0}_{1\times(\m1-1)}\\[-2ex]
& \\ \hdashline
& \\[-2ex]
 \vv{0}_{(\m1-1)\times 1}& \tM\end{array}\right]\,\Rb\right)\\[1ex]
 & = & \kron\left(\Id_2,\, (\tR)^{\tra}\,\tM\,\tR\right)\,,
\end{eqnarray*}
which $\tR$ is defined in \eqref{mat:Rb-0}.
Next, proceeding similar to $\AbK$, we obtain that
$\DbK\in\R^{2n_k^V\times 2n_k^V}$ can be written as:
$$\DbK\ :=\ \left[\mathcal{S}_V^K\big(\PsiV_i - \Rk(\PsiV_i), \PsiV_j - \Rk(\PsiV_j)\big)\right]\ =\ \kron(\Id_2,\, \Hd^{\tra}\;\Hd)\,,$$
where $\Hd\in\R^{n_k^V\times n_k^V}$ is given by:
$$\Hd\ :=\ \Id_{n_k^V}\, -\, \left[\begin{array}{c}
\PV\\[-2ex]
\\ \hdashline 
\\[-2ex]
\vv{H}_3
\end{array}\right]\Rb\,.$$
The matrix $\vv{H}_3\in\R^{\mo\times\m1}$ is the submatrix of $\MASS$,
whose index range is $[1,\mo]\times[1,\m1]$. Furthermore, let $\vv{v}_1^v,\vv{v}_2^v,\ldots,\vv{v}_{d_K}^v$
be the vertices of $K$ and, for each $e\in\partial K$, let $\vv{v}_1^e,\vv{v}_2^e,\ldots,\vv{v}_k^e$
be the $k$ uniformly spaced points on $e$ (sorted with the respective orientation), which define
$m^{V}_{e}$ (cf. \eqref{eqn:dof-Vkv}). Then, gathering all these points in a matrix $\vv{p}:=[p_{ij}]\in\R^{2\times ((k+1)d_K)}$
as follow:
$$\vv{p}\ :=\ \left[\begin{array}{c:c:c:c:c}
\vv{v}_1^v,\, \ldots,\, \vv{v}_{d_K}^v & \vv{v}_1^{e_1},\, \ldots,\, \vv{v}_k^{e_1} & \vv{v}_1^{e_2},\, \ldots,\, \vv{v}_k^{e_2} & \cdots & \vv{v}_1^{e_{d_K}},\, \ldots,\, \vv{v}_k^{e_{d_K}}
\end{array}\right]\,,$$
we introduce the matrix $\PV := \big[\varphi_j^K(p_{1i},\,p_{2i})\big]\in\R^{((k+1)d_K)\times\m1}$,
which complete the previous definition of $\Hd$.

On the other hand, the operator $\DcK\in\R^{2n_k^V\times 2n_k^V}$
is given by
$$\DcK\ :=\ \left[\int_{\partial K\cap\Gamma}\PsiV_i\cdot\PsiV_j\right]\ =\kron\left(\Id_2,\, \sum_{e\in\partial K\cap\Gamma}\left[\int_{e}\psiV_i\,\psiV_j\right]\right)\,,$$
where, employing \eqref{eqn:interpLagrange}, it follows that
\begin{eqnarray*}
\DcK & = & \kron\left(\Id_2,\, \sum_{e\in\partial K\cap\Gamma}(\CLe{})^{\tra}\,\left[\int_{e}\mathcal{L}_i^e\,\mathcal{L}_j^e\right]\,\CLe{}\right)\\
     & = & \kron\left(\Id_2,\, \sum_{e\in\partial K\cap\Gamma}h_e\,(\CLe{})^{\tra}\,\MassL\,\CLe{}\right).
\end{eqnarray*}

We end this section by considering the operator $\vv{b}_3^K\in\R^{2n_k^V\times 1}$,
which employing \eqref{eqn:pg} is given by:
\begin{eqnarray*}
\vv{b}_3^K & := & \left[\int_{\partial K\cap\Gamma}\PsiV_i\cdot\vv{g}\right]\ =\ \sum_{e\in\partial K\cap\Gamma}\left[\begin{array}{c}
\displaystyle\left[\int_e\psiV_i\,g_1\right]\\[-2ex]
\\ \hdashline
\\[-2ex]
\displaystyle\left[\int_e\psiV_i\,g_2\right]
\end{array}\right]\\[1ex]
& = & \sum_{e\in\partial K\cap\Gamma}\left[\begin{array}{c}
\displaystyle(\CLe{})^{\tra}\left[\int_e\mathcal{L}_i^e\,\phi_j^e\right]\Pg{1}\\[-2ex]
\\ \hdashline
\\[-2ex]
\displaystyle(\CLe{})^{\tra}\left[\int_e\mathcal{L}_i^e\,\phi_j^e\right]\Pg{2}
\end{array}\right]\ =\ \sum_{e\in\partial K\cap\Gamma}\left[\begin{array}{c}
\displaystyle (\CLe{})^{\tra}\,\MCBe\,\Pg{1}\\[-2ex]
\\ \hdashline
\\[-2ex]
(\CLe{})^{\tra}\,\MCBe\,\Pg{2}
\end{array}\right].
\end{eqnarray*}
Furthermore, using \eqref{eqn:pf}, we describe the implementation of the operator:
\begin{eqnarray*}
\vv{b}_4^K & := & \left[\int_K\PsiV_i\cdot\Pk(\vv{f})\right]\ =\ \left[\begin{array}{c}
\displaystyle \left[\int_KP_k^K(\psiV_i)\,\varphi_j^K\right]\Pf{1}\\
\\[-2ex] \hdashline
\\[-2ex]
\displaystyle \left[\int_KP_k^K(\psiV_i)\,\varphi_j^K\right]\Pf{2}
\end{array}\right]\\
& = & \left[\begin{array}{c}
\displaystyle (\PU)^{\tra}\left[\int_K\varphi_i^K\,\varphi_j^K\right]\Pf{1}\\
\\[-2ex] \hdashline
\\[-2ex]
\displaystyle (\PU)^{\tra}\left[\int_K\varphi_i^K\,\varphi_j^K\right]\Pf{2}
\end{array}\right]\ =\ \left[\begin{array}{c}
\displaystyle (\PU)^{\tra}\,\Mass\,\Pf{1}\\
\\[-2ex] \hdashline
\\[-2ex]
\displaystyle (\PU)^{\tra}\,\Mass\,\Pf{2}
\end{array}\right]\\
& = & \left[\begin{array}{c}
\displaystyle (\PU)^{\tra}\left[\int_K\varphi_i^K\,f_1\right]\\
\\[-2ex] \hdashline
\\[-2ex]
\displaystyle (\PU)^{\tra}\left[\int_K\varphi_i^K\,f_2\right]
\end{array}\right]\in\R^{2n_k^V\times 1}\,.
\end{eqnarray*}

\subsection{Operators related to the elements of $\Hm_k^K$ and $\vv{V}_k^K$}

Next, we define the following operators:
	\begin{itemize}
	\item $\displaystyle\BK\ :=\ \left[\int_K\bdiv(\PsiH_i)\cdot\PsiV_j\right]\in\R^{2n_k^H\times 2n_k^V}$
	\item $\displaystyle\EK\ :=\ \left[\int_K\bPk(\nabla\PsiV_i):\big[\bPk(\PsiH_j)\big]^{\td}\right]\in\R^{2n_k^V\times 2n_k^H}$
	\end{itemize}
where, using the matrices $\PU\in\R^{m\times n_k^V}$ and $\PGU\in\R^{(4m)\times(2n_k^V)}$
defined in Section \ref{sec:PL2u}, we can assemble $\BK\in\R^{2n_k^H\times 2n_k^V}$
and $\EK\in\R^{2n_k^V\times 2n_k^H}$ similar to previous operators.
Indeed, we have that
\begin{eqnarray*}
\BK & := & \left[\int_K\bdiv(\PsiH_i)\cdot\PsiV_j\right]\ =\ \kron\left(\Id_2,\, \left[\int_K\div(\psiH_i)\,\psiV_j\right]\right)\\[1ex]
& = & \kron\left(\Id_2,\, (\Di)^{\tra}\,\left[\int_K\varphi_i^K\,\psiV_j\right]\right)\ =\ \kron\left(\Id_2,\, (\Di)^{\tra}\,\left[\int_K\varphi_i^K\,P_k^K(\psiV_j)\right]\right)\\[1ex]
& = & \kron\left(\Id_2,\, (\Di)^{\tra}\,\left[\int_K\varphi_i^K\,\varphi_j^K\right]\,\PU\right)\ =\ \kron\left(\Id_2,\, (\Di)^{\tra}\,\Mass\,\PU\right),
\end{eqnarray*}
and
\begin{eqnarray*}
\EK & := & \left[\int_K\bPk(\nabla\PsiV_i):\big[\bPk(\PsiH_j)\big]^{\td}\right]\ =\ (\PGU)^{\tra}\,\left[\int_K\Phi_i^K:(\Phi_j^K)^{\td}\right]\,\Pb\\[1ex]
& = & (\PGU)^{\tra}\,\kron\left(\vv{M}_{\text{dev}},\, \Mass\right)\,\Pb\,,
\end{eqnarray*}
where the matrix $\vv{M}_{\text{dev}}\in\R^{4\times 4}$ is defined in \eqref{mat:Mdev}.

\section{A particular example: the Navier-Stokes problem}\label{sec:navier-stokes}

In this section we show how to use the previously defined discrete operators in
a particular formulation. More precisely, we present some specific aspects
on the computational implementation of a mixed virtual element method
for the two-dimensional pseudostress-velocity formulation of the Navier-Stokes
equations with Dirichlet boundary conditions. Indeed, the formulation used
below was originally proposed and analyzed in \cite{gms-vem-ns}. Here we
recall the continuous and discrete formulations, and propose an algorithm
for the assembly of the associated global linear system for the Newton's
iteration. Finally, a numerical example illustrating the performance
of the mixed-VEM scheme and confirming these theoretical rates is presented.

\subsection{The continuous problem}

We begin, recalling from \cite[Section 2]{gms-vem-ns}, the boundary value problem
of interest. Indeed, letting $\Omega\subset\R^2$ be a bounded polygonal domain with
boundary $\Gamma$, we consider the stationary Navier-Stokes equations with nonhomogeneous
Dirichlet boundary conditions. More precisely, given a volume force $\fv\in\Lv(\Omega)$
and a Dirichlet datum $\bg\in\Hv^{1/2}(\Gamma)$, we seek a vector field (the velocity)
$\bu$ and a scalar field (the pressure) $p$ of a fluid on $\Omega$, such that
\begin{eqnarray}\label{eqn:modelproblem}
\begin{array}{c}
-\mu\Delta\bu\, +\, (\nabla\bu)\,\bu\, +\, \nabla p\ =\ \fv \qin \Omega\,,\qquad \div(\bu)\ =\ 0 \qin \Omega\,,\\[1ex]
\bu\ =\ \bg \qon \Gamma\,,\qan\displaystyle\int_{\Omega}p\ =\ 0\,,
\end{array}
\end{eqnarray}
where $\mu > 0$ is the viscosity constant. In addition, it is important to recall
here that the incompressibility condition given by the second equation of
\eqref{eqn:modelproblem}, establishes that the datum $\bg$ satisfies the compatibility
condition $\int_{\Gamma}\bg\cdot\nv = 0$, where $\nv$ stands for the unit outward
normal at $\Gamma$.

On the other hand, defining the constant $c := -\frac{1}{2|\Omega|}\|\bu\|_{0,\Omega}^2$
and $\Id_2$ the identity matrix of $\R^{2\times 2}$, we introduce the pseudostress tensor
(see \cite{cgot-SINUM-2016,cot-MC-2017,cgo-NMPDE-2016})
\begin{equation}\label{eqn:pseudostress}
\bsig\ :=\ \mu\nabla\bu\, -\, \bu\otimes\bu\, -\, (p+c)\,\Id_2\qin\Omega\,,
\end{equation}
which allows us to arrive at the equivalent system: Find the pseudostress $\bsig$ and the velocity
$\bu$ such that
\begin{equation}\label{eqn:modelproblem2}
\begin{array}{c}
\bsig^{\td}\ =\ \mu\nabla\bu\, -\, (\bu\otimes\bu)^{\td}\qin\Omega\,,\qquad -\bdiv(\bsig)\ =\ \fv\qin\Omega\,,\\[1ex]
\bu\ =\ \bg\qon\Gamma\,, \qan \displaystyle\int_\Omega\tr(\bsig)\ =\ 0\,,
\end{array}
\end{equation}
where the pseudostress variable has eliminated the pressure from the original model
\eqref{eqn:modelproblem}, which can be recover through the postprocessing formula:
\begin{equation*}\label{eqn:pressure}
p\ =\ -\frac{1}{2}\,\big\{\tr(\bsig)\, +\, \tr(\bu\otimes\bu)\big\}\, -\, c\qin\Omega\,,
\end{equation*}
which is obtained from \eqref{eqn:pseudostress} and the incompressibility condition.

Next, we let $\Hm_0(\bdiv;\Omega) := \left\{\btau\in\Hm(\bdiv;\Omega)\, :\, \int_{\Omega}\tr(\btau) = 0\right\}$,
and recall from \cite[Section 2]{gms-vem-ns} the following redundant terms:
\begin{equation*}\label{eqn:aumentos}
\begin{array}{rcll}
\displaystyle \kappa_1\int_{\Omega}\bdiv(\bsig)\cdot\bdiv(\btau) & = & \displaystyle -\kappa_1\int_{\Omega}\fv\cdot\bdiv(\btau)\quad & \forall\ \btau\in\Hm_0(\bdiv;\Omega)\,,\\[2ex]
\displaystyle \kappa_2\int_{\Omega}\Big\{\mu\nabla\bu - \bsig^{\td} - (\bu\otimes\bu)^{\td}\Big\}:\nabla\bv & = & 0 & \forall\ \bv\in\Hv^1(\Omega)\,,\\[2ex]
\displaystyle \kappa_3\int_{\Gamma}\bu\cdot\bv & = & \displaystyle\kappa_3\int_{\Gamma}\bg\cdot\bv &  \forall\ \bv\in\Hv^1(\Omega)\,,
\end{array}
\end{equation*}
where, according to \cite[Theorem 2.1]{gms-vem-ns}, the parameters $\kappa_1,\kappa_2,\kappa_3$
must satisfy that $\kappa_1,\kappa_3 >0$ and $0 < \kappa_2 < 2\mu$. Then, we consider the
continuous formulation of \eqref{eqn:modelproblem2} introduced in \cite[Section 2]{gms-vem-ns},
whose well-posedness has been established in \cite[Theorem 2.1]{gms-vem-ns}. More precisely, we
seek $\vec{\bsig} := (\bsig,\bu)\in\mathbbm{X} := \Hm_0(\bdiv;\Omega)\times\Hv^1(\Omega)$ such that
\begin{equation}\label{eqn:continuous-scheme}
a(\vec{\bsig}, \vec{\btau})\ +\ b(\bu; \vec{\bsig}, \vec{\btau})\ =\ F(\vec{\btau})\qquad\forall\ \vec{\btau} := (\btau,\bv)\in\mathbbm{X}\,,
\end{equation}
where $a:\mathbbm{X}\times\mathbbm{X}\to\R$ is the bilinear form
\begin{equation*}\label{eqn:operadorA}
\begin{array}{rcl}
a(\vec{\bzeta},\vec{\btau}) & := & \displaystyle\int_{\Omega}\bzeta^{\td}:\btau^{\td}\, +\, \kappa_1\int_{\Omega}\bdiv(\bzeta)\cdot\bdiv(\btau)\, +\, \kappa_2\mu\int_{\Omega}\nabla\bw:\nabla\bv\, +\, \kappa_3\int_{\Gamma}\bw\cdot\bv\\[2ex]
&& \displaystyle -\ \mu\int_{\Omega}\bv\cdot\bdiv(\bzeta)\, +\, \mu\int_{\Omega}\bw\cdot\bdiv(\btau)\, -\, 
\kappa_2\int_{\Omega}\bzeta^\td:\nabla\bv
\end{array}
\end{equation*}
for all $\vec{\bzeta} := (\bzeta,\bw)$, $\vec{\btau} := (\btau,\bv)\in\mathbbm{X}$,
$F:\mathbbm{X}\to\R$ is the linear functional
\begin{equation*}\label{eqn:funcional}
F(\vec{\btau})\ :=\ \mu\inn{\btau\nv,\bg}_{\Gamma}\, -\, \kappa_1\int_{\Omega}\fv\cdot\bdiv(\btau)\, +\, \mu\int_{\Omega}\fv\cdot\bv\, +\, \kappa_3\int_{\Gamma}\bg\cdot\bv\,,
\end{equation*}
for all $\vec{\btau} := (\btau,\bv)\in\mathbbm{X}$, and given $\bz\in\Hv^1(\Omega)$, 
$b(\bz;\,\cdot\,,\,\cdot\,):\mathbbm{X}\times\mathbbm{X}\to\R$ is the bilinear form
\begin{equation*}\label{eqn:operadorB}
b(\bz;\vec{\bzeta}, \vec{\btau})\ :=\ \int_{\Omega}(\bw\otimes\bz)^\td : \big\{\btau - \kappa_2\nabla\bv\big\}\,,
\end{equation*}
for all $\vec{\bzeta} := (\bzeta,\bw)$, $\vec{\btau} := (\btau,\bv)\in\mathbbm{X}$.

\subsection{The mixed-VEM formulation}\label{sec:model}

Now, given an integer $k \geq 0$, we consider the virtual element subspace
$\mathbbm{X}_h$ of $\mathbbm{X} := \Hm_0(\bdiv;\Omega)\times\Hv^1(\Omega)$
given by
\begin{equation}\label{eqn:defXh}
\mathbbm{X}_h\ :=\ \Big\{(\btau,\vv{v})\in\mathbbm{X}\,:\quad  \btau|_K\in\Hm_k^K \;\text{ and }\; \vv{v}|_K\in\vv{V}_k^K\quad \forall\ K\in\Th\Big\}\,,
\end{equation}
where, for each $K\in\Th$, the virtual subspaces $\Hm_k^K$ and $\vv{V}_k^K$
are defined in \eqref{eqn:defHh} and \eqref{eqn:defVh}, respectively.
In turn, we now aim to define the nonlinear mixed virtual element scheme
associated with \eqref{eqn:continuous-scheme} and introduced in \cite[Section 5]{gms-vem-ns}.
That is, we seek $\vec{\bsig}_h := (\bsig_h,\bu_h)\in\mathbbm{X}_h$ such that
\begin{equation}\label{eqn:discrete-scheme}
a_h(\vec{\bsig_h}, \vec{\btau}_h)\ +\ b_h(\bu_h; \vec{\bsig}_h, \vec{\btau}_h)\ =\ F_h(\vec{\btau}_h)\qquad\forall\ \vec{\btau}_h := (\btau_h,\bv_h)\in\mathbbm{X}_h\,.
\end{equation}
The equation \eqref{eqn:discrete-scheme} will be one of the few expression where we use the
subscript $h$ for the elements in $\mathbbm{X}_h$. In what follows we mostly omit that
subscript in order to simplify the notation. Thus, we have that $a_h:\mathbbm{X}_h\times\mathbbm{X}_h\to\R$
is the bilinear form defined by:
\begin{eqnarray*}
a_h(\vec{\bzeta}, \vec{\btau}) & := & \sum_{K\in\Th}\Bigg\{\int_K\big[\bPk(\bzeta)\big]^{\td}:\big[\bPk(\btau)\big]^{\td}\ +\ \mathcal{S}_H^K\big(\bzeta - \bPk(\bzeta), \btau - \bPk(\btau)\big)\\[2ex]
& & +\ \kappa_1\int_K\bdiv(\bzeta)\cdot\bdiv(\btau)\ +\ \mu\int_K\bw\cdot\bdiv(\btau)\ -\ \mu\int_K\bv\cdot\bdiv(\bzeta)\\[2ex]
& & -\ \kappa_2\int_K\big[\bPk(\bzeta)\big]^{\td}:\bPk(\nabla\bv)\ +\ \kappa_2\mu\int_K\nabla\Rk(\bw):\nabla\Rk(\bv)\\[2ex]
& & +\ \mathcal{S}_V^K\big(\bw - \Rk(\bw), \bv - \Rk(\bv)\big)\ +\ \kappa_3\int_{\partial K\cap\Gamma}\bw\cdot\bv\Bigg\}
\end{eqnarray*}
for all $\vec{\bzeta} := (\bzeta,\bw)$, $\vec{\btau} := (\btau,\bv)\in\mathbbm{X}_h$,
$F_h:\mathbbm{X}_h\to\R$ is the linear functional
\begin{eqnarray*}\label{funcional}
F_h(\vec{\btau}) & := & \sum_{K\in\Th}\Bigg\{\mu\int_{\partial K\cap\Gamma}\btau\nv\cdot\vv{g}\ -\ \kappa_1\int_K\vv{f}\cdot\bdiv(\btau)\ +\ \kappa_3\int_{\partial K\cap\Gamma}\vv{g}\cdot\bv\ +\ \mu\int_K\Pk(\vv{f})\cdot\bv\Bigg\}
\end{eqnarray*}
for all $\vec{\btau} := (\btau,\bv)\in\mathbbm{X}_h$, and given $\vv{z}\in\Hv^1(\Omega)$ such
that $\vv{z}|_K\in\vv{V}_k^K$ for all $K\in\Th$, the bilinear form $b_h(\vv{z};\,\cdot\,,\,\cdot\,):\mathbbm{X}_h\times\mathbbm{X}_h\to\R$
is given by
\begin{eqnarray}\label{eqn:operadorBh}
b_h(\vv{z}; \vec{\bzeta}, \vec{\btau}) & := & \sum_{K\in\Th}\Bigg\{\int_K\big[\Pk(\bw)\otimes\Pk(\vv{z})\big]^{\td}:\big[\bPk(\btau) - \kappa_2\bPk(\nabla\bv)\big]\Bigg\}
\end{eqnarray}
for all $\vec{\bzeta} := (\bzeta, \bw)$, $\vec{\btau} := (\btau,\bv)\in\mathbbm{X}_h$.

The bilinear forms $\mathcal{S}_H^K$ and $\mathcal{S}_V^K$ are defined in \eqref{eqn:defSH}
and \eqref{eqn:defSV}, respectively. Moreover, we recall here that $\Pk:\Lv(K)\to\Pv_k(K)$
and $\bPk:\Lm(K)\to\Pm_k(K)$ are the corresponding $\Lv(K)$ and $\Lm(K)$ orthogonal projections
(see at the end of Section \ref{sec:dofH}).
On the other hand, under suitable assumptions, the discrete scheme \eqref{eqn:discrete-scheme}
has a unique solution, which was proved in \cite[Theorem 5.1]{gms-vem-ns}, whereas in
\cite[Theorem 5.3]{gms-vem-ns} the respective a priori error estimates were established.

\subsection{The Newton's iteration and the linear system assembly}

In this section, Newton's method is described as an option to solve
the discrete scheme \eqref{eqn:discrete-scheme}, which as usual requires
the assembly and resolution of a series of linear systems. Thus, we now
aim to propose the following Newton's iteration for the discrete scheme
\eqref{eqn:discrete-scheme}: Given $\vec{\bsig}^{(0)}_h := (\bsig^{(0)}_h,\bu^{(0)}_h)\in\widetilde{\mathbbm{X}}_h$
and $\xi^{(0)}_h\in\R$, for each integer $s \geq 0$, we apply the iteration:
\begin{enumerate}
\item Find $\vec{\bzeta}^{(s)}_h := (\bzeta^{(s)}_h,\bw^{(s)}_h)\in\widetilde{\mathbbm{X}}_h$ and $\eta^{(s)}_h\in\R$ such that
	\begin{equation}\label{eqn:newton-scheme}
	\begin{array}{rcl}
	\displaystyle a_h(\vec{\bzeta}^{(s)}_h, \vec{\btau}_h)\, +\, \mathcal{D}b_h(\bu^{(s)}_h; \vec{\bzeta}^{(s)}_h, \vec{\btau}_h)
	\, +\, \eta^{(s)}_h\int_{\Omega}\tr(\btau_h) & = & F_h(\vec{\btau}_h)\, -\, a_h(\vec{\bsig}^{(s)}_h, \vec{\btau}_h)\\
& & -\ b_h(\bu^{(s)}_h; \vec{\bsig}^{(s)}_h, \vec{\btau}_h)\\[1ex]
& & -\ \displaystyle\xi^{(s)}_h\int_{\Omega}\tr(\btau_h)\,,\\[2ex]
	\displaystyle\lambda_h\int_{\Omega}\tr(\bzeta^{(s)}_h) & = & \displaystyle -\lambda_h\int_{\Omega}\tr(\bsig^{(s)}_h)\,,
	\end{array}
	\end{equation}
for all $\vec{\btau}_h:=(\btau_h,\bv_h)\in\widetilde{\mathbbm{X}}_h$ and for all
$\lambda_h\in\R$, where $\vec{\bsig}^{(s)}_h := (\bsig^{(s)}_h,\bu^{(s)}_h)\in\widetilde{\mathbbm{X}}_h$.
\item Compute $\vec{\bsig}^{(s+1)}_h := \vec{\bsig}^{(s)}_h\, +\, \vec{\bzeta}^{(s)}_h$ and $\xi^{(s+1)}_h := \xi^{(s)}_h\, +\, \eta^{(s)}_h$.
\end{enumerate}
Here, for each $\vv{z}\in\Hv^1(\Omega)$ such that $\vv{z}|_K\in\vv{V}_k^K$ for all $K\in\Th$,
the bilinear form $\, \mathcal{D}b_h(\vv{z};\,\cdot\,,\,\cdot\,):\widetilde{\mathbbm{X}}_h\times\widetilde{\mathbbm{X}}_h\to\R$
is the G$\hat{\text{a}}$teaux derivative of $b_h(\vv{z};\,\cdot\,,\,\cdot\,)$
(cf. \eqref{eqn:operadorBh}), given by
\begin{eqnarray*}\label{eqn:operadorDB}
\mathcal{D}b_h(\vv{z}; \vec{\bzeta}, \vec{\btau}) & := & \sum_{K\in\Th}\bigg\{\int_K\big[\Pk(\bw)\otimes\Pk(\vv{z}) + \Pk(\vv{z})\otimes\Pk(\bw)\big]^{\td}:\big[\bPk(\btau) - \kappa_2\bPk(\nabla\bv)\big]\bigg\}\,,
\end{eqnarray*}
for all $\vec{\bzeta} := (\bzeta, \bw)$, $\vec{\btau} := (\btau,\bv)\in\widetilde{\mathbbm{X}}_h$.
Furthermore, in order to relax the restrictions of $\mathbbm{X}_h$ (cf. \eqref{eqn:defXh}),
$\xi\in\R$ is introduced as the Lagrange multiplier that allows us to extract the 
condition $\int_{\Omega}\tr(\bsig^{(s)}_h) = 0$ (see, e.g., \cite[eq. (5.3)]{gms-vem-ns}).
Then, we replace the space $\mathbbm{X}_h$ as follow:
\begin{equation*}
\widetilde{\mathbbm{X}}_h\ :=\ \Big\{(\btau,\vv{v})\in\Hm(\bdiv;\Omega)\times\Hv^1(\Omega)\,:\quad  \btau|_K\in\Hm_k^K \;\text{ and }\; \vv{v}|_K\in\vv{V}_k^K\quad \forall\ K\in\Th\Big\}\,,
\end{equation*}
which, is identical to $\mathbbm{X}_h$ (cf. \eqref{eqn:defXh}), except that the
null trace integral condition is no longer imposed on its elements.

Next, the global linear system associated to \eqref{eqn:newton-scheme} has the
matrix structure:
\begin{equation}\label{eqn:linearSystem}
\mathcal{D}\vvm{\mathcal{A}}^{(s)}\,
\left[\begin{array}{c}
\bzeta^{(s)}_h\\
\\[-2ex] \hdashline
\\[-2ex]
\bw^{(s)}_h\\
\\[-2ex] \hdashline
\\[-2ex]
\eta^{(s)}_h
\end{array}\right]\ =\ 
\vvm{b}\, -\, \vvm{\mathcal{A}}^{(s)}\,
\left[\begin{array}{c}
\bsig^{(s)}_h\\
\\[-2ex] \hdashline
\\[-2ex]
\bu^{(s)}_h\\
\\[-2ex] \hdashline
\\[-2ex]
\xi^{(s)}_h
\end{array}\right],
\end{equation}
where $\mathcal{D}\vvm{\mathcal{A}}^{(s)},\vvm{\mathcal{A}}^{(s)}\in\R^{N\times N}$
and $\vvm{b}\in\R^{N\times 1}$, with $s$ indicating the dependence of $\bu_h^{(s)}$. In
addition, $N$ is the size of the system \eqref{eqn:linearSystem} which is given
by (see \eqref{eqn:dof-Hkv} and \eqref{eqn:dof-Vkv}):
\begin{eqnarray}
N & \!\!:=\!\! & 2\cdot\underbrace{(k+1)}_{m_{q,\nv}^H}\cdot(\text{\# of edges in }\Th)\ +\ 2\!\cdot\!\bigg(\underbrace{\frac{(k+1)(k+2)}{2}-1}_{m_{q,\div}^H}\, +\, \underbrace{\frac{k(k+1)}{2}}_{m_{\vv{q},\rot}^H}\bigg)\!\cdot\!(\text{\# of elements in }\Th)\nonumber\\
  &    & + \ 2\cdot\underbrace{1}_{m^{V}_{i,v}}\cdot(\text{\# of nodes in }\Th)\ +\ 2\cdot\underbrace{k}_{m^{V}_{e}}\cdot(\text{\# of edges in }\Th)\nonumber\\
  &    & +\ 2\!\cdot\!\bigg(\underbrace{\frac{k(k+1)}{2}}_{m^{V}_{q,K}}\bigg)\!\cdot\!(\text{\# of elements in }\Th)\ +\ \underbrace{1}_{\xi}\,,\label{eqn:N}
\end{eqnarray}
which indicates that there is $N$ unknowns associated with the degrees of freedom.

On the other hand, as is usual in finite element methods, the explicit
construction of the coefficient matrix and the right-hand side vector
in the system \eqref{eqn:linearSystem}, is done by assembling local
discrete operators from each element $K\in\Th$. More precisely, for each
$K\in\Th$, consider the local version (i.e., the contribution of $K$) of
the system \eqref{eqn:linearSystem} given by:
\begin{eqnarray*}
\mathcal{D}\vvm{\mathcal{A}}^{(s)}\big|_K & = & \left[\begin{array}{c:c:c}
\AaK + \AbK + \kappa_1\,\AcK & \mu\,\BK + \mathcal{D}\vv{G}_{1,K}^{(s)} & \;\;\CK\;\;\\
& & \\[-2ex] \hdashline
& & \\[-2ex]
-\mu\,(\BK)^{\tra} - \kappa_2\,\EK & \kappa_2\mu\,\DaK + \DbK + \kappa_3\,\DcK + \mathcal{D}\vv{G}_{2,K}^{(s)} & \vv{0}\\
& & \\[-2ex] \hdashline
& & \\[-2ex]
(\CK)^{\tra} & \vv{0} & 0
\end{array}\right]\,,\\[2ex]
\vvm{\mathcal{A}}^{(s)}\big|_K & = & \left[\begin{array}{c:c:c}
\AaK + \AbK + \kappa_1\,\AcK & \mu\,\BK + \vv{G}_{1,K}^{(s)} & \;\;\CK\;\;\\
& & \\[-2ex] \hdashline
& & \\[-2ex]
-\mu\,(\BK)^{\tra} - \kappa_2\,\EK & \kappa_2\mu\,\DaK + \DbK + \kappa_3\,\DcK + \vv{G}_{2,K}^{(s)} & \vv{0}\\
& & \\[-2ex] \hdashline
& & \\[-2ex]
(\CK)^{\tra} & \vv{0} & 0
\end{array}\right]\,,
\end{eqnarray*}
and
$$\vvm{b}\big|_K\ =\ \left[\begin{array}{c}
\mu\,\vv{b}_1^K - \kappa_1\,\vv{b}_2^K\\
\\[-2ex] \hdashline
\\[-2ex]
\kappa_3\,\vv{b}_3^K + \mu\,\vv{b}_4^K\\
\\[-2ex] \hdashline
\\[-2ex]
0
\end{array}\right]\,.$$
The explicit construction of the discrete operators $\AaK$, $\AbK$, $\AcK$,
$\BK$, $\CK$, $\DaK$, $\DbK$, $\DcK$, $\EK$, $\vv{b}_1^K$, $\vv{b}_2^K$, $\vv{b}_3^K$,
and $\vv{b}_4^K$ are detailed in Section \ref{sec:operators}. Conversely,
the operators:
	\begin{itemize}
	\item $\displaystyle\vv{G}_{1,K}^{(s)}\ :=\ \left[\int_K\bPk(\PsiH_i):\big[\Pk(\PsiV_j)\otimes\Pk(\vv{u}^{(s)}_h)\big]^{\td}\right]\in\R^{2n_k^H\times 2n_k^V}$
	\item $\displaystyle\vv{G}_{2,K}^{(s)}\ :=\ \left[-\kappa_2\int_K\bPk(\nabla\PsiV_i):\big[\Pk(\PsiV_j)\otimes\Pk(\vv{u}^{(s)}_h)\big]^{\td}\right]\in\R^{2n_k^V\times 2n_k^V}$
	\item $\displaystyle\mathcal{D}\vv{G}_{1,K}^{(s)}\ :=\ \left[\int_K\bPk(\PsiH_i):\big[\Pk(\PsiV_j)\otimes\Pk(\vv{u}^{(s)}_h) + \Pk(\vv{u}^{(s)}_h)\otimes\Pk(\PsiV_j)\big]^{\td}\right]\in\R^{2n_k^H\times 2n_k^V}$
	\item $\displaystyle\mathcal{D}\vv{G}_{2,K}^{(s)}\ :=\ \left[-\kappa_2\int_K\bPk(\nabla\PsiV_i):\big[\Pk(\PsiV_j)\otimes\Pk(\vv{u}^{(s)}_h) + \Pk(\vv{u}^{(s)}_h)\otimes\Pk(\PsiV_j)\big]^{\td}\right]\in\R^{2n_k^V\times 2n_k^V}$
	\end{itemize}
will be defined in Section \ref{sec:opNoLineal}. At the moment, for each
$K\in\Th$, we assume that these were already calculated, in order to
describe below the assembly of the global linear system \eqref{eqn:linearSystem}.

According to the previous discussion, the global matrices $\mathcal{D}\vvm{\mathcal{A}}^{(s)}\in\R^{N\times N}$
and $\vvm{\mathcal{A}}^{(s)}\in\R^{N\times N}$, along with the vector $\vvm{b}\in\R^{N\times 1}$,
can be assembled through the following algorithm:
\begin{enumerate}[\quad\footnotesize 1.]
\item Define $\mathcal{D}\vvm{\mathcal{A}}^{(s)} := \vv{0}$, $\vvm{\mathcal{A}}^{(s)} := \vv{0}$,
      and $\vvm{b} := \vv{0}$
\item Define $m^H := \frac{1}{2}(k+1)(k+2)-1 + \frac{1}{2}k(k+1)$,\, and\, $m^V := \frac{1}{2}k(k+1)$
\item Define $w_0 := (\text{\# of edges in }\Th)\cdot 2(k+1)\, +\, (\text{\# of elements in }\Th)\cdot 2m^H$
\item For each $K\in\Th$ do:
\item\qquad Construct the discrete operators: $\AaK$, $\AbK$, $\AcK$, $\BK$,
            $\CK$, $\DaK$, $\DbK$, $\DcK$, $\EK$, $\vv{G}_{1,K}^{(s)}$,\\
     $\phantom{\qquad Construct the discrete operators:}$\; $\vv{G}_{2,K}^{(s)}$, $\mathcal{D}\vv{G}_{1,K}^{(s)}$, $\mathcal{D}\vv{G}_{2,K}^{(s)}$,
            $\vv{b}_1^K$, $\vv{b}_2^K$, $\vv{b}_3^K$, and $\vv{b}_4^K$
\item\qquad Define $n^H := (k + 1)d_K + m^H$,\, and\, $n^V := (k + 1)d_K + m^V$
\item\qquad Define $\vv{p}^H:=(p_i^H)\in\R^{n^H\times 1}$,\, and\, $\vv{q}^H:=(q_i^H)\in\R^{n^H\times 1}$
\item\qquad Define $\vv{p}^V:=(p_i^V)\in\R^{n^V\times 1}$,\, and\, $\vv{q}^V:=(q_i^V)\in\R^{n^V\times 1}$
\item\qquad For $e = 1$ until $d_K$ do:
\item\qquad\qquad Let $I_e$ be the global index of the local edge $e$ in $\Th$
\item\qquad\qquad Define $w := (I_e-1)\cdot 2(k+1)$
\item\qquad\qquad For $r = 1$ until $k+1$ do:
\item\qquad\qquad\qquad Define $i := (e-1)(k+1) + r$
\item\qquad\qquad\qquad Set $p_i^H := w + r$
\item\qquad\qquad\qquad Set $q_i^H := w + (k+1) + r$
\item\qquad\qquad End for of $r$
\item\qquad\qquad Let $I_n$ be the global index of the $e$th local node in $\Th$
\item\qquad\qquad Set $p_e^V := w_0 + 2(I_n-1) + 1$
\item\qquad\qquad Set $q_e^V := w_0 + 2(I_n-1) + 2$
\item\qquad\qquad Define $w := w_0\, +\, (\text{\# of nodes in }\Th)\cdot 2\, +\, (I_e-1)\cdot 2k$
\item\qquad\qquad If edge $e$ has positive orientation in $\Th$ do: 
\item\qquad\qquad\qquad For $r = 1$ until $k$ do:
\item\qquad\qquad\qquad\qquad Define $i := d_K + (e-1)\cdot k + r$
\item\qquad\qquad\qquad\qquad Set $p_i^V := w + r$
\item\qquad\qquad\qquad\qquad Set $q_i^V := w + k + r$
\item\qquad\qquad\qquad End for of $r$
\item\qquad\qquad Else
\item\qquad\qquad\qquad For $r = 1$ until $k$ do:
\item\qquad\qquad\qquad\qquad Define $i := d_K + (e-1)\cdot k + r$
\item\qquad\qquad\qquad\qquad Set $p_i^V := w + (k+1-r)$
\item\qquad\qquad\qquad\qquad Set $q_i^V := w + k + (k+1-r)$
\item\qquad\qquad\qquad End for of $r$
\item\qquad\qquad End if
\item\qquad End for of $e$
\item\qquad Let $I_K$ be the global index of the element $K$ in $\Th$ 
\item\qquad Define $w^H := (\text{\# of edges in }\Th)\cdot 2(k+1)\, +\, (I_K-1)\cdot 2m^H$
\item\qquad For $r = 1$ until $m^H$ do:
\item\qquad\qquad Define $i := (k+1)d_K + r$
\item\qquad\qquad Set $p_i^H := w^H + r$
\item\qquad\qquad Set $q_i^H := w^H + m^H + r$
\item\qquad End for of $r$
\item\qquad Define $w^V := w_0\, +\, (\text{\# of nodes in }\Th)\cdot 2\, +\, (\text{\# of edges in }\Th)\cdot 2k\, +\, (I_K-1)\cdot 2m^V$
\item\qquad For $r = 1$ until $m^V$ do:
\item\qquad\qquad Define $i := (k+1)d_K + r$
\item\qquad\qquad Set $p_i^V := w^V + r$
\item\qquad\qquad Set $q_i^V := w^V + m^V + r$
\item\qquad End for of $r$
\item\qquad Define $\vv{u}^H := \mat{\begin{smallmatrix}\vv{p}^H\\[0.5ex]\vv{q}^H\end{smallmatrix}}\in\R^{2n^H\times 1}$,
            and\, $\vv{u}^V := \mat{\begin{smallmatrix}\vv{p}^V\\[0.5ex]\vv{q}^V\end{smallmatrix}}\in\R^{2n^V\times 1}$
\item\qquad $\mathcal{D}\vvm{\mathcal{A}}^{(s)}(\vv{u}^H,\vv{u}^H)\, :=\, \mathcal{D}\vvm{\mathcal{A}}^{(s)}(\vv{u}^H,\vv{u}^H)\, +\, \AaK\, +\, \AbK\, +\, \kappa_1\,\AcK$
\item\qquad $\vvm{\mathcal{A}}^{(s)}(\vv{u}^H,\vv{u}^H)\, :=\, \vvm{\mathcal{A}}^{(s)}(\vv{u}^H,\vv{u}^H)\, +\, \AaK\, +\, \AbK\, +\, \kappa_1\,\AcK$
\item\qquad $\mathcal{D}\vvm{\mathcal{A}}^{(s)}(\vv{u}^H,\vv{u}^V)\, :=\, \mathcal{D}\vvm{\mathcal{A}}^{(s)}(\vv{u}^H,\vv{u}^V)\, +\, \mu\,\BK\, +\, \mathcal{D}\vv{G}_{1,K}^{(s)}$
\item\qquad $\vvm{\mathcal{A}}^{(s)}(\vv{u}^H,\vv{u}^V)\, :=\, \vvm{\mathcal{A}}^{(s)}(\vv{u}^H,\vv{u}^V)\, +\, \mu\,\BK\, +\, \vv{G}_{1,K}^{(s)}$
\item\qquad $\mathcal{D}\vvm{\mathcal{A}}^{(s)}(\vv{u}^H,N)\, :=\, \mathcal{D}\vvm{\mathcal{A}}^{(s)}(\vv{u}^H,N)\, +\, \CK$
\item\qquad $\vvm{\mathcal{A}}^{(s)}(\vv{u}^H,N)\, :=\, \vvm{\mathcal{A}}^{(s)}(\vv{u}^H,N)\, +\, \CK$
\item\qquad $\mathcal{D}\vvm{\mathcal{A}}^{(s)}(\vv{u}^V,\vv{u}^H)\, :=\, \mathcal{D}\vvm{\mathcal{A}}^{(s)}(\vv{u}^V,\vv{u}^H)\, -\, \mu\,(\BK)^{\tra}\, -\, \kappa_2\,\EK$
\item\qquad $\vvm{\mathcal{A}}^{(s)}(\vv{u}^V,\vv{u}^H)\, :=\, \vvm{\mathcal{A}}^{(s)}(\vv{u}^V,\vv{u}^H)\, -\, \mu\,(\BK)^{\tra}\, -\, \kappa_2\,\EK$
\item\qquad $\mathcal{D}\vvm{\mathcal{A}}^{(s)}(\vv{u}^V,\vv{u}^V)\, :=\, \mathcal{D}\vvm{\mathcal{A}}^{(s)}(\vv{u}^V,\vv{u}^V)\, +\, \kappa_2\mu\,\DaK\, +\, \DbK\, +\, \kappa_3\,\DcK\, +\, \mathcal{D}\vv{G}_{2,K}^{(s)}$
\item\qquad $\vvm{\mathcal{A}}^{(s)}(\vv{u}^V,\vv{u}^V)\, :=\, \vvm{\mathcal{A}}^{(s)}(\vv{u}^V,\vv{u}^V)\, +\, \kappa_2\mu\,\DaK\, +\, \DbK\, +\, \kappa_3\,\DcK\, +\, \vv{G}_{2,K}^{(s)}$
\item\qquad $\mathcal{D}\vvm{\mathcal{A}}^{(s)}(N,\vv{u}^H)\, :=\, \mathcal{D}\vvm{\mathcal{A}}^{(s)}(N,\vv{u}^H)\, +\, (\CK)^{\tra}$
\item\qquad $\vvm{\mathcal{A}}^{(s)}(N,\vv{u}^H)\, :=\, \vvm{\mathcal{A}}^{(s)}(N,\vv{u}^H)\, +\, (\CK)^{\tra}$
\item\qquad $\vvm{b}(\vv{u}^H)\, :=\, \vvm{b}(\vv{u}^H)\, +\, \mu\,\vv{b}_1^K\, -\, \kappa_1\,\vv{b}_2^K$
\item\qquad $\vvm{b}(\vv{u}^V)\, :=\, \vvm{b}(\vv{u}^V)\, +\, \kappa_3\,\vv{b}_3^K\, +\, \mu\,\vv{b}_4^K$
\item End for of $K$
\end{enumerate}

Regarding the procedure described previously, it is important to recall that the
construction of discrete operators of line 5 is described in Sections \ref{sec:operators}
and \ref{sec:opNoLineal}.
On the other hand, in lines from 7 to 48, we construct vectors $\vv{u}^H$ and $\vv{u}^V$,
in order to map the local degrees of freedom of $K$ to their corresponding location in the
global system \eqref{eqn:linearSystem}. In fact, $\vv{p}^H$ maps the first row of the tensor
$\bsig_h$, whereas $\vv{q}^H$ maps the second one. Similarly, $\vv{p}^V$ maps the first component
of the vector $\vv{u}_h$ and $\vv{q}^H$ maps its second component (see Figure \ref{fig:dof}).
Furthermore, the global assembly is performed in lines 49 through 62, where, in particular, the
notation $\vvm{\mathcal{A}}(\vv{i},\vv{j})$ is related to access of the block of $\vvm{\mathcal{A}}$
obtained through rows with indexes stored in $\vv{i}$ and columns with indexes stored in $\vv{j}$.
In particular, we remark here that MATLAB allows us to perform these operations in
natural way. Finally, observe that only the matrices $\mathcal{D}\vv{G}_1^{(s)}$,
$\mathcal{D}\vv{G}_2^{(s)}$, $\vv{G}_1^{(s)}$, and $\vv{G}_2^{(s)}$, depend on the previous
approximation and they are the main difference between the matrices $\mathcal{D}\vvm{\mathcal{A}}^{(s)}$
and $\vvm{\mathcal{A}}^{(s)}$, which can be used to improve the efficiency of the foregoing
assembly.

\subsection{The operators \boldmath $\vv{G}_{1,K}^{(s)}$, $\vv{G}_{2,K}^{(s)}$, $\mathcal{D}\vv{G}_{1,K}^{(s)}$, and $\mathcal{D}\vv{G}_{2,K}^{(s)}$}\label{sec:opNoLineal}

We now aim to describe the implementation of the last four operators
related to the nonlinearly of our problem. To do that, let $s\geq 0$ be
an integer representing the current iteration of the Newton's method.
In turn, let $\bu^{(s)}_{h,K}\in\vv{V}_k^K$ be the local approximation of
$\vv{u}$ in the $s$th iteration of the Newton's method, which satisfies that
\begin{equation}\label{eqn:uhK}
\bu^{(s)}_{h,K}\ :=\ \bu^{(s)}_h\big|_K\ =\ \sum_{i=1}^{2n_k^V}\beta_i^K\,\PsiV_i\,,
\end{equation}
where $\{\beta_i^K\}_{i=1}^{2n_k^V}$ are the respective local degrees of freedom
of $\vv{u}_h^{(s)}$ on $K$. In addition, we define
$$\vv{z}^{(s)}\ :=\ \Pk(\bu^{(s)}_{h,K})\ =\ \sum_{i=1}^{2m}\gamma_i^K\,\vvm{\varphi}_i^K\,.$$
It is important to remark here that, defining $\vvm{\beta}^K:=(\beta_1^K,\ldots,\beta_{2n_k^V}^K)^{\tra}$
and $\vvm{\gamma}^K:=(\gamma_1^K,\ldots,\gamma_{2m}^K)^{\tra}$, there holds
$$\vvm{\gamma}^K\ =\ \kron(\Id_2,\,\PU)\,\vvm{\beta}^K\,,$$
where $\PU$ is defined in \eqref{eqn:matPV}.

Now, employing the notation $\vv{z}^{(s)} := (z_1^{(s)},z_2^{(s)})^{\tra}$, such that
$$z_{\ell}^{(s)}\ :=\ \sum_{r=1}^mz_{\ell,r}^{(s)}\,\varphi_r^K\,,\qquad\text{for }\ell=1,2\,,$$
we introduce the matrices:
$$\MZ{\ell}\ :=\ \left[\int_Kz_{\ell}^{(s)}\varphi_i^K\varphi_j^K\right]\ =\ \sum_{r=1}^mz_{\ell,r}^{(s)}\,\left[\int_K\varphi_r^K\varphi_i^K\varphi_j^K\right]\in\R^{m\times m}\,,$$
with $\ell=1,2$, whose entries can be calculated by using the formula \eqref{eqn:formulaIntegral}.
On the other hand, using matrices $\MZ{1}$ and $\MZ{2}$, we now introduce the last
auxiliary matrices:
$$\MdPU\ :=\ \kron\left(\mat{\begin{smallmatrix}
\phantom{-}\frac{1}{2} & \phantom{-}0\\
\phantom{-}0 & \phantom{-}1\\
\phantom{-}0 & \phantom{-}0\\
-\frac{1}{2} & \phantom{-}0
\end{smallmatrix}},\, \MZ{1}\,\PU\right)\, +\, \kron\left(\mat{\begin{smallmatrix}
0 & -\frac{1}{2}\\
0 & \phantom{-}0\\
1 & \phantom{-}0\\
0 & \phantom{-}\frac{1}{2}
\end{smallmatrix}},\, \MZ{2}\,\PU\right)\in\R^{(4m)\times(2n_k^V)}\,,$$
and
$$\IMd0I\ :=\ \Id_{4m}\, +\, \kron\left(\mat{\begin{smallmatrix}
1 & \,0\, & \,0\, & 0\\[0.5ex]
0 & 0 & 1 & 0\\[0.5ex]
0 & 1 & 0 & 0\\[0.5ex]
0 & 0 & 0 & 1
\end{smallmatrix}},\, \Id_m\right).$$

Next, according to the above notation, the operator $\vv{G}_{1,K}^{(s)}\in\R^{2n_k^H\times 2n_k^V}$
is given by:
\begin{eqnarray*}
\vv{G}_{1,K}^{(s)} & := & \left[\int_K\bPk(\PsiH_i):\big[\Pk(\PsiV_j)\otimes\Pk(\vv{u}^{(s)}_h)\big]^{\td}\right]\ =\ (\Pb)^{\tra}\left[\int_K(\Phi_i^K)^{\td}:\big[\Pk(\PsiV_j)\otimes\vv{z}^{(s)}\big]\right]\\[1ex]
& = & (\Pb)^{\tra}\left[\begin{array}{c:c}
\displaystyle\frac{1}{2}\left[\int_Kz_1\varphi_i^K\,P_k^K(\psiV_j)\right] & \displaystyle-\frac{1}{2}\left[\int_Kz_2\varphi_i^K\,P_k^K(\psiV_j)\right]\\[-2ex]
\\ \hdashline
\\[-2ex]
\displaystyle\left[\int_Kz_2\varphi_i^K\,P_k^K(\psiV_j)\right] & \vv{0}_{m\times n_k^V}\\[-2ex]
\\ \hdashline
\\[-2ex]
\vv{0}_{m\times n_k^V} & \displaystyle\left[\int_Kz_1\varphi_i^K\,P_k^K(\psiV_j)\right]\\[-2ex]
\\ \hdashline
\\[-2ex]
\displaystyle-\frac{1}{2}\left[\int_Kz_1\varphi_i^K\,P_k^K(\psiV_j)\right] & \displaystyle\frac{1}{2}\left[\int_Kz_2\varphi_i^K\,P_k^K(\psiV_j)\right]
\end{array}\right]\ =\ (\Pb)^{\tra}\,\MdPU\,.
\end{eqnarray*}
Furthermore, in a similar way, we deduce that
$$\vv{G}_{2,K}^{(s)}\ :=\ \left[-\kappa_2\int_K\bPk(\nabla\PsiV_i):\big[\Pk(\PsiV_j)\otimes\Pk(\vv{u}^{(s)}_h)\big]^{\td}\right]\ =\ -\kappa_2\,(\PGU)^{\tra}\,\MdPU\in\R^{2n_k^V\times 2n_k^V}\,.$$

Finally, following the previous analysis, it is not difficult to obtain that 
\begin{eqnarray*}
\mathcal{D}\vv{G}_{1,K}^{(s)} & := & \left[\int_K\bPk(\PsiH_i):\big[\Pk(\PsiV_j)\otimes\Pk(\vv{u}^{(s)}_h) + \Pk(\vv{u}^{(s)}_h)\otimes\Pk(\PsiV_j)\big]^{\td}\right]\\[1ex]
& = & (\Pb)^{\tra}\,\IMd0I\,\MdPU\in\R^{2n_k^H\times 2n_k^V}\,,
\end{eqnarray*}
and
\begin{eqnarray*}
\mathcal{D}\vv{G}_{2,K}^{(s)} & := & \left[-\kappa_2\int_K\bPk(\nabla\PsiV_i):\big[\Pk(\PsiV_j)\otimes\Pk(\vv{u}^{(s)}_h) + \Pk(\vv{u}^{(s)}_h)\otimes\Pk(\PsiV_j)\big]^{\td}\right]\\[1ex]
& = & -\kappa_2\,(\PGU)^{\tra}\,\IMd0I\,\MdPU\in\R^{2n_k^V\times 2n_k^V}\,.
\end{eqnarray*}

\subsection{Calculable approximations of $\bsig$, $\bu$, and $p$}\label{sec:calcApp}

Once Newton's iteration is over, we can use the discrete operators described
in previous sections to find non-virtual approximations of all the unknowns
(see \cite[Sections 5.3 and 5.4]{gms-vem-ns}). Indeed, given $(\bsig_h^{(s)},\bu_h^{(s)})\in\mathbb{X}_h$
the final approximation of the solution of \eqref{eqn:discrete-scheme}, we
consider, for each $K\in\Th$, the vector $\mat{\begin{smallmatrix}
\vvm{\alpha}^K\\ \vvm{\beta}^K\\ \xi_h
\end{smallmatrix}}$ containing the local degrees of freedom of
$\bsig_h^{(s)}$, $\bu_h^{(s)}$ and $\xi_h^{(s)}$, respectively,
sorted as indicated in Figure \ref{fig:dof}. More precisely,
following \eqref{eqn:uhK}, we have
$$\bsig^{(s)}_h\big|_K\ =\ \sum_{i=1}^{2n_k^H}\alpha_i^K\,\PsiH_i \qan \bu^{(s)}_h\big|_K\ =\ \sum_{i=1}^{2n_k^V}\beta_i^K\,\PsiV_i\,,$$
where $\vvm{\alpha}^K = (\alpha_1^K, \ldots, \alpha_{2n_k^H}^K)^{\tra}$
and $\vvm{\beta}^K = (\beta_1^K, \ldots, \beta_{2n_k^V}^K)^{\tra}$.
Now, according to \cite[eq. (5.32)]{gms-vem-ns}), we introduce the fully computable
local approximations of $\bsig_h^{(s)}$ and $\bu_h^{(s)}$ given by
$$\widehat{\bsig}_h\big|_K\ :=\ \bPk(\bsig^{(s)}_h\big|_K)\ =\ \sum_{i=1}^{4m}a_i^K\,\Phi_i^K\,,$$
and
$$\widehat{\bu}_h\big|_K\ :=\ \Pk(\bu^{(s)}_h\big|_K)\ =\ \sum_{i=1}^{2m}b_i^K\,\vvm{\varphi}_i^K\,,$$
respectively. Moreover, it is quite simple to see that
$$\vvm{a}^K\ =\ \Pb\vvm{\alpha}^K\qan \vvm{b}^K\ =\ \kron(\Id_2,\,\Rb)\,\vvm{\beta}^K\,,$$
where $\vvm{a}^K = (a_1^K, \ldots, a_{4m}^K)^{\tra}$
and $\vvm{b}^K = (b_1^K, \ldots, b_{2m}^K)^{\tra}$. In addition, we
now present the following computable approximation of the pressure:
$$\widehat{p}_h\big|_K\ :=\ -\frac{1}{2}\,\tr\Big(\widehat{\bsig}_h\big|_K\, +\, \widehat{c}_h\,\Id_2\, +\, \widehat{\bu}_h\big|_K\otimes\widehat{\bu}_h\big|_K\Big)\,,$$
with $\widehat{c}_h := -\frac{1}{2|\Omega|}\,\|\widehat{\bu}_h\|_{0,\Omega}^2$.
Finally, in \cite[Section 5.4]{gms-vem-ns}, a second approximation $\widetilde{\bsig}_h$
of the pseudostress $\bsig$, which yields optimal rate of convergence in the broken
$\Hm(\bdiv;\Omega)$-norm is presented. Here, the calculation of $\widetilde{\bsig}_h$
is not presented, since it follows similarly as the previous operators
described before.

\subsection{Numerical results}

In this section we present a numerical experiment in order to illustrate the performance
of the mixed virtual element scheme \eqref{eqn:discrete-scheme} employing the
Newton's iteration introduced in \eqref{eqn:newton-scheme}. It allows us to validate
the operators introduced in Section \ref{sec:operators}, together with the numerical
experiments presented in recent papers about mixed-VEM schemes, which utilized our
implementation approach (see \cite{cgs-1,cgs-2,cgs-3,gms-brinkman,gms-boussinesq,ms-aposteriori-2020}).
We begin by recalling from \eqref{eqn:N} that $N$ stands for the total number of
degrees of freedom (unknowns) of \eqref{eqn:newton-scheme}. In addition, the individual
errors are defined by
$${\tt e}(\bsig)\ :=\ \|\bsig - \widehat{\bsig}_h\|_{0,\Omega}\,,\qquad{\tt e}(\bu)\ :=\ \|\bu - \widehat{\bu}_h\|_{0,\Omega}\,,\qquad {\tt e}(\widehat{\bu})\ :=\ \Bigg\{\sum_{K\in\Th}\|\bu - \widehat{\bu}_h\|_{1,K}^2\Bigg\}^{1/2}\,,$$
$${\tt e}(p)\ :=\ \|p - \widehat{p}_h\|_{0,\Omega}\,, \qan {\tt e}(\widetilde{\bsig})\ :=\ \Bigg\{\sum_{K\in\Th}\|\bsig - \widetilde{\bsig}_h\|_{\bdiv;K}^2\Bigg\}^{1/2}\,,$$
where $\widehat{\bsig}_h$, $\widehat{\bu}_h$, $\widehat{p}_h$, and $\widetilde{\bsig}_h$
are introduced in Section \ref{sec:calcApp}. In turn, the associated experimental rates
of convergence are given by
$${\tt r}(\cdot)\ :=\ \displaystyle\frac{\log\big({\tt e}(\cdot)\, /\, {\tt e}^{\prime}(\cdot)\big)}{\log(h\, /\, h^{\prime})}\,,$$
where ${\tt e}$ and ${\tt e}^{\prime}$ denote the corresponding errors for two consecutive
meshes with sizes $h$ and $h^{\prime}$, respectively.

The Newton method \eqref{eqn:newton-scheme} is solved by using a tolerance of $10^{-6}$
and taking as initial iteration the solution of the associated linear Stokes problem,
where four iterations were required to achieve the given tolerance. On the other hand,
the numerical results presented below were obtained using a MATLAB code, where the
corresponding linear systems were solved using its instruction ``{\textbackslash}"
as main solver.

\begin{figure}[h!t]\centering
\scalebox{0.45}{\includegraphics{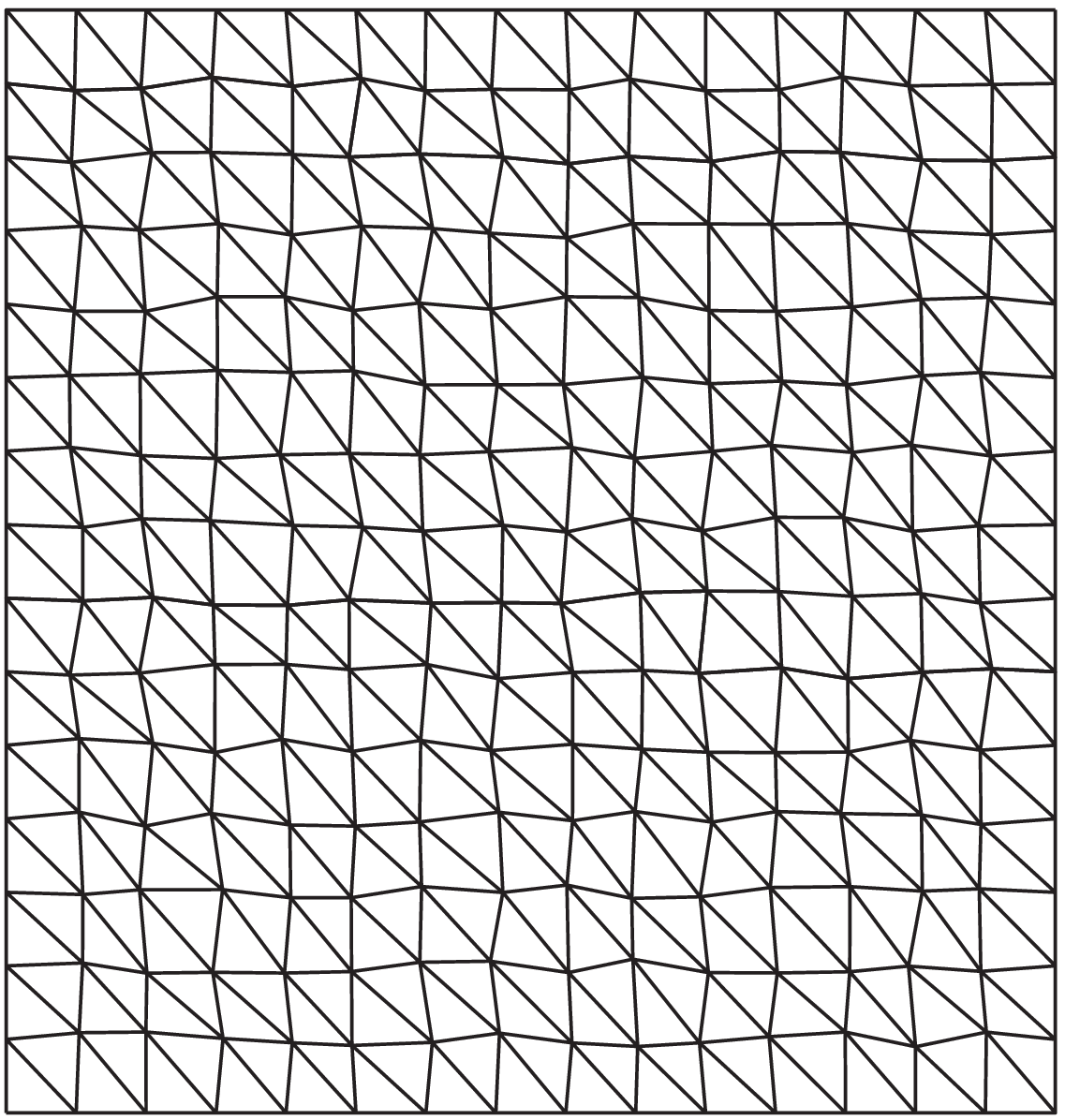}\qquad
                \includegraphics{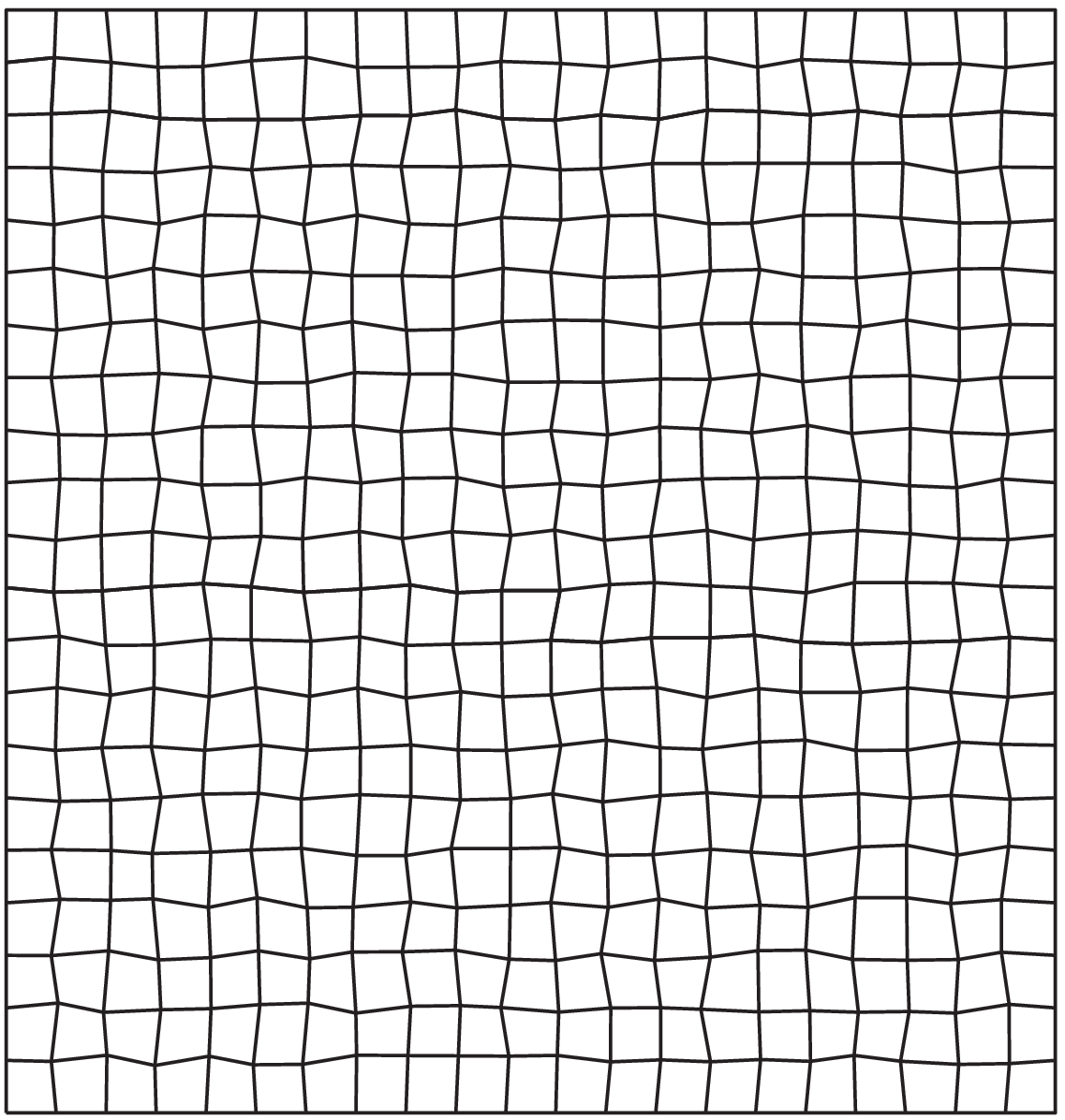}\qquad
                \includegraphics{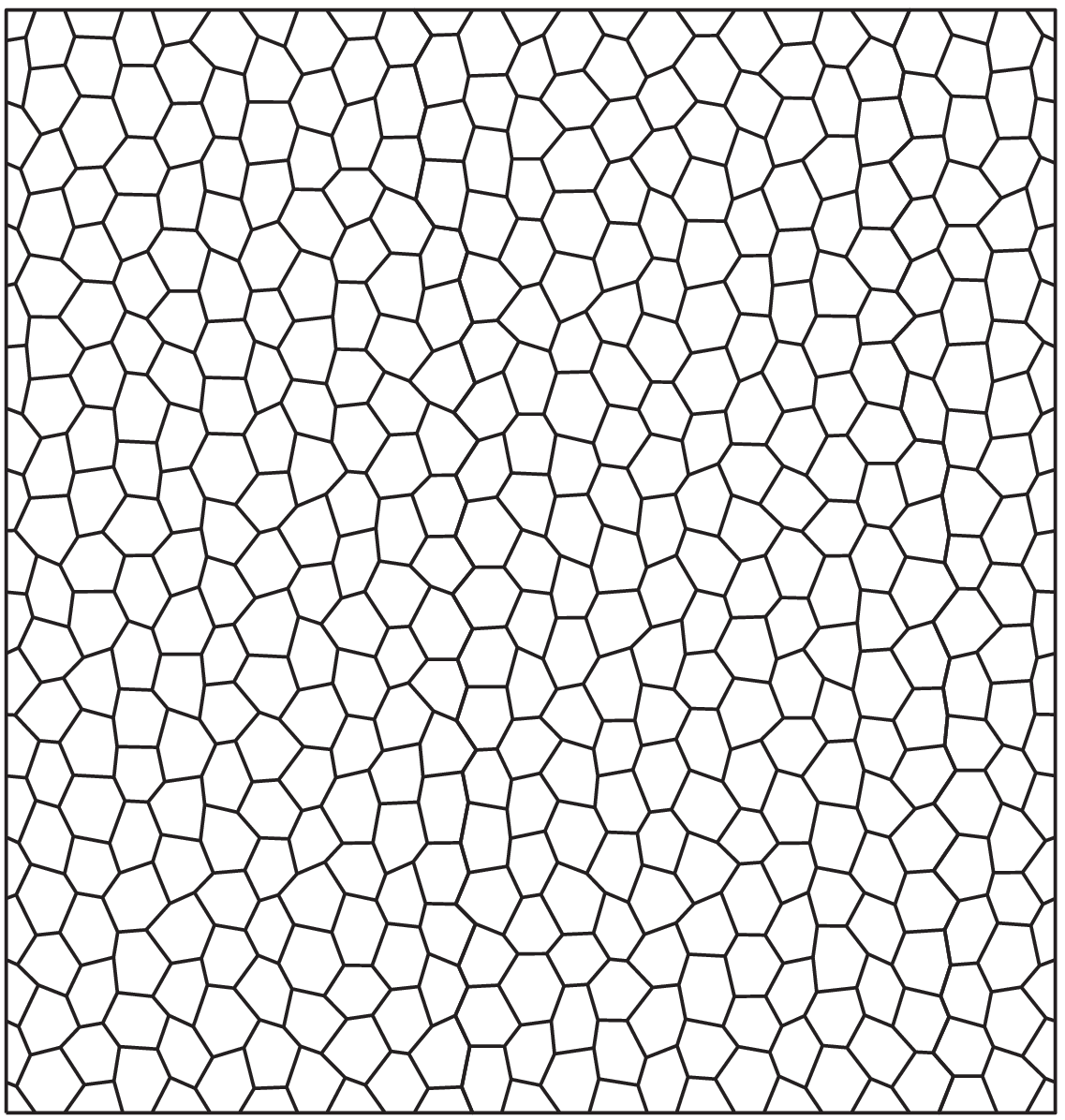}}
\caption{Examples of the meshes to be used in the calculations.}\label{fig:meshes}
\end{figure}

Next, we consider $\Omega := (-0.5, 1.5)\times(0,2)$, $\mu = 0.1$, and choose
the data $\vv{f}$ and $\vv{g}$ so that the exact solution is given by the flow from
\cite{kovaszany48}, that is,
$$\vv{u}(\vv{x})\ =\ \left(\begin{array}{c}
\displaystyle 1 - \exp(\lambda x_1)\cos(2\pi x_2)\\[1ex]
\displaystyle \frac{\lambda}{2\pi}\exp(\lambda x_1)\sin(2\pi x_2)
\end{array}\right)\qan
p(\vv{x})\ =\ \frac{1}{2}\exp(2\lambda x_1)\, -\, \frac{1}{8\lambda}\big\{\exp(3\lambda) - \exp(-\lambda)\big\}\,,$$
for all $\vv{x} := (x_1,x_2)^{\tra} \in \Omega$, where $\lambda := \frac{Re}{2} - \sqrt{\frac{Re^2}{4} + 4\pi^2}$
and $Re := \mu^{-1} = 10$ is the Reynolds number. Moreover, according to \cite[Theorem 2.1]{gms-vem-ns},
we set the parameters $\kappa_1=\kappa_2=\kappa_3 = 0.1$. In addition, we employ polynomial
degrees $k\in\{0,1,2\}$, and for the decompositions of $\Omega$ used in our computations, we
consider triangles, distorted squares and distorted hexagons, as illustrated in the
Figure \ref{fig:meshes}.

\begin{figure}[h!t]\centering
\scalebox{0.45}{\includegraphics{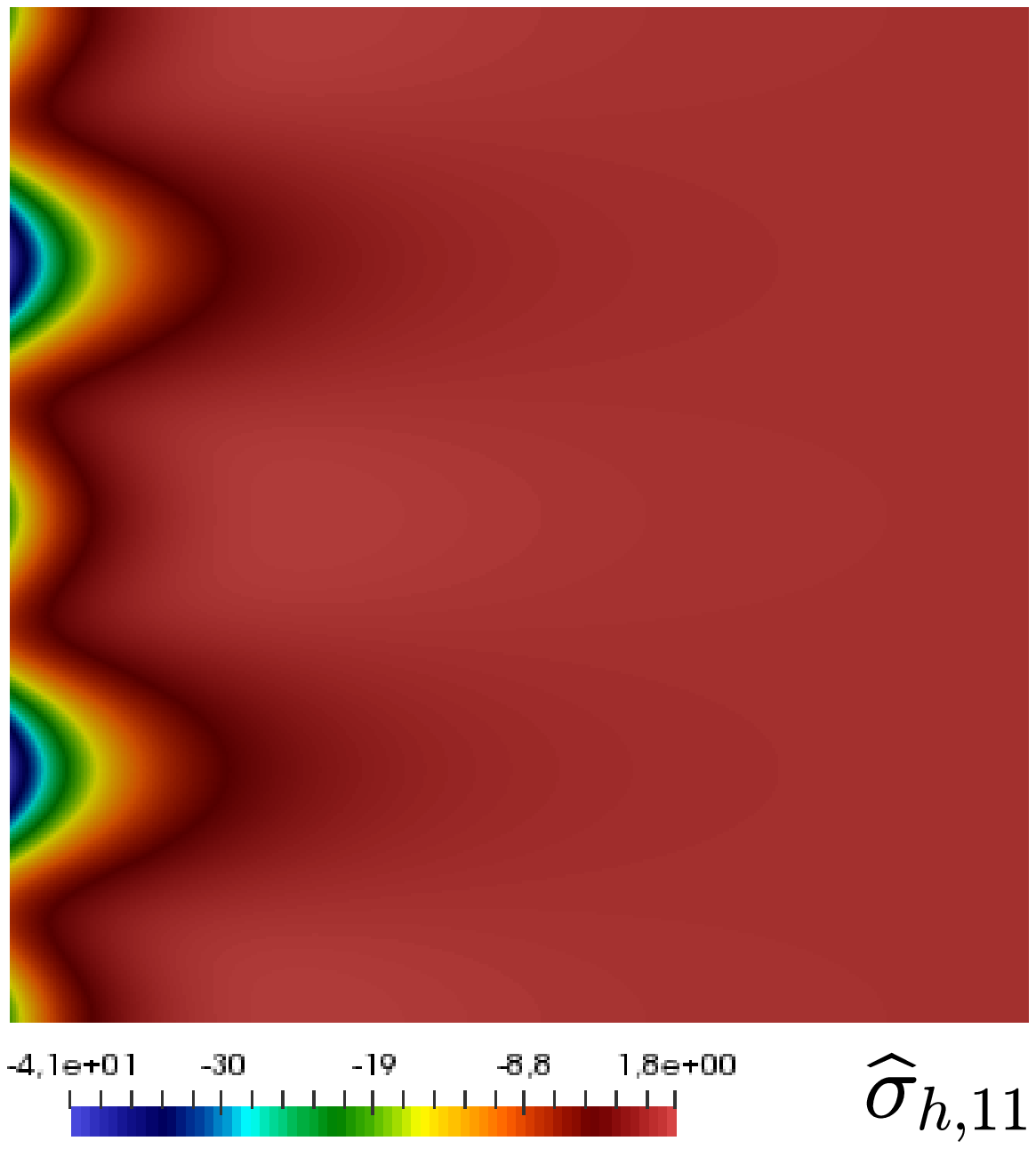}\qquad
                \includegraphics{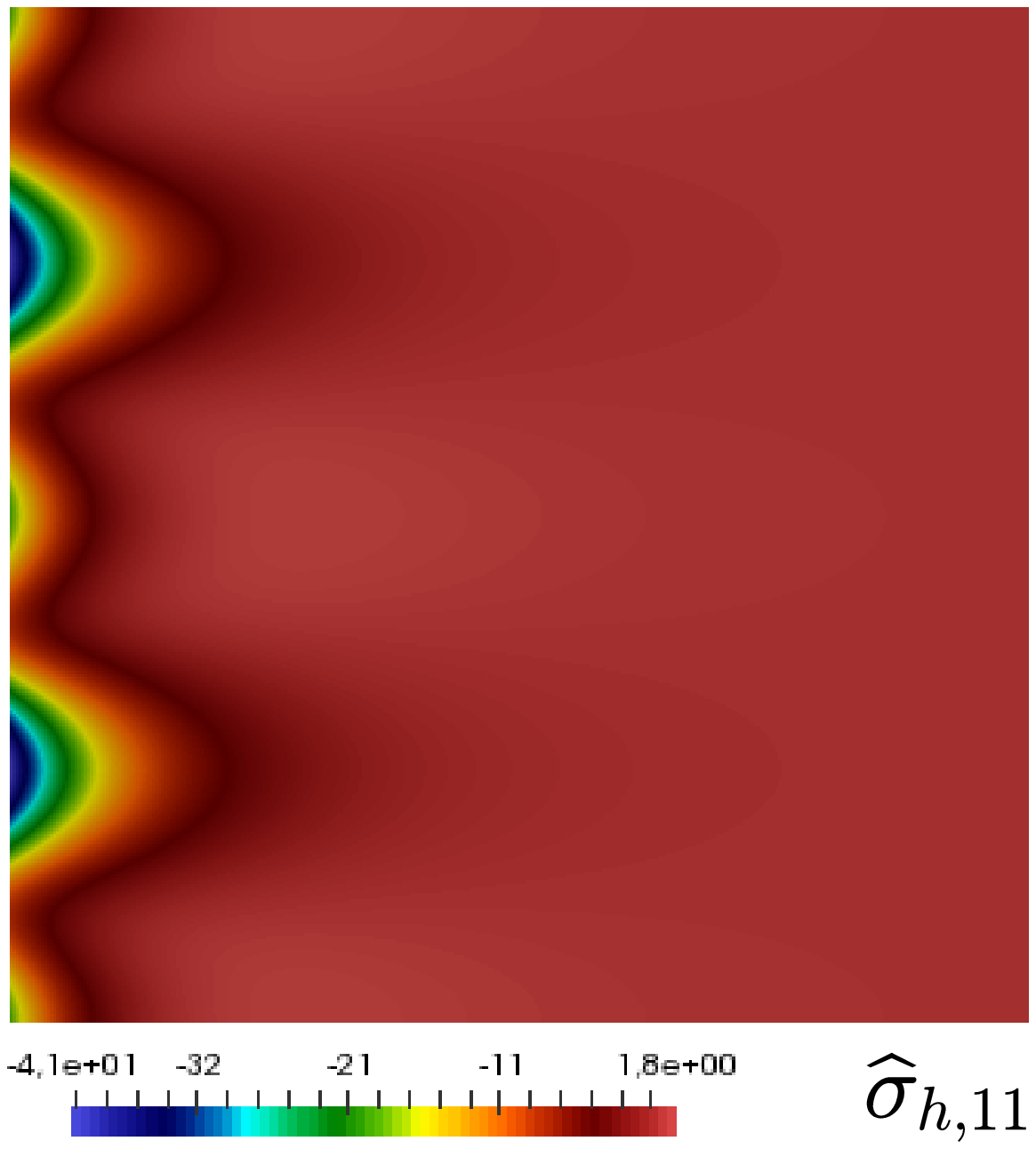}\qquad
                \includegraphics{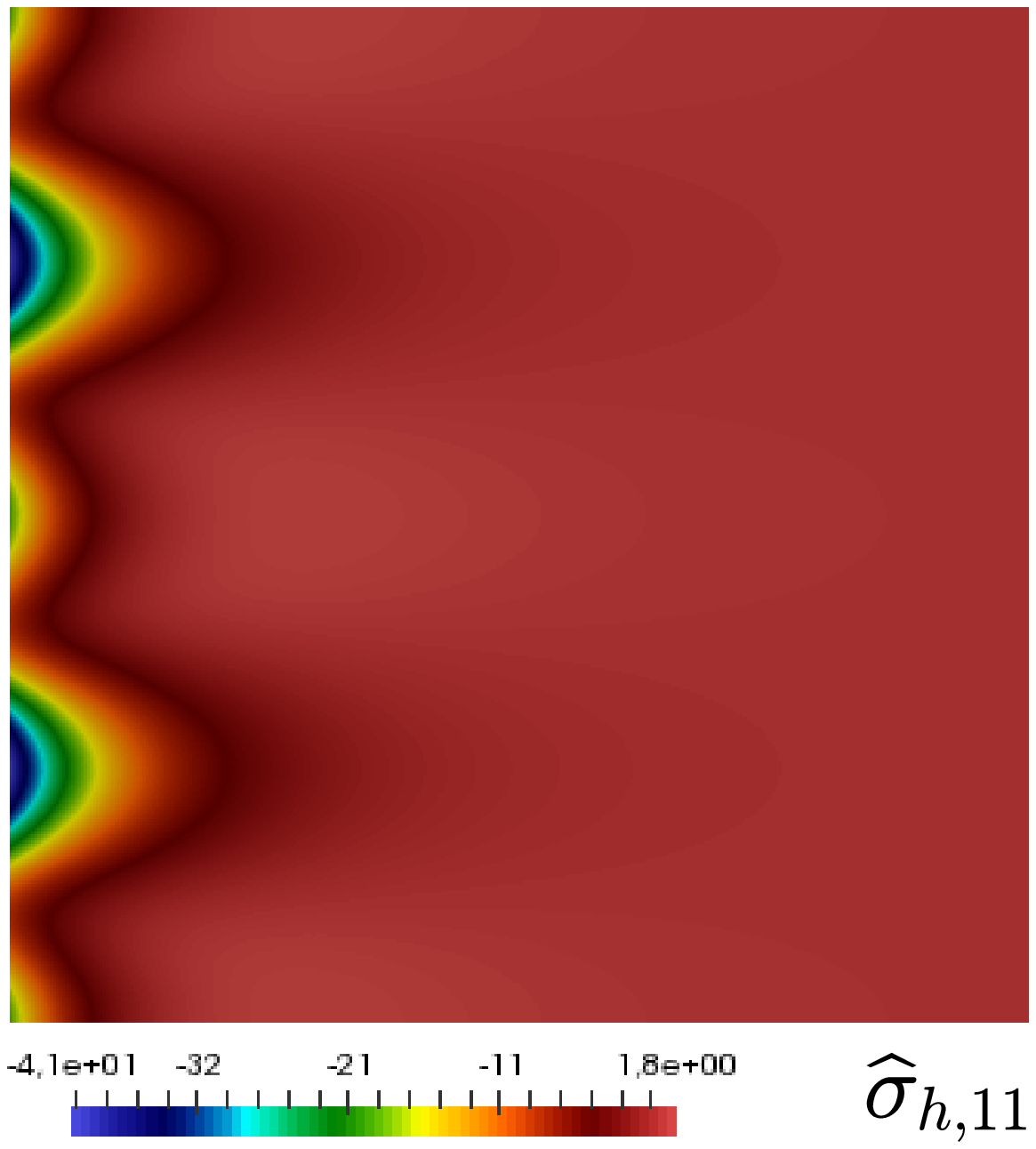}}\\[2ex]
\scalebox{0.45}{\includegraphics{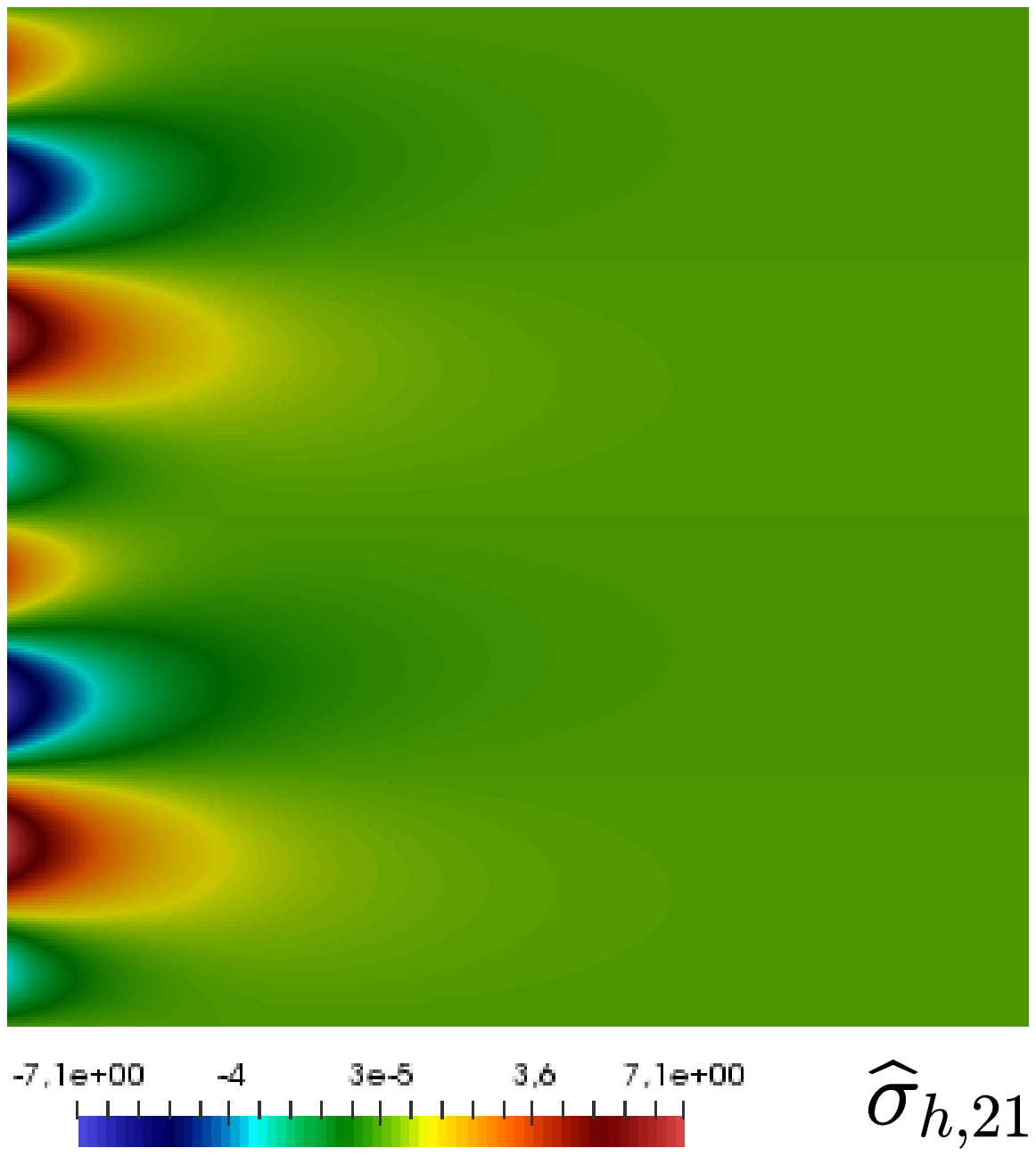}\qquad
                \includegraphics{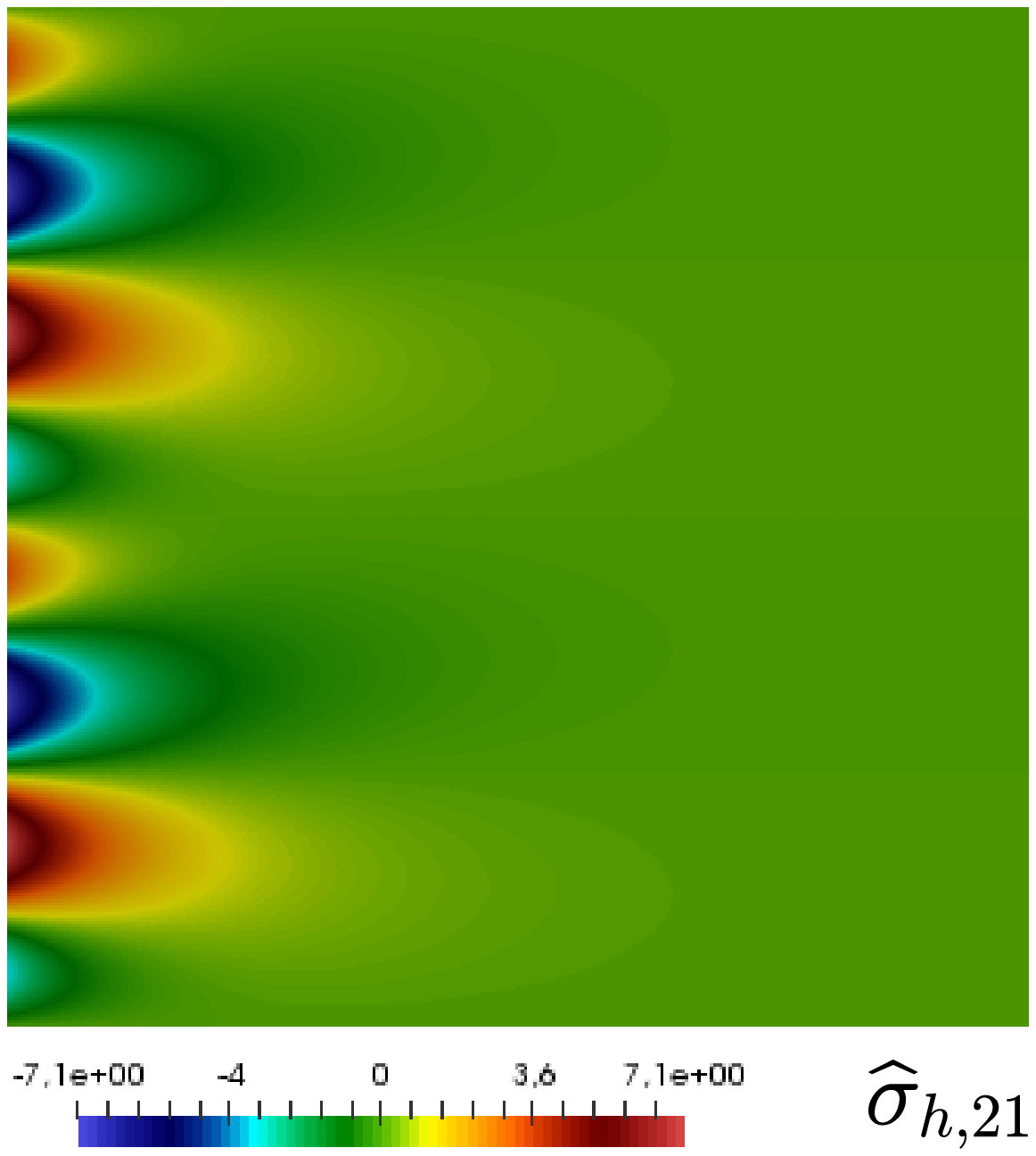}\qquad
                \includegraphics{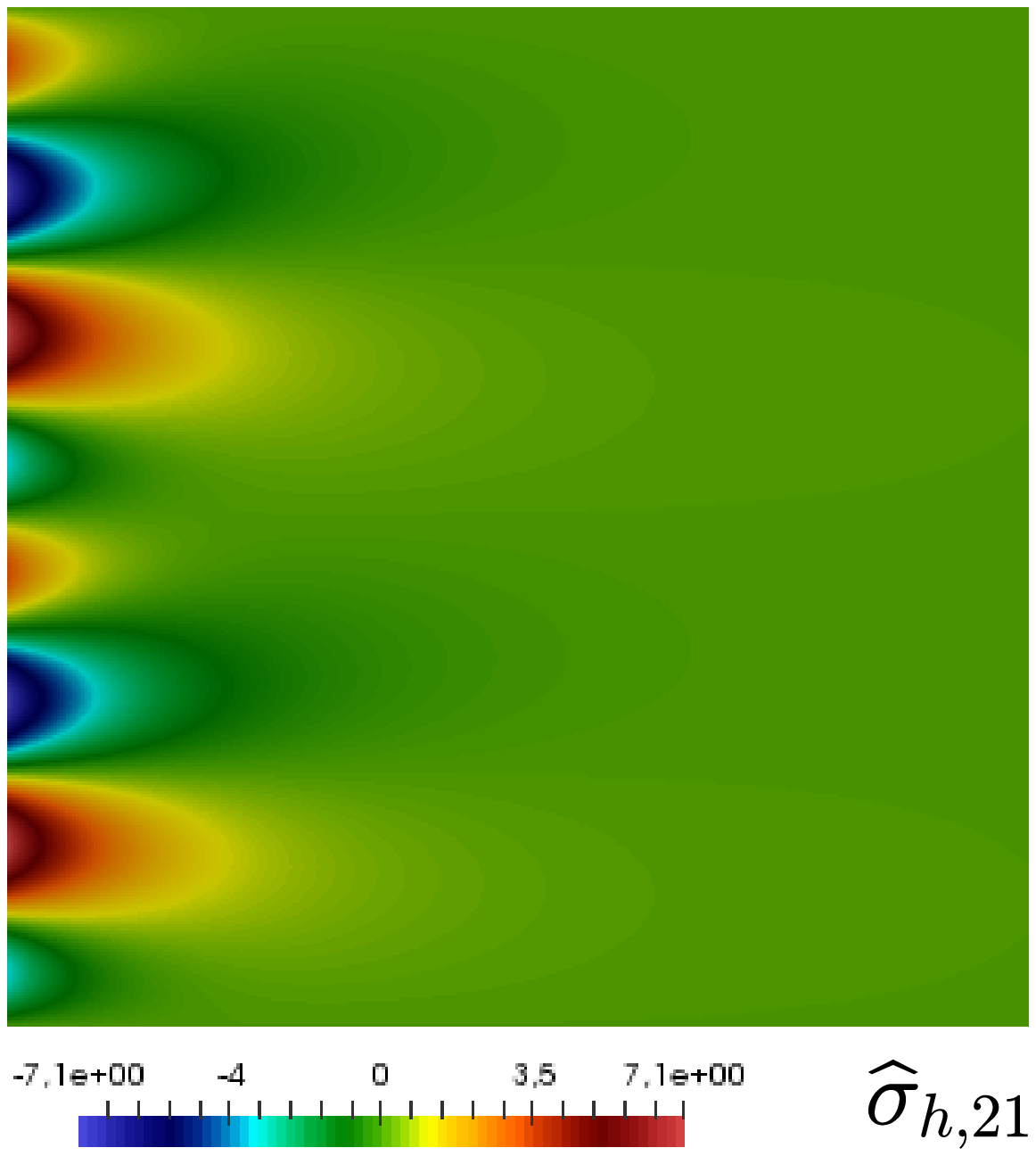}}\\[2ex]
\scalebox{0.45}{\includegraphics{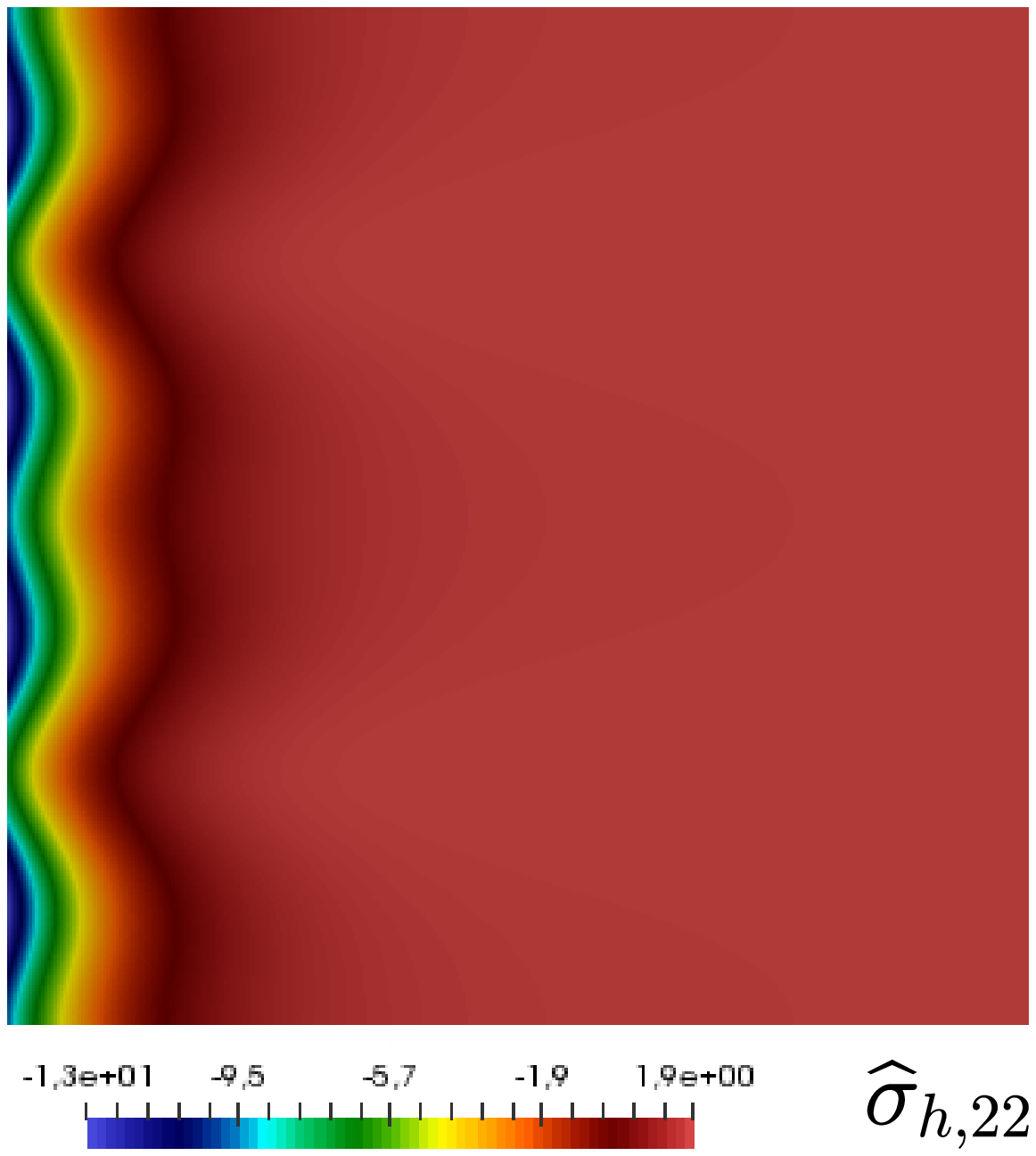}\qquad
                \includegraphics{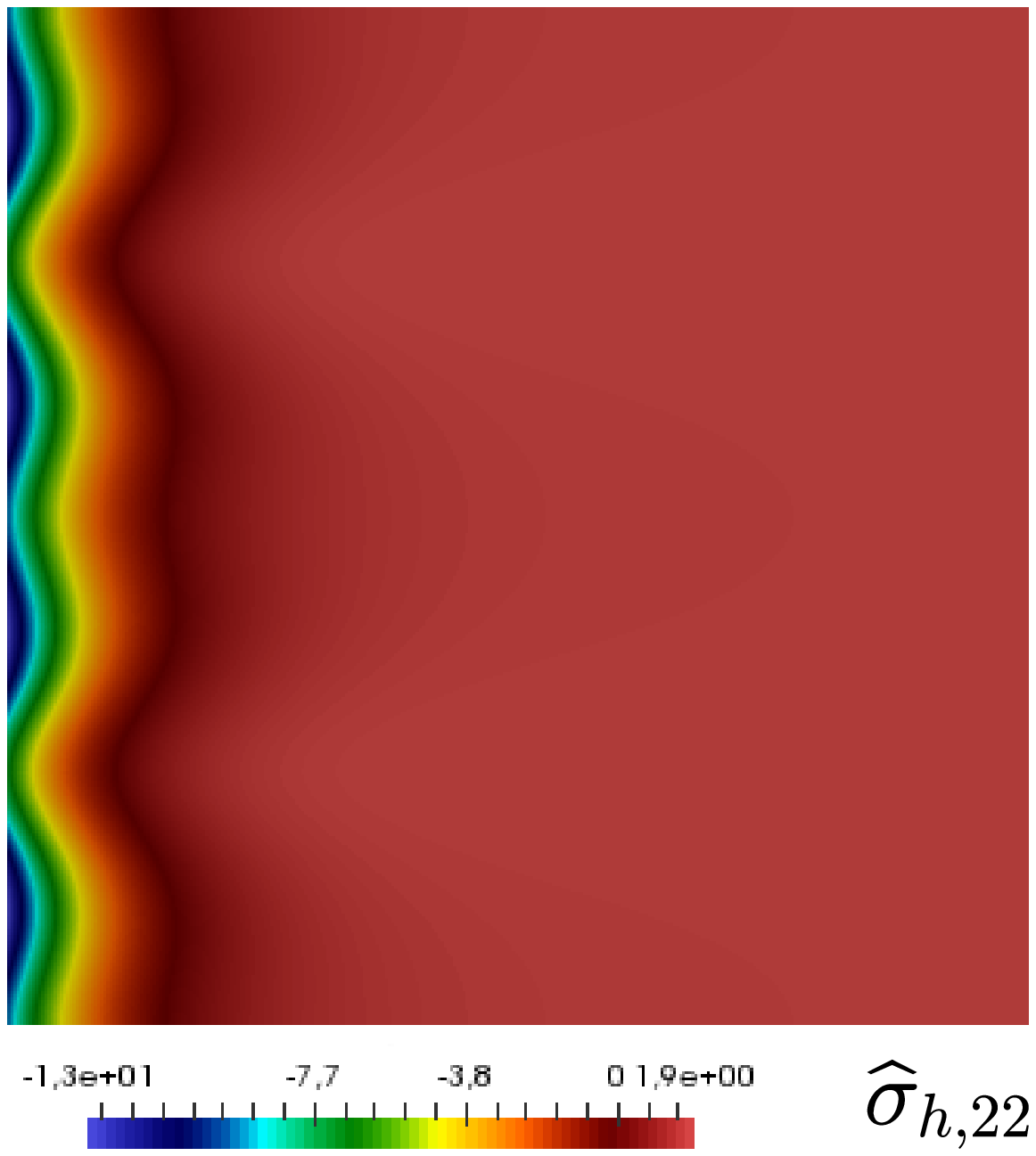}\qquad
                \includegraphics{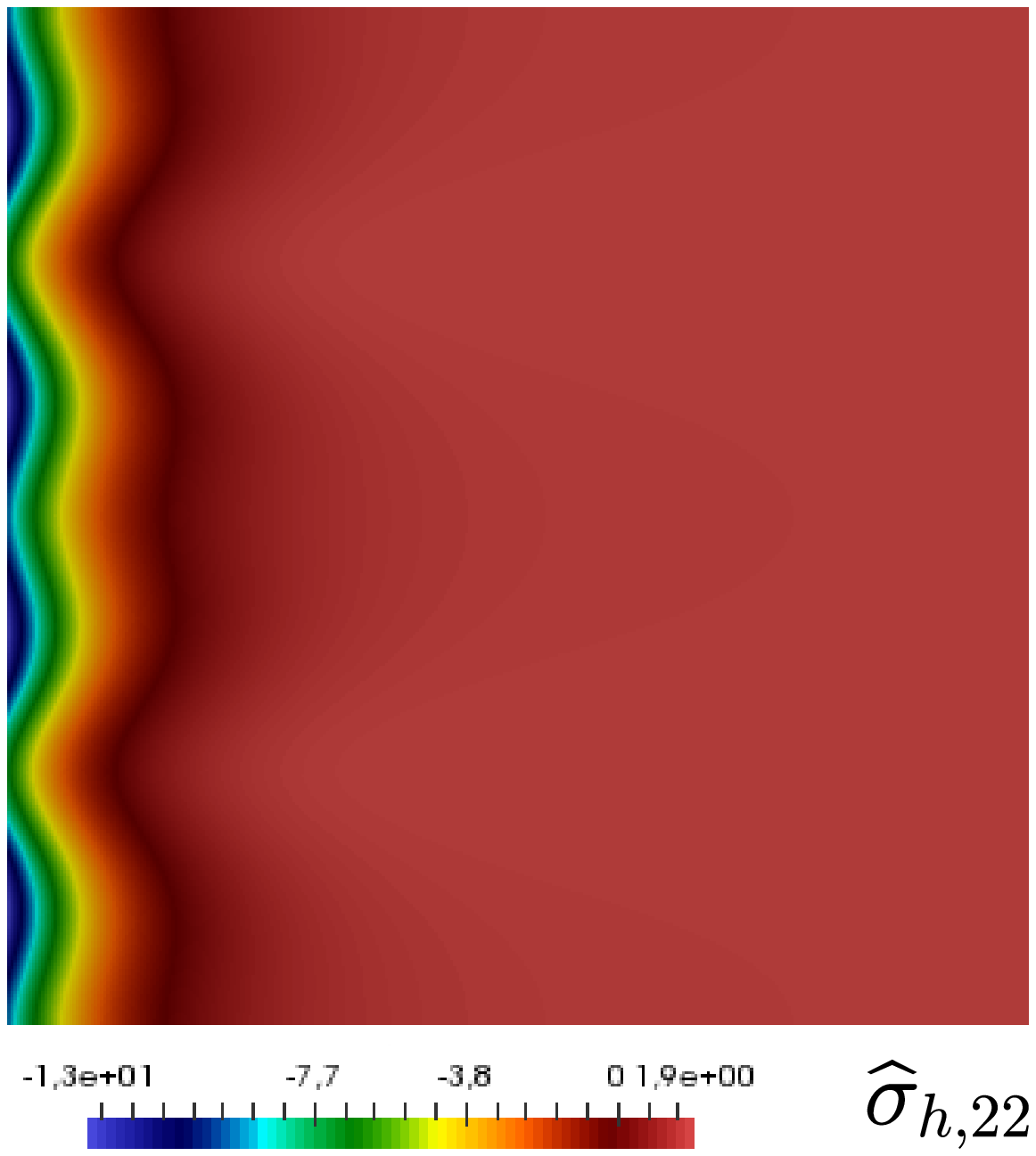}}
\caption{Example 2, $\widehat{\sigma}_{h,11}$ (top), $\widehat{\sigma}_{h,21}$ (center)
and $\widehat{\sigma}_{h,22}$ (bottom), using $k=2$ and the fifth mesh of triangles
(left column), quadrilaterals (center column) and hexagons (right column).}\label{fig:exa02-1}
\end{figure}

\begin{figure}[h!t]\centering
\scalebox{0.45}{\includegraphics{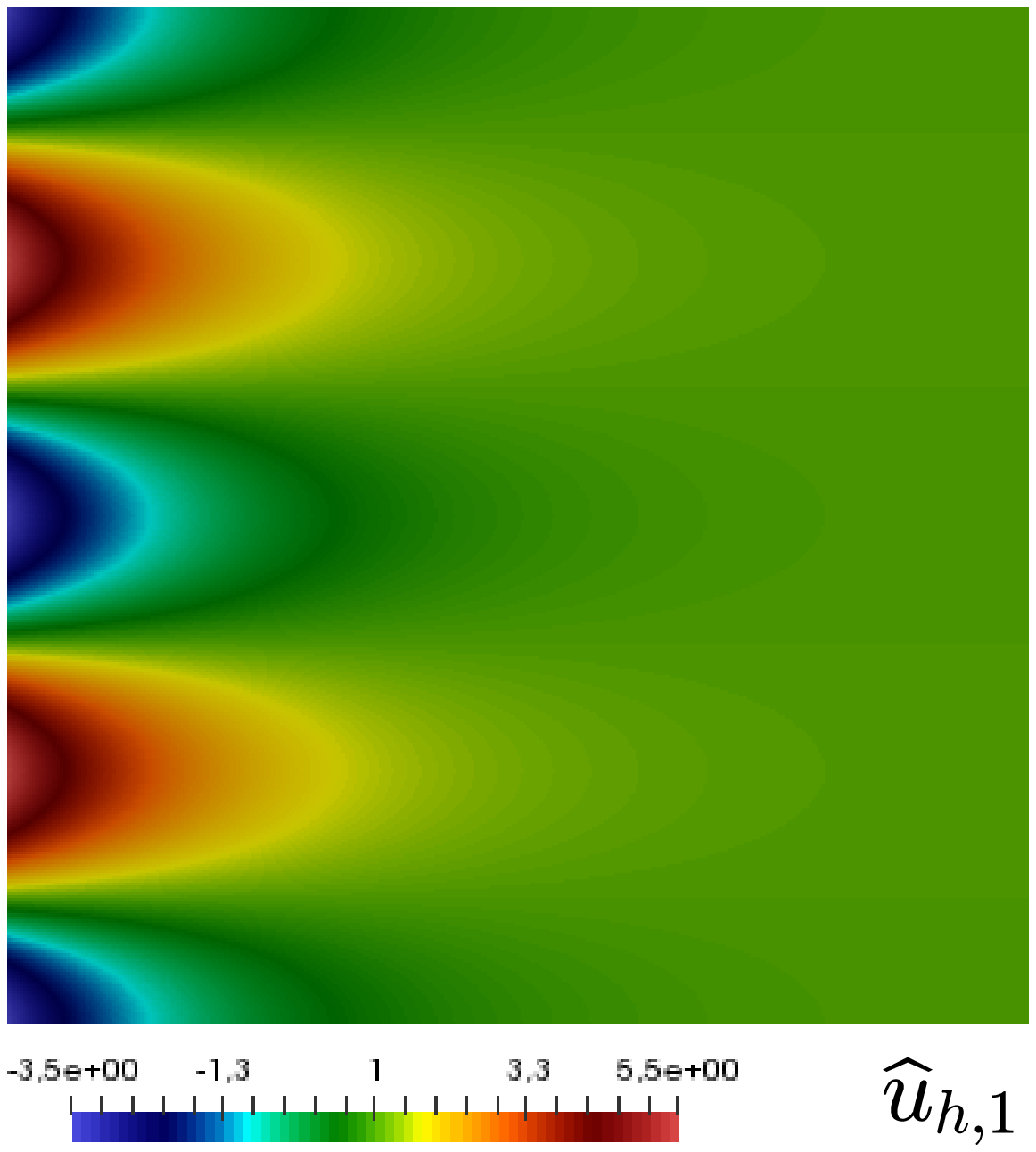}\qquad
                \includegraphics{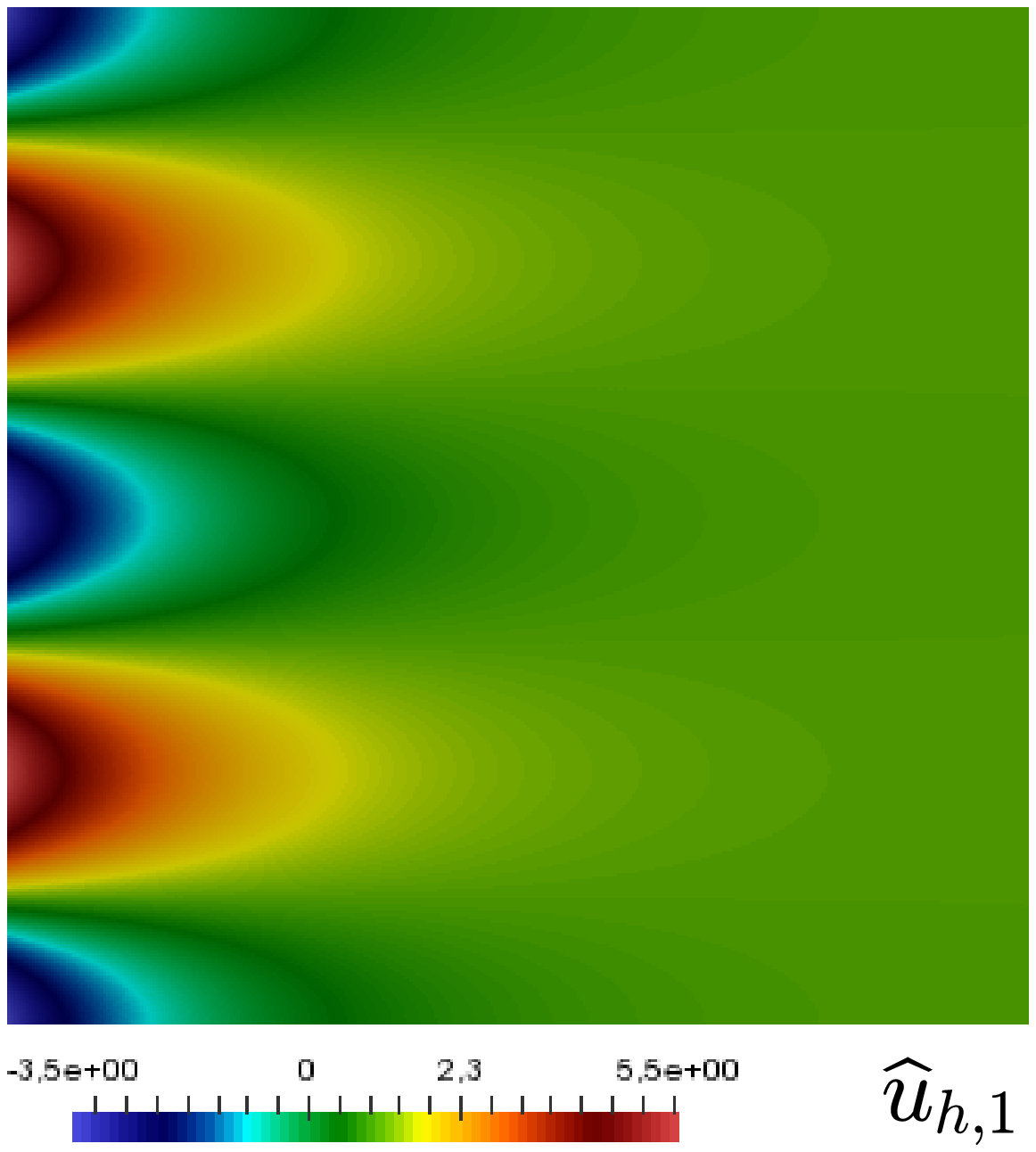}\qquad
                \includegraphics{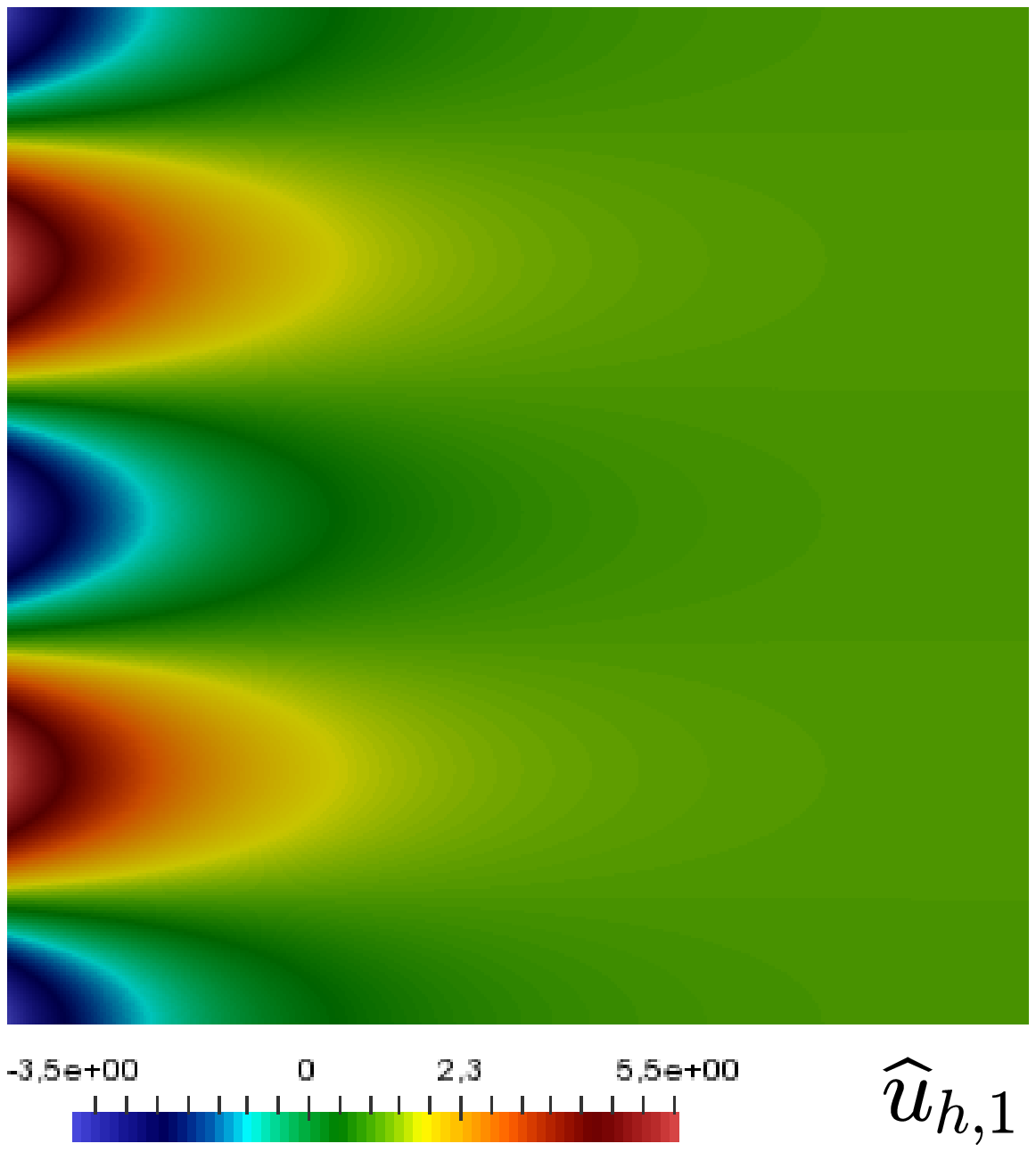}}\\[2ex]
\scalebox{0.45}{\includegraphics{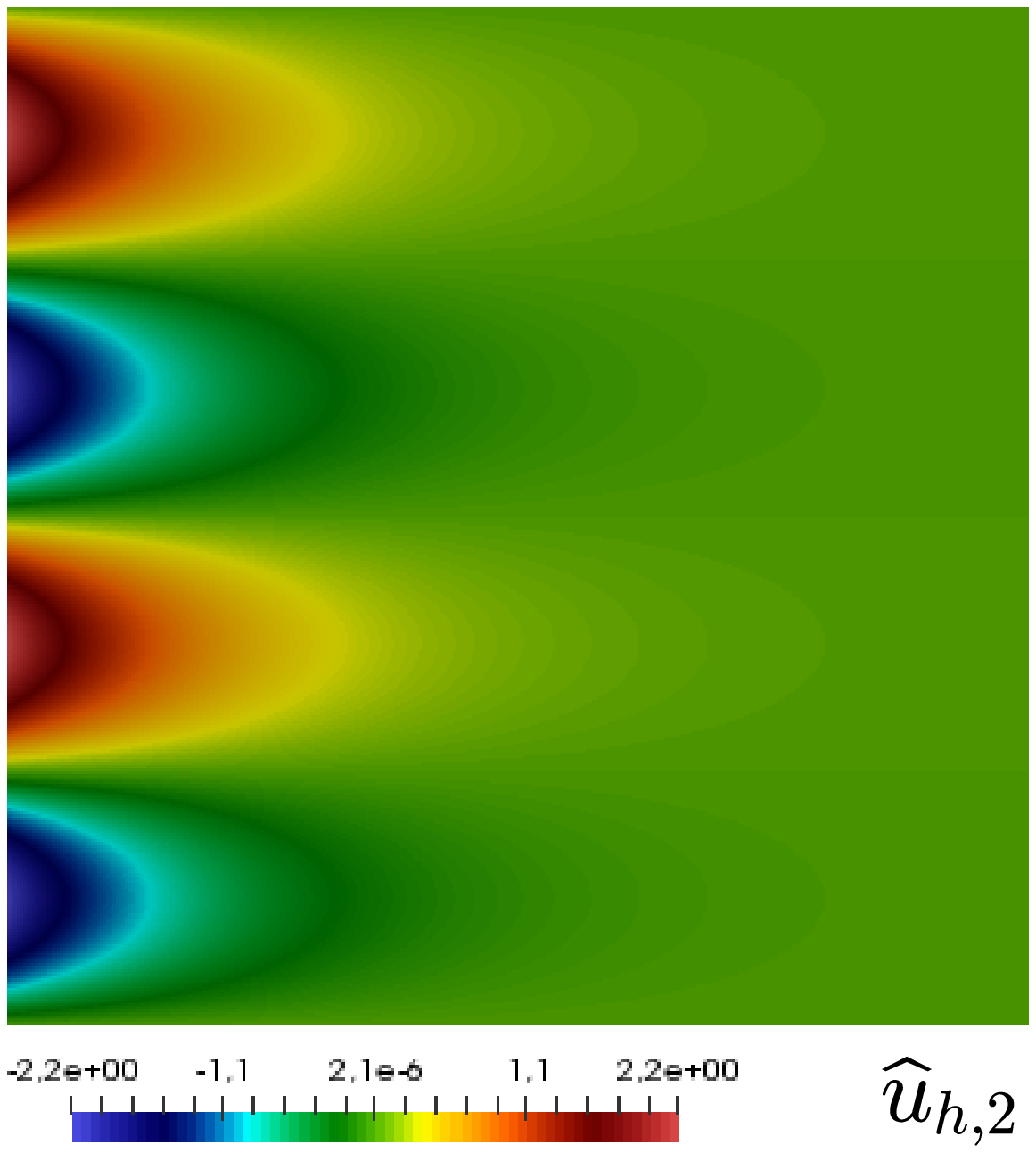}\qquad
                \includegraphics{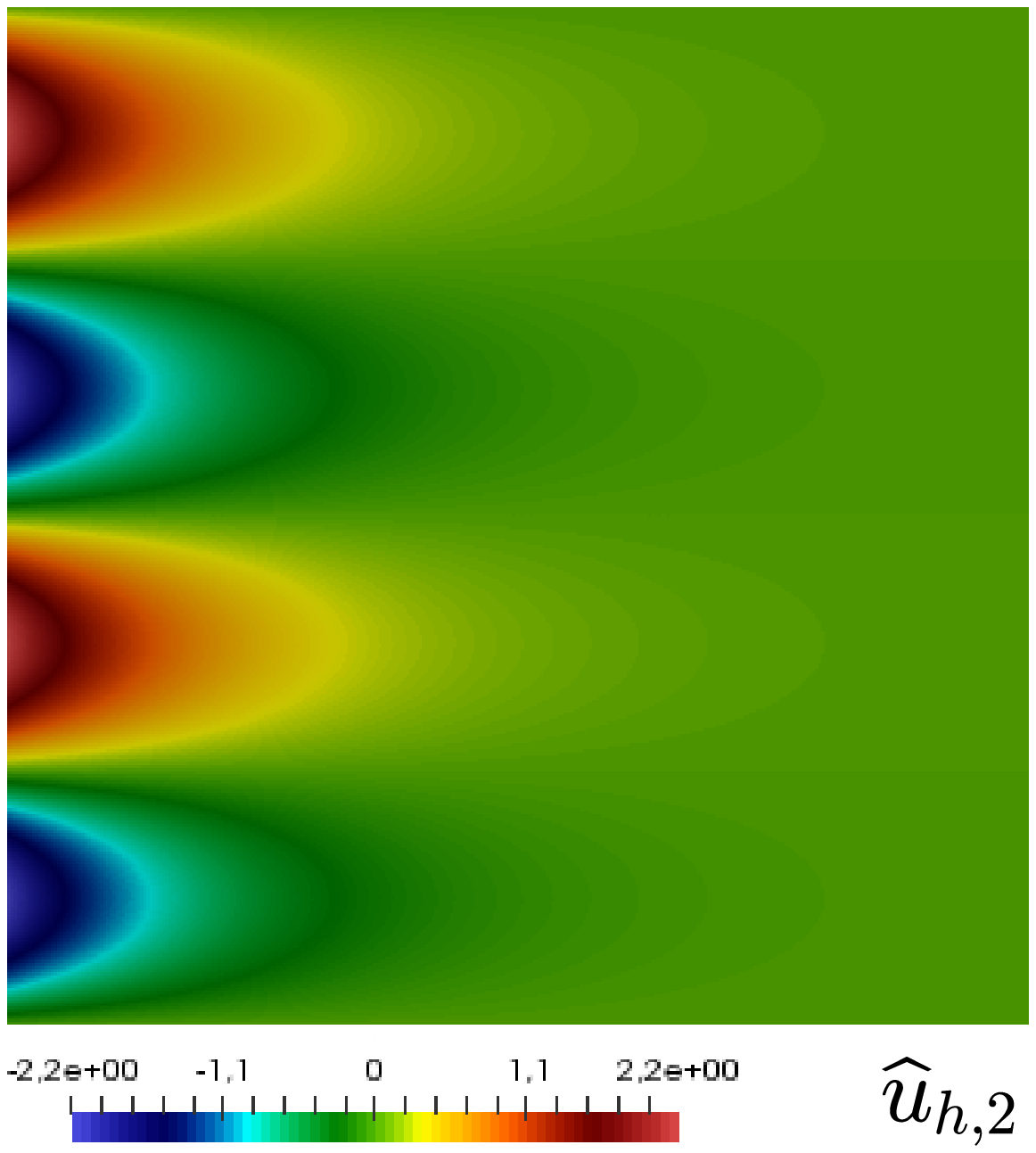}\qquad
                \includegraphics{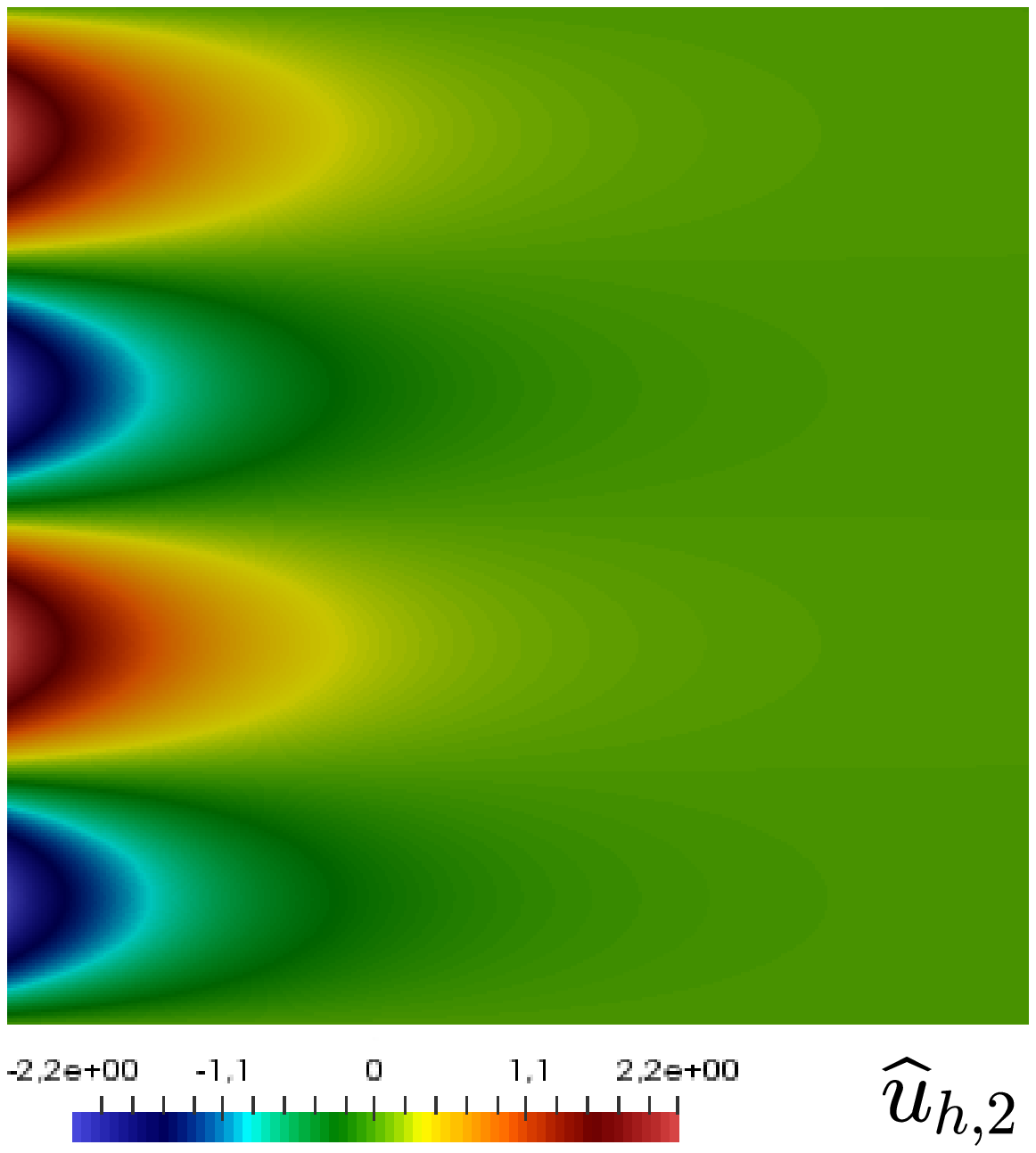}}\\[2ex]
\scalebox{0.45}{\includegraphics{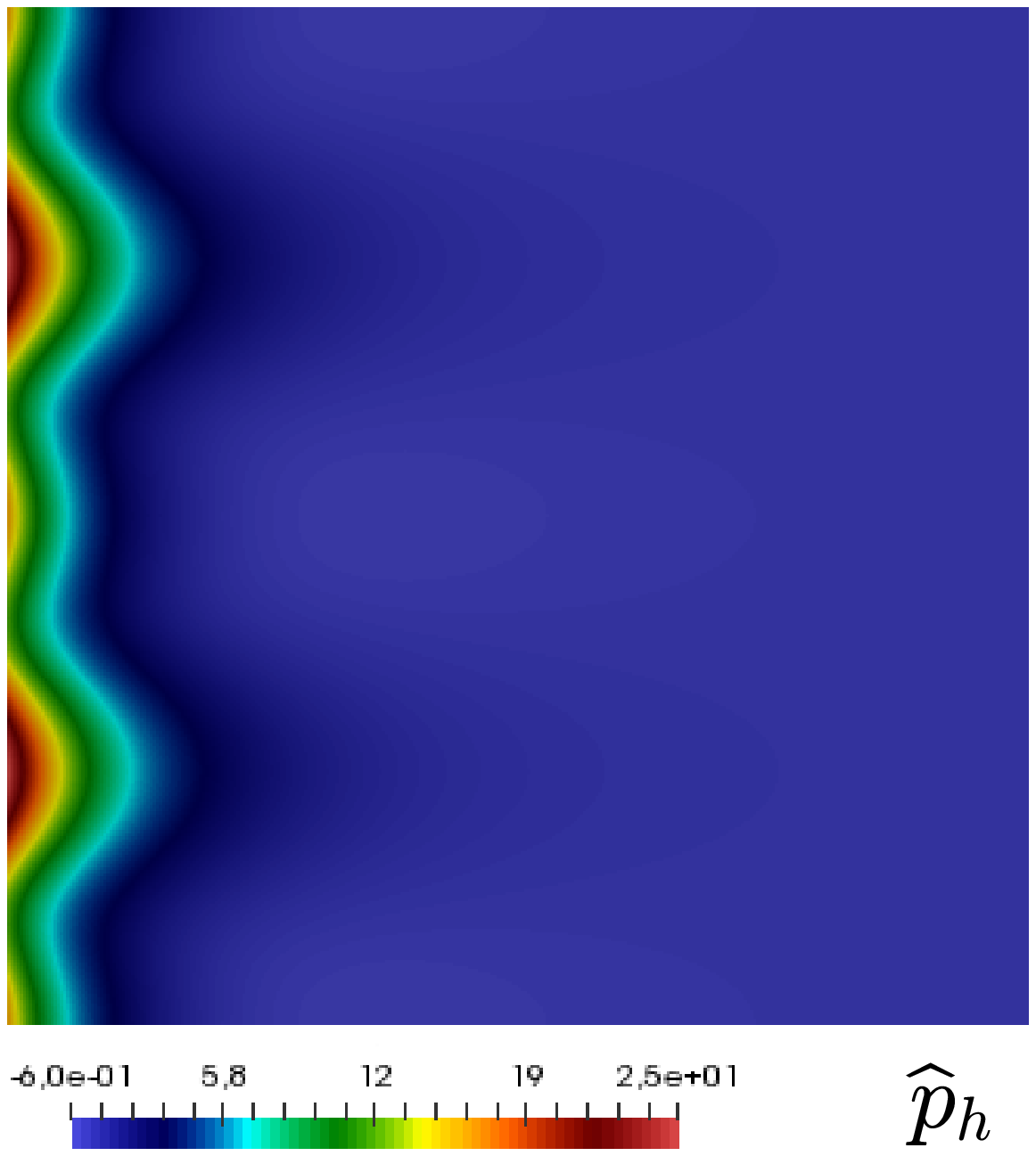}\qquad
                \includegraphics{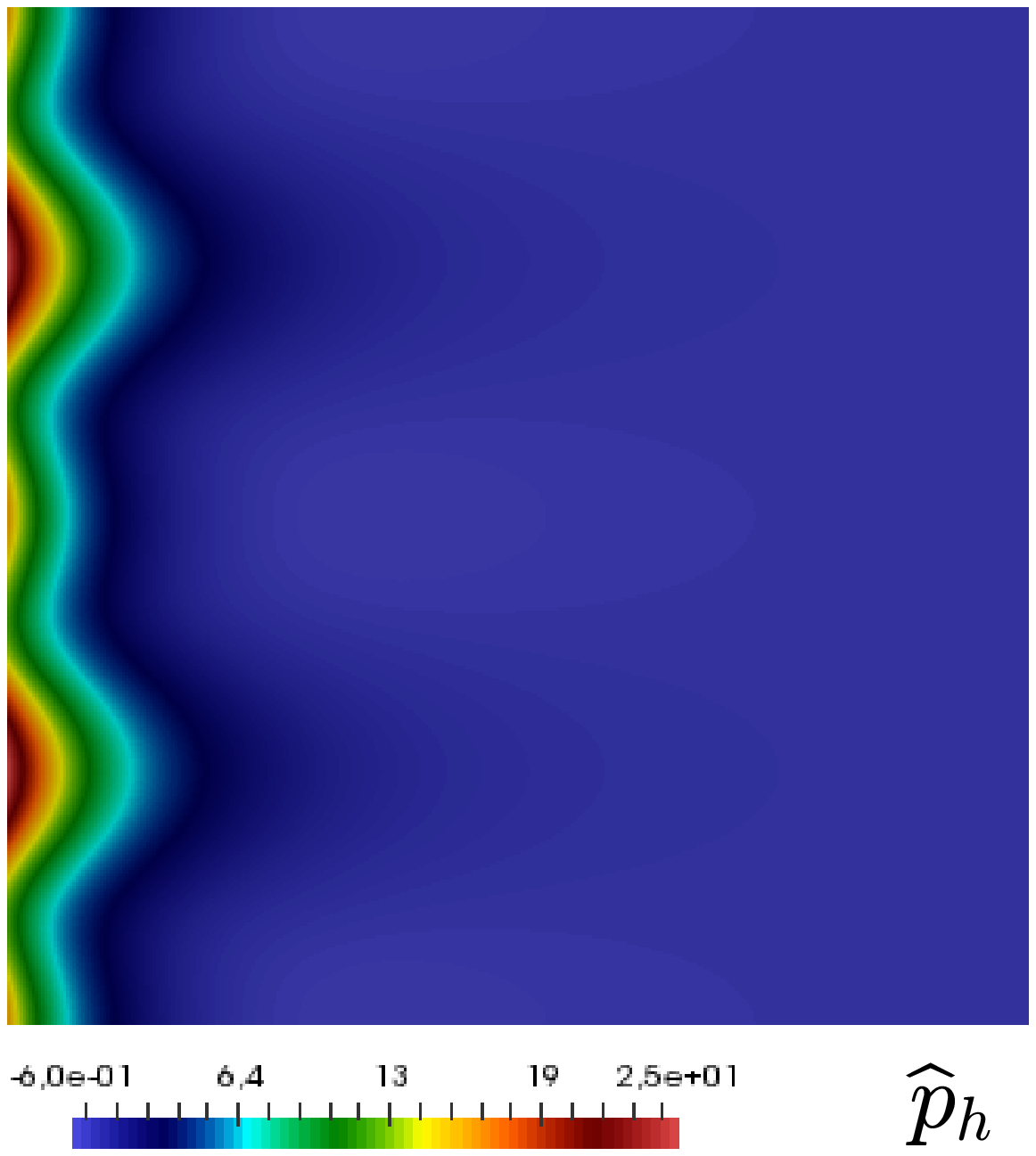}\qquad
                \includegraphics{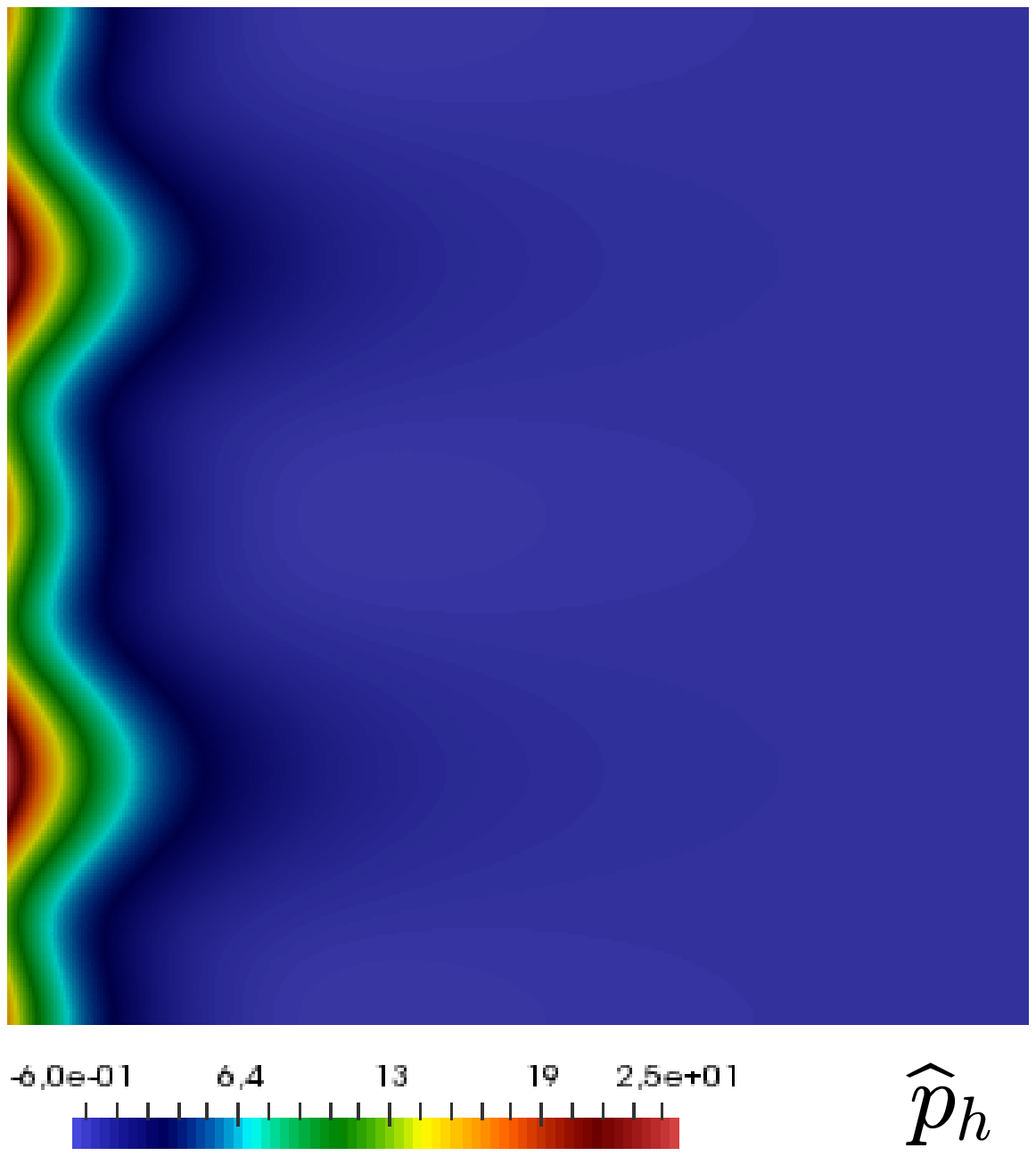}}
\caption{Example 2, $\widehat{u}_{h,1}$ (top), $\widehat{u}_{h,2}$ (center)
and $\widehat{p}_h$ (bottom), using $k=2$ and the fifth mesh of triangles
(left column), quadrilaterals (center column) and hexagons (right column).}\label{fig:exa02-2}
\end{figure}

In Tables \ref{tab:exa01-tri}, \ref{tab:exa01-qua} and \ref{tab:exa01-hex}, we summarize the
convergence history of the mixed virtual element scheme \eqref{eqn:discrete-scheme} as
applied to the present example. It follows from \cite[Theorem 5.6]{gms-vem-ns} that the rate
of convergence $O(h^{k+1})$ is attained by ${\tt e}(\bsig)$, ${\tt e}(\bu)$, ${\tt e}(p)$,
and ${\tt e}(\widetilde{\bsig})$, whereas the rate of convergence $O(h^k)$ is attained by 
${\tt e}(\widehat{\bu})$ in this smooth example, for triangular as well as for quadrilateral
and hexagonal meshes. Hence, the results validate the suitable behavior of our computational
implementation, along with the analysis carried out in \cite{gms-vem-ns}.
Finally, in Figures \ref{fig:exa02-1} and \ref{fig:exa02-2} we display some components of the
approximate solutions obtained in this section. They all correspond to those obtained with the
last mesh of each kind (triangles, quadrilaterals and hexagons, respectively) and for the
polynomial degree $k = 2$.

We end this paper by remarking some possible future directions.
It would be interesting to explain some computational aspects
about the adaptivity and the a posteriori error estimates for mixed virtual
element schemes (see, e.g., \cite{cgps-aposteriori-2017,ms-aposteriori-2020}).
In addition, in this work we used direct solvers for solving each
lineal system. However, we are interesting to develop an algebraic
or semi-algebraic multilevel preconditioner in order to employ iterative
solvers (see, e.g., \cite{castillo11,grss-uniciencia-2020}).

\begin{table}[h!t]\centering\footnotesize
\begin{tabular}{|c|c|c|c@{\hspace{1.5ex}}c|c@{\hspace{1.5ex}}c|c@{\hspace{1.5ex}}c|c@{\hspace{1.5ex}}c|c@{\hspace{1.5ex}}c|}\hline
$k$ & $h$ & $N$ & ${\tt e}(\bsig)$ & ${\tt r}(\bsig)$ & ${\tt e}(\bu)$ & ${\tt r}(\bu)$ & ${\tt e}(\widehat{\bu})$ & ${\tt r}(\widehat{\bu})$ & ${\tt e}(p)$ & ${\tt r}(p)$ & ${\tt e}(\widetilde{\bsig})$ & ${\tt r}(\widetilde{\bsig})$\\ \hline
  & 0.1230 &   4419 & 3.02e+00 & $--$ & 4.99e-01 & $--$ & 1.43e+01 & $--$ & 1.39e+00 & $--$ & 6.90e+00 & $--$\\
  & 0.0943 &   7443 & 2.24e+00 & 1.13 & 3.50e-01 & 1.34 & 1.43e+01 & 0.00 & 9.76e-01 & 1.33 & 5.28e+00 & 1.00\\
0 & 0.0488 &  27379 & 1.05e+00 & 1.15 & 1.44e-01 & 1.34 & 1.43e+01 & 0.00 & 3.98e-01 & 1.36 & 2.70e+00 & 1.02\\
  & 0.0354 &  51843 & 7.36e-01 & 1.10 & 9.63e-02 & 1.25 & 1.43e+01 & 0.00 & 2.65e-01 & 1.26 & 1.95e+00 & 1.01\\
  & 0.0283 &  80803 & 5.79e-01 & 1.07 & 7.40e-02 & 1.18 & 1.43e+01 & 0.00 & 2.03e-01 & 1.19 & 1.56e+00 & 1.01\\ \hline
  & 0.1230 &  19415 & 2.19e-01 & $--$ & 2.29e-02 & $--$ & 2.19e+00 & $--$ & 9.98e-02 & $--$ & 4.55e-01 & $--$\\
  & 0.0943 &  32883 & 1.29e-01 & 1.98 & 1.29e-02 & 2.15 & 1.68e+00 & 1.00 & 5.85e-02 & 2.01 & 2.69e-01 & 1.97\\
1 & 0.0488 & 122035 & 3.52e-02 & 1.98 & 3.35e-03 & 2.05 & 8.71e-01 & 1.00 & 1.56e-02 & 2.00 & 7.27e-02 & 1.99\\
  & 0.0354 & 231683 & 1.86e-02 & 1.98 & 1.76e-03 & 2.01 & 6.32e-01 & 1.00 & 8.23e-03 & 2.00 & 3.83e-02 & 1.99\\
  & 0.0283 & 361603 & 1.19e-02 & 1.98 & 1.12e-03 & 2.00 & 5.05e-01 & 1.00 & 5.27e-03 & 2.00 & 2.45e-02 & 1.99\\ \hline
  & 0.1230 &  40759 & 1.85e-02 & $--$ & 1.07e-03 & $--$ & 1.75e-01 & $--$ & 7.84e-03 & $--$ & 2.52e-02 & $--$\\
  & 0.0943 &  69123 & 8.39e-03 & 2.98 & 4.72e-04 & 3.07 & 1.03e-01 & 2.00 & 3.54e-03 & 2.99 & 1.14e-02 & 2.98\\
2 & 0.0488 & 257059 & 1.16e-03 & 3.00 & 6.45e-05 & 3.02 & 2.75e-02 & 2.00 & 4.86e-04 & 3.01 & 1.59e-03 & 2.99\\
  & 0.0354 & 488323 & 4.43e-04 & 3.00 & 2.46e-05 & 3.00 & 1.45e-02 & 2.00 & 1.84e-04 & 3.02 & 6.07e-04 & 3.00\\
  & 0.0283 & 762403 & 2.27e-04 & 3.00 & 1.26e-05 & 3.00 & 9.27e-03 & 2.00 & 9.38e-05 & 3.02 & 3.11e-04 & 3.00\\ \hline
\end{tabular}
\caption{History of convergence using triangles.}\label{tab:exa01-tri}
\end{table}

\begin{table}[h!t]\centering\footnotesize
\begin{tabular}{|c|c|c|c@{\hspace{1.5ex}}c|c@{\hspace{1.5ex}}c|c@{\hspace{1.5ex}}c|c@{\hspace{1.5ex}}c|c@{\hspace{1.5ex}}c|}\hline
$k$ & $h$ & $N$ & ${\tt e}(\bsig)$ & ${\tt r}(\bsig)$ & ${\tt e}(\bu)$ & ${\tt r}(\bu)$ & ${\tt e}(\widehat{\bu})$ & ${\tt r}(\widehat{\bu})$ & ${\tt e}(p)$ & ${\tt r}(p)$ & ${\tt e}(\widetilde{\bsig})$ & ${\tt r}(\widetilde{\bsig})$\\ \hline
  & 0.1008 & 6403 & 3.03e+00 & $--$ & 5.44e-01 & $--$ & 1.43e+01 & $--$ & 1.79e+00 & $--$ & 6.30e+00 & $--$\\
  & 0.0787 & 10417 & 2.20e+00 & 1.30 & 3.85e-01 & 1.39 & 1.43e+01 & 0.00 & 1.25e+00 & 1.45 & 4.88e+00 & 1.03\\
0 & 0.0404 & 39043 & 8.69e-01 & 1.39 & 1.49e-01 & 1.42 & 1.43e+01 & 0.00 & 4.25e-01 & 1.61 & 2.40e+00 & 1.06\\
  & 0.0307 & 66993 & 6.04e-01 & 1.34 & 1.02e-01 & 1.40 & 1.43e+01 & 0.00 & 2.73e-01 & 1.63 & 1.81e+00 & 1.04\\
  & 0.0229 & 120417 & 4.16e-01 & 1.27 & 6.89e-02 & 1.32 & 1.43e+01 & 0.00 & 1.71e-01 & 1.58 & 1.33e+00 & 1.03\\ \hline
  & 0.1008 & 23043 & 1.57e-01 & $--$ & 1.81e-02 & $--$ & 1.82e+00 & $--$ & 5.86e-02 & $--$ & 3.16e-01 & $--$\\
  & 0.0787 & 37641 & 9.59e-02 & 2.00 & 1.10e-02 & 2.00 & 1.42e+00 & 0.99 & 3.50e-02 & 2.08 & 1.95e-01 & 1.95\\
1 & 0.0404 & 142083 & 2.49e-02 & 2.02 & 2.88e-03 & 2.01 & 7.29e-01 & 1.00 & 8.78e-03 & 2.07 & 5.12e-02 & 2.00\\
  & 0.0307 & 244233 & 1.44e-02 & 2.02 & 1.67e-03 & 2.00 & 5.56e-01 & 1.00 & 5.04e-03 & 2.04 & 2.97e-02 & 2.01\\
  & 0.0229 & 439641 & 7.95e-03 & 2.01 & 9.28e-04 & 2.00 & 4.14e-01 & 1.00 & 2.77e-03 & 2.03 & 1.65e-02 & 2.00\\ \hline
  & 0.1008 & 45827 & 1.34e-02 & $--$ & 8.58e-04 & $--$ & 1.30e-01 & $--$ & 4.55e-03 & $--$ & 1.60e-02 & $--$\\
  & 0.0787 & 74951 & 6.41e-03 & 2.98 & 4.08e-04 & 3.00 & 7.96e-02 & 1.99 & 2.17e-03 & 3.00 & 7.69e-03 & 2.96\\
2 & 0.0404 & 283523 & 8.55e-04 & 3.02 & 5.45e-05 & 3.01 & 2.09e-02 & 2.00 & 2.81e-04 & 3.05 & 1.03e-03 & 3.01\\
  & 0.0307 & 487623 & 3.77e-04 & 3.01 & 2.41e-05 & 3.00 & 1.21e-02 & 2.00 & 1.23e-04 & 3.04 & 4.54e-04 & 3.01\\
  & 0.0229 & 878151 & 1.55e-04 & 3.00 & 9.96e-06 & 3.00 & 6.74e-03 & 2.00 & 5.08e-05 & 3.00 & 1.88e-04 & 3.00\\ \hline
\end{tabular}
\caption{History of convergence using quadrilaterals.}\label{tab:exa01-qua}
\end{table}

\begin{table}[h!t]\centering\footnotesize
\begin{tabular}{|c|c|c|c@{\hspace{1.5ex}}c|c@{\hspace{1.5ex}}c|c@{\hspace{1.5ex}}c|c@{\hspace{1.5ex}}c|c@{\hspace{1.5ex}}c|}\hline
$k$ & $h$ & $N$ & ${\tt e}(\bsig)$ & ${\tt r}(\bsig)$ & ${\tt e}(\bu)$ & ${\tt r}(\bu)$ & ${\tt e}(\widehat{\bu})$ & ${\tt r}(\widehat{\bu})$ & ${\tt e}(p)$ & ${\tt r}(p)$ & ${\tt e}(\widetilde{\bsig})$ & ${\tt r}(\widetilde{\bsig})$\\ \hline
  & 0.0959 & 10535 & 2.42e+00 & $--$ & 4.43e-01 & $--$ & 1.43e+01 & $--$ & 1.27e+00 & $--$ & 5.24e+00 & $--$\\
  & 0.0732 & 17897 & 1.83e+00 & 1.05 & 3.22e-01 & 1.18 & 1.43e+01 & 0.00 & 8.92e-01 & 1.31 & 4.07e+00 & 0.94\\
0 & 0.0527 & 34143 & 1.16e+00 & 1.39 & 1.98e-01 & 1.48 & 1.43e+01 & 0.00 & 5.15e-01 & 1.67 & 2.92e+00 & 1.01\\
  & 0.0390 & 61887 & 7.77e-01 & 1.34 & 1.32e-01 & 1.36 & 1.43e+01 & 0.00 & 3.17e-01 & 1.62 & 2.17e+00 & 0.99\\
  & 0.0301 & 103495 & 5.66e-01 & 1.22 & 9.56e-02 & 1.23 & 1.43e+01 & 0.00 & 2.12e-01 & 1.54 & 1.67e+00 & 1.01\\ \hline
  & 0.0959 & 31707 & 1.55e-01 & $--$ & 1.90e-02 & $--$ & 1.87e+00 & $--$ & 5.23e-02 & $--$ & 2.84e-01 & $--$\\
  & 0.0732 & 53681 & 9.10e-02 & 1.97 & 1.12e-02 & 1.97 & 1.44e+00 & 0.97 & 3.02e-02 & 2.03 & 1.68e-01 & 1.94\\
1 & 0.0527 & 102623 & 4.76e-02 & 1.98 & 5.84e-03 & 1.97 & 1.04e+00 & 0.99 & 1.55e-02 & 2.04 & 9.00e-02 & 1.90\\
  & 0.0390 & 185651 & 2.64e-02 & 1.97 & 3.23e-03 & 1.98 & 7.72e-01 & 0.99 & 8.46e-03 & 2.01 & 5.02e-02 & 1.94\\
  & 0.0301 & 310827 & 1.56e-02 & 2.01 & 1.93e-03 & 1.99 & 5.97e-01 & 0.99 & 4.99e-03 & 2.03 & 2.98e-02 & 2.01\\ \hline
  & 0.0959 & 59263 & 1.55e-02 & $--$ & 9.46e-04 & $--$ & 1.39e-01 & $--$ & 5.96e-03 & $--$ & 1.78e-02 & $--$\\
  & 0.0732 & 100199 & 7.20e-03 & 2.85 & 4.30e-04 & 2.93 & 8.18e-02 & 1.95 & 2.79e-03 & 2.82 & 8.22e-03 & 2.88\\
2 & 0.0527 & 191707 & 2.75e-03 & 2.93 & 1.61e-04 & 2.98 & 4.25e-02 & 1.99 & 1.04e-03 & 2.99 & 3.17e-03 & 2.90\\
  & 0.0390 & 346543 & 1.14e-03 & 2.94 & 6.57e-05 & 2.99 & 2.34e-02 & 1.99 & 4.36e-04 & 2.91 & 1.31e-03 & 2.94\\
  & 0.0301 & 580463 & 5.26e-04 & 2.98 & 3.03e-05 & 2.99 & 1.39e-02 & 1.99 & 2.01e-04 & 2.97 & 6.02e-04 & 2.99\\ \hline
\end{tabular}
\caption{History of convergence using hexagons.}\label{tab:exa01-hex}
\end{table}

\section*{Acknowledgements}

The authors would like to thank Gabriel N. Gatica, CI$^2$MA and Departamento de
Ingenier\'ia Matem\'a\-tica, Universidad de Concepci\' on, Chile, for his suggestions
that significantly influenced the organization of this paper. On the other hand,
the work of Fil\'ander A. Sequeira was partially supported by Universidad Nacional,
Costa Rica, through the project 0103-18.


\end{document}